\documentclass[3p,times]{elsarticle}

\usepackage{graphicx}
\usepackage{psfrag}
\usepackage{color}

\usepackage{amssymb}
\usepackage{amsmath}
\usepackage{dsfont}
\usepackage{rotating}

\newcommand{\U}{\mathbf{U}}

\newcommand{\Q}{\mathbf{Q}} 
\newcommand{\q}{\mathbf{q}} 
\newcommand{\f}{\mathbf{f}} 
\newcommand{\g}{\mathbf{g}}
\newcommand{\h}{\mathbf{h}}

\newcommand{\F}{\mathbf{F}} 
 
\newcommand{\x}{\mathbf{x}} 
\newcommand{\xxi}{\boldsymbol{\xi}}
\newcommand{\kk}{\boldsymbol{\kappa}}

\renewcommand{\v}{\mathbf{v}}

\newcommand{\quotew}[1]{``#1''}

\begin{document}

\begin{frontmatter}

\journal{Journal of Computational Physics} 



\title{Arbitrary-Lagrangian-Eulerian discontinuous Galerkin schemes 
with a posteriori subcell finite volume limiting on moving unstructured meshes} 


\author[UNITN]{Walter Boscheri}
\ead{walter.boscheri@unitn.it}
\author[UNITN]{Michael Dumbser}
\ead{michael.dumbser@unitn.it}
\address[UNITN]{Laboratory of Applied Mathematics \\ Department of Civil, Environmental and Mechanical Engineering \\ University of Trento, Via Mesiano 77, I-38123 Trento, Italy}

\begin{abstract}
We present a new family of high order accurate \textit{fully discrete} one-step Discontinuous Galerkin (DG) finite element schemes on moving unstructured meshes for the solution of nonlinear hyperbolic PDE in multiple space dimensions, which may also include parabolic terms in order to model \textit{dissipative} transport processes, like molecular viscosity or heat conduction. 
High order piecewise polynomials of degree $N$ are adopted to represent the discrete solution at each time level and within each spatial control volume of the computational grid, while high order of accuracy in time is achieved by the ADER approach, making use of an element-local space-time Galerkin finite element predictor. 
A novel nodal solver algorithm based on the HLL flux is derived to compute the velocity for each nodal degree of freedom that describes the current mesh geometry. 
In our algorithm the spatial mesh configuration can be defined in two different ways: either by an \textit{isoparametric} approach  that generates \textit{curved} control volumes, or by a \textit{piecewise linear} decomposition of each spatial control volume into simplex sub-elements. Each technique generates a corresponding number of geometrical degrees of freedom needed to describe the current mesh configuration and which must be considered by the nodal solver for determining the grid velocity.  

The connection of the old mesh configuration at time $t^{n}$ with the new one at time $t^{n+1}$ provides the space-time control volumes on which the governing equations have to be integrated in order to obtain the time evolution of the discrete solution. 
Our numerical method belongs to the category of so-called \textit{direct} Arbitrary-Lagrangian-Eulerian (ALE) schemes, where a space-time conservation formulation of the governing PDE system is considered and which already takes into account the new grid geometry 
(including a possible rezoning step) \textit{directly} during the computation of the numerical fluxes. 
We emphasize that our method is a \textit{moving mesh} method, as opposed to \textit{total Lagrangian} formulations that are based on a \textit{fixed} computational grid and which instead evolve the \textit{mapping} of the reference configuration to the current one. 

Our new Lagrangian-type DG scheme adopts the novel \textit{a posteriori} sub-cell finite volume limiter method recently developed in 
\cite{DGLimiter3} for fixed unstructured grids. In this approach, the validity of the \textit{candidate solution} produced in each cell by an \textit{unlimited} ADER-DG scheme is verified against a set of physical and numerical detection criteria, such as the positivity of pressure and density, the absence of floating point errors (NaN) and the satisfaction of a relaxed discrete maximum principle (DMP) in the sense of polynomials. Those cells which do not satisfy all of the above criteria are flagged as 
\textit{troubled cells} and are \textit{recomputed} at the aid of a more robust second order TVD finite volume scheme. To preserve the subcell resolution capability of the original DG scheme, the FV limiter is run on a sub-grid that is $2N+1$ times finer compared to the mesh of the original unlimited DG scheme. The new subcell averages are then gathered back into a high order DG polynomial by a usual conservative finite volume reconstruction operator.  

The numerical convergence rates of the new ALE ADER-DG schemes are studied up to fourth order in space and time and several test problems are simulated in order to check the accuracy and the robustness of the proposed numerical method in the context of the 
Euler and Navier-Stokes equations for compressible gas dynamics, considering both inviscid and viscous fluids. Finally, an 
application inspired by Inertial Confinement Fusion (ICF) type flows is considered by solving the Euler equations and the PDE of viscous and resistive magnetohydrodynamics (VRMHD). 

\end{abstract}
\end{frontmatter}


\section{Introduction}
\label{sec.introduction}

Lagrangian algorithms have become very popular in the last decades \cite{vonneumann50,Benson1992,munz94,Caramana1998,Smith1999,Maire2007,Despres2009} due to to the excellent properties achieved by these numerical methods in the resolution of moving material interfaces and contact waves. Since the computational mesh is moving with the local fluid velocity, Lagrangian methods are typically affected by much less numerical dissipation compared to classical Eulerian approaches on fixed grids, hence obtaining a more accurate approximation of the solution.

As governing equations we consider nonlinear systems of hyperbolic conservation laws combined with parabolic terms, which cover a wide range of phenomena, such as environmental and meteorological flows, hydrodynamic and thermodynamic problems, plasma flows as well as the dynamics of many industrial and mechanical processes. A widespread technique for the solution of nonlinear hyperbolic systems of PDE is given by Godunov-type finite volume schemes \cite{GodunovRS,vanLeerRS}. In this approach the numerical solution is stored under the form of piecewise constant cell averages within each control volume of the computational mesh, and the time evolution is obtained by considering the integral form of the conservation laws. A lot of work has been done in the development of Lagrangian finite volume schemes \cite{munz94,Depres2012,ShashkovCellCentered,Maire2009,Maire2009b,Maire2010,Maire2011} achieving up to second order of accuracy in space and time. Higher order Lagrangian-type schemes based on ENO reconstruction have been introduced for the first time by Cheng and Shu in \cite{chengshu1,chengshu2}. Since all variables are located at the cell barycenter, these methods are also referred to as \textit{cell-centered} Lagrangian algorithms, contrarily to the \textit{staggered mesh} schemes \cite{StagLag,LoubereSedov3D,maire_loubere_vachal10}, where the velocity is defined at the grid vertices and the other variables are considered at the cell center.

Lagrangian methods either directly move the mesh, or they evolve the mapping of a reference configuration onto the current one. 
In any case, they may produce highly distorted elements in the current configuration, depending on the flow motion. Problems arise in particular for strong shear flows. This can lead to highly deformed and distorted cells, which inevitably lead to very small time steps in the case of explicit schemes due to the CFL stability condition. Even invalid elements with negative volume can be generated in the worst case. To overcome this problem, the so-called indirect cell-centered Arbitrary-Lagrangian-Eulerian (ALE) algorithms have been developed \cite{ShashkovCellCentered,ShashkovRemap1,ShashkovRemap3,ShashkovRemap4,ShashkovRemap5,MaireMM2}, where the mesh velocity can be chosen \textit{independently} from the local fluid velocity, therefore the grid nodes can be arbitrarily moved. The mesh quality is optimized during the simulation using a remeshing strategy, where a new mesh with better quality is generated, followed by a remapping procedure in which the numerical solution is projected from the old mesh to the new one. Multi-phase and multi-material flow problems are typically solved relying on this approach \cite{ShashkovMultiMat1,ShashkovMultiMat2,ShashkovMultiMat3,ShashkovMultiMat4,Hirt74,Peery2000,Smith1999}.

In a recent series of papers \cite{Lagrange1D,ALELTS1D,Lagrange2D,Lagrange3D,LagrangeMHD,LagrangeMDRS,ALEMOOD1,ALEMOOD2,LagrangeQF,ALEMQF,ALELTS2D} a new family of high order accurate ADER finite volume schemes has been proposed in the ALE context on moving meshes in one and multiple space dimensions. These methods are addressed with \textit{direct ALE} schemes, because the mesh motion is taken into account \textit{directly} in the numerical flux computation of the finite volume scheme, therefore without needing any remeshing plus remapping strategy. High order of accuracy in space is achieved either by the use of a WENO reconstruction technique \cite{DumbserEnauxToro,Dumbser2007204,Dumbser2007693,JiangShu1996,HuShuVortex1999,ZhangShu3D} or by the recent \textit{a posteriori} MOOD paradigm \cite{ALEMOOD1,ALEMOOD2,CDL1,CDL2,CDL3}, while the schemes are allowed to be high order accurate also in time by adopting a local space-time Galerkin predictor method introduced in \cite{DumbserEnauxToro,HidalgoDumbser}, that derives from the ADER approach proposed by Toro et al. \cite{mill,toro3,titarevtoro,DumbserEnauxToro,Dumbser20088209,BalsaraRumpf,CastroToro}. Unstructured \textit{curvilinear} meshes have been recently considered in \cite{LagrangeISO}, while in \cite{LagrangeHPR} such methods have been successfully applied to the equations of nonlinear hyperelasticity. For direct ALE schemes on moving polygonal and polyhedral meshes, see also the very interesting work of Springel \cite{Springel}. 

Another option for the numerical solution of hyperbolic conservation laws is given by Discontinuous Galerkin (DG) methods, first applied to neutron transport equations \cite{reed} and later extended to general nonlinear systems of hyperbolic conservation laws in one a multiple space dimensions in a well-known series of papers by Cockburn and Shu and coworkers \cite{cbs0,cbs1,cbs2,cbs3,cbs4}. Here, the numerical solution is approximated by polynomials within each control volume, hence leading to a natural piecewise high order data representation. Thus, DG schemes do not need any reconstruction procedure, unlike high order finite volume schemes.  These methods are widely used to solve fluid dynamics problems, even in the Lagrangian framework. Finite element algorithms for Lagrangian hydrodynamics and the equations of nonlinear elasto-plasticity have been proposed in \cite{scovazzi1,scovazzi2,Dobrev1,Dobrev2012,Dobrev2013}, while Lagrangian DG methods have been presented for the first time in \cite{Vilar1,Vilar2,Vilar3,Yuetal}. 
In \cite{Vilar1,Vilar2,Vilar3} a so-called \textit{total Lagrangian} approach was chosen, i.e. the computational grid is kept 
\textit{fixed} and the equations of gas dynamics have been written by means of the Lagrangian coordinates related to the initial configuration of the flow. However, as a consequence, within the governing equations one has to take into account also the time evolution of the the Jacobian matrix associated with the mapping of the current configuration to the reference configuration. 
Explicit DG methods as the ones listed so far suffer from a very severe time step restriction, therefore a high order \textit{implicit} time discretization for DG schemes has been proposed in \cite{BassiImplicit1,BassiImplicit2,BassiImplicit3,BassiImplicit4}, while \textit{semi-implicit} DG schemes can be found in \cite{GiraldoRestelli,Dolejsi1,Dolejsi2,Dolejsi3,STINS2D,STINS3D}. 

In this paper we present a new family of high order accurate explicit ADER-DG schemes based on the algorithm proposed in \cite{Lagrange2D,Lagrange3D}, where the computational mesh is \textit{moved} according to the fluid flow and not mapped to the initial configuration as done in \cite{Vilar1,Vilar2,Vilar3}. The method is designed for moving unstructured triangular and tetrahedral meshes in the ALE framework. The use of the local space-time Galerkin predictor naturally permits the development of a one-step algorithm, that is more efficient compared to explicit TVD Runge-Kutta schemes typically adopted for the time integration in the DG context \cite{CBS-book,cbs0,cbs1,cbs2,cbs3,cbs4}. Since DG schemes need some sort of nonlinear limiting to avoid the Gibbs phenomenon at shock waves or other discontinuities, we rely on the recently proposed \textit{a posteriori} sub-cell limiting procedure \cite{DGLimiter1,DGLimiter2,DGLimiter3} which is based on the MOOD paradigm \cite{CDL1,CDL2,CDL3} that has already been used on moving unstructured meshes, see \cite{ALEMOOD1,ALEMOOD2}. \\
The outline of this article is as follows: all the details regarding the proposed numerical method are contained in Section \ref{sec.numethod}, while in Section \ref{sec.validation} we show numerical convergence rates up to fourth order of 
accuracy in space and time for a smooth problem as well as a wide set of benchmark test problems considering inviscid and viscous compressible flows. An application close to Inertial Confinement Fusion (ICF) simulation is also presented at the end of this manuscript and the ideal classical and viscous resistive (magnetohydrodynamics) MHD equations are considered. Finally, we give some concluding remarks and an outlook to possible future work in Section \ref{sec.concl}. 

\section{The ADER Discontinuous Galerkin method on moving unstructured meshes}
\label{sec.numethod} 

\begin{figure}[!htbp]
\begin{center}
\begin{tabular}{cc} 
\includegraphics[width=0.4\textwidth]{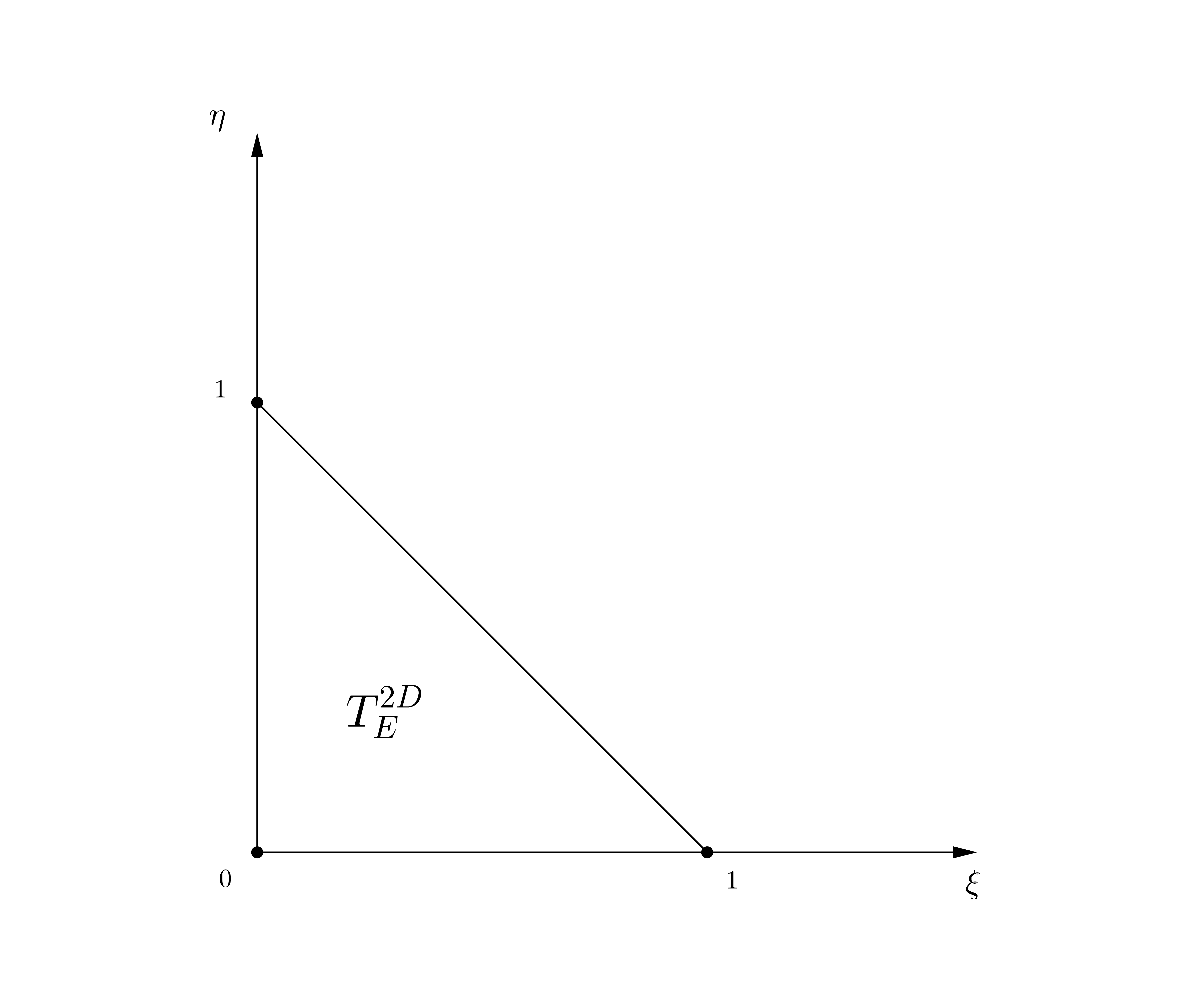}  &           
\includegraphics[width=0.4\textwidth]{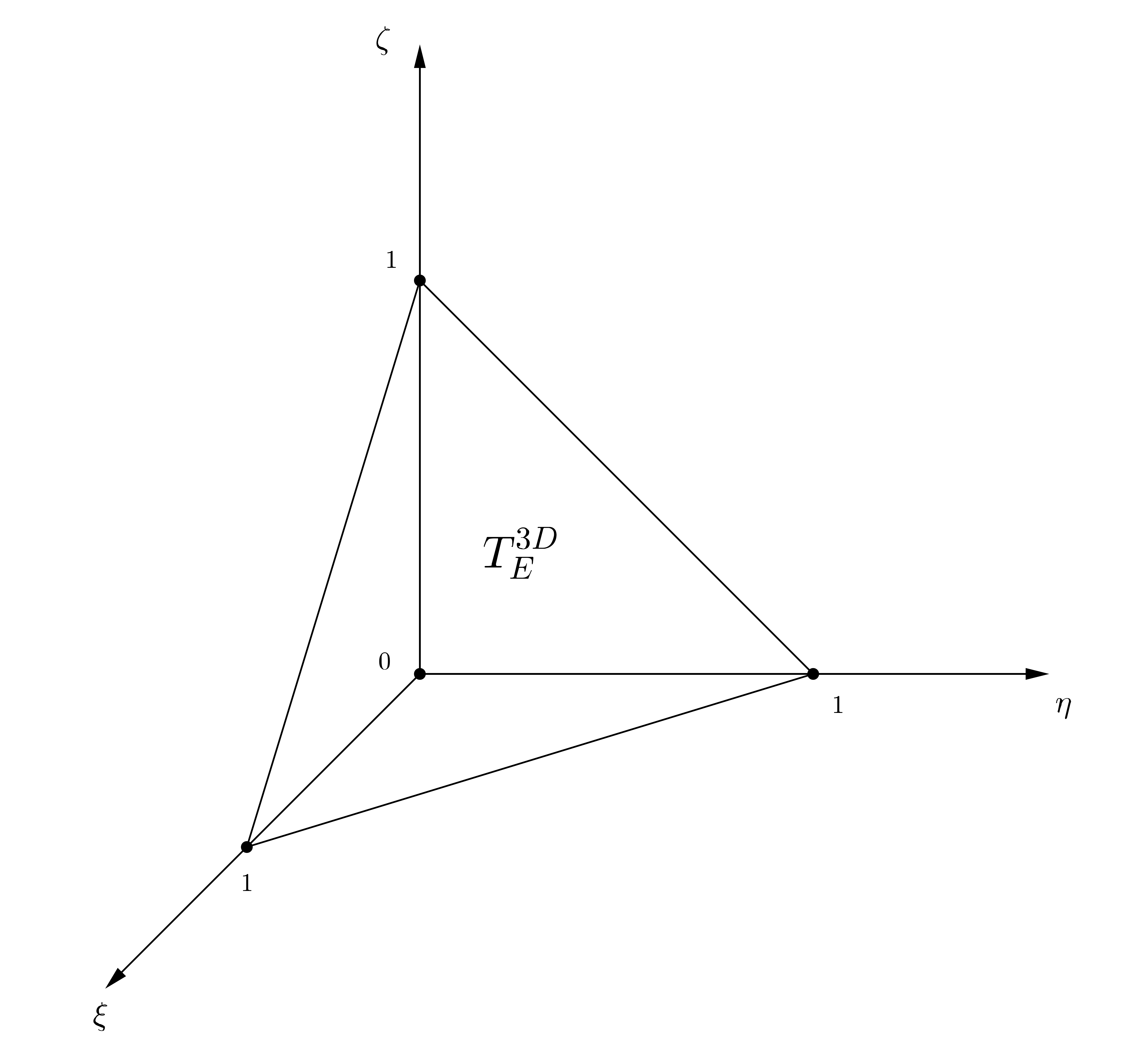} \\           
\end{tabular} 
\caption{Reference element in 2D (left) and in 3D (right) used to define the Dubiner-type basis functions $\phi_l$ in \eqref{eqn.uh}.}
\label{fig.Te}
\end{center}
\end{figure}

In this paper we consider nonlinear homogeneous systems of conservation laws of the form 
\begin{equation}
\label{PDE}
 \frac{\partial \Q}{\partial t} + \nabla \cdot \F(\Q,\nabla \Q) = \mathbf{0}, \qquad \x \in \Omega(t) \subset \mathds{R}^d, \quad t \in \mathds{R}_0^+, \quad \Q \in \Omega_{\Q} \subset \mathds{R}^\nu,     
\end{equation} 
with $\Q$ denoting the vector of conserved variables defined in the space of the admissible states $\Omega_{\Q}\subset \mathds{R}^\nu$ and $\F(\Q,\nabla \Q)=\left( \f(\Q,\nabla \Q),\g(\Q,\nabla \Q),\h(\Q,\nabla \Q) \right)$ representing the nonlinear flux tensor which depends on the state $\Q$ and its gradient $\nabla \Q$. The computational domain $\Omega(t)$ is defined in $d\in[2,3]$ space dimensions by the spatial coordinate vector $\x=(x,y,z)$ and in the ALE framework it is \textit{time-dependent}, hence continuously changing its configuration. At the current time $t^n$ a total number $N_E$ of non-overlapping unstructured control volumes $T_i^n$ is used to discretize the domain $\Omega$, yielding the \textit{current mesh configuration} $\mathcal{T}^n_{\Omega}$:
\begin{equation}
\mathcal{T}^n_{\Omega} = \bigcup \limits_{i=1}^{N_E}{T^n_i}. 
\label{trian}
\end{equation} 
The elements are chosen to be piecewise straight or curved simplex control volumes, i.e generalized triangles and tetrahedra in two and three space dimensions, respectively. \\
The numerical solution for the state vector $\Q$ in \eqref{PDE} is represented within each cell $T_i^n$ at the current time $t^n$ by piecewise polynomials of degree $N \geq 0$ denoted by $\mathbf{u}_h(\x,t^n)$ and defined in the space $\mathcal{U}_h$. Thus, the discrete representation of the solution is written as
\begin{equation}
\mathbf{u}_h(\x,t^n) = \sum \limits_{l=1}^\mathcal{N} \phi_l(\x) \hat{\mathbf{u}}^{n}_{l} \qquad \x \in T_i^n,
\label{eqn.uh}
\end{equation}
where $\phi_l(\x)$ is a set of spatial basis functions used to span the space $\mathcal{U}_h$ up to degree $N$. In the rest of the paper we will use classical tensor index notation based on the Einstein summation convention, which implies summation over two equal indices. The total number $\mathcal{N}$ of expansion coefficients (degrees of freedom) $\hat{\mathbf{u}}^{n}_{l}$ for the basis functions depends on the polynomial degree $N$ and is given by
\begin{equation}
  \mathcal{N} = \mathcal{N}(N,d) = \frac{1}{d!} \prod \limits_{m=1}^{d} (N+m).
\label{eqn.nDOF}
\end{equation}
The  Dubiner-type basis functions \cite{Dubiner,orth-basis,CBS-book} are used as basis functions $\phi_l$ in \eqref{eqn.uh} and they are defined on the \textit{reference element} $T_E$ in the reference coordinate system $\boldsymbol{\xi}=(\xi,\eta,\zeta)$. The reference element is depicted in Figure \ref{fig.Te} and it is the unit triangle in 2D, defined by vertices $\boldsymbol{\xi}_{E,1}=(\xi_{E,1},\eta_{E,1})=(0,0)$, $\boldsymbol{\xi}_{E,2}=(\xi_{E,2},\eta_{E,2})=(1,0)$ and $\boldsymbol{\xi}_{E,3}=(\xi_{E,3},\eta_{E,3})=(0,1)$, or the unit tetrahedron in 3D with nodes $\boldsymbol{\xi}_{E,1}=(\xi_{E,1},\eta_{E,1},\zeta_{E,1})=(0,0,0)$, $\boldsymbol{\xi}_{E,2}=(\xi_{E,2},\eta_{E,2},\zeta_{E,2})=(1,0,0)$, $\boldsymbol{\xi}_{E,3}=(\xi_{E,3},\eta_{E,3},\zeta_{E,3})=(0,1,0)$ and $\boldsymbol{\xi}_{E,4}=(\xi_{E,4},\eta_{E,4},\zeta_{E,4})=(0,0,1)$.

The governing equations \eqref{PDE} are solved at the aid of a high order ADER-DG (Discontinuous Galerkin) algorithm \cite{QiuDumbserShu,ADERNC} which is based on a one-step predictor-corrector method presented in \cite{Dumbser2008}. The ADER predictor step solves system \eqref{PDE} \textit{locally} (in the small) by considering the space-time evolution of the conservation law within each space-time element, while the corrector step is given by directly integrating a weak form of the governing PDE on a set of space-time control volumes. The scheme provides high order of accuracy in space and time in one single time step $\Delta t$, which is evaluated under a classical (global) Courant-Friedrichs-Levy number (CFL) stability condition as
\begin{equation}
\Delta t < \frac{\textnormal{CFL}}{2N+1} \, \min \limits_{T_i^n} \frac{h_i}{|\lambda_{\max,i}|}, \qquad \forall T_i^n \in \Omega^n. 
\label{eq:timestep}
\end{equation}
The characteristic element size $h_i$ is taken to be either the incircle or the insphere diameter for triangles or tetrahedra, respectively, while $|\lambda_{\max,i}|$ is given by the maximum absolute value of the eigenvalues computed from the current solution $\Q_i^n$ in $T_i^n$. On unstructured meshes the CFL stability condition requires the inequality $\textnormal{CFL} \leq \frac{1}{d}$ to be satisfied. For high order Lagrangian schemes with time-accurate local time stepping (LTS), see \cite{ALELTS1D,ALELTS2D}. 

It is well known, the DG method needs some sort of nonlinear limiting to avoid the Gibbs phenomenon at shock waves or other discontinuities. In our approach we rely on the \textit{a posteriori} sub-cell finite volume limiter recently developed in \cite{DGLimiter1,DGLimiter2,DGLimiter3}. It is based on a low order finite volume scheme that acts on a fine sub-grid onto which the numerical solution $\mathbf{u}_h(\x,t^n)$ is scattered when needed. First, we illustrate how the sub-grid is built, then we briefly recall the ADER predictor step and we present a novel strategy to move the mesh to the next time level. Then, the ADER-DG corrector strategy is described and finally we provide an overview of the \textit{a posteriori} sub-cell limiter adopted in the ALE context on moving unstructured meshes.

\subsection{Piecewise linear sub-cell element description needed for the limiter} 
\label{sec.subgrid}

\begin{figure}[!htbp]
  \begin{center}
    \begin{tabular}{cc} 
      \includegraphics[width=0.4\textwidth]{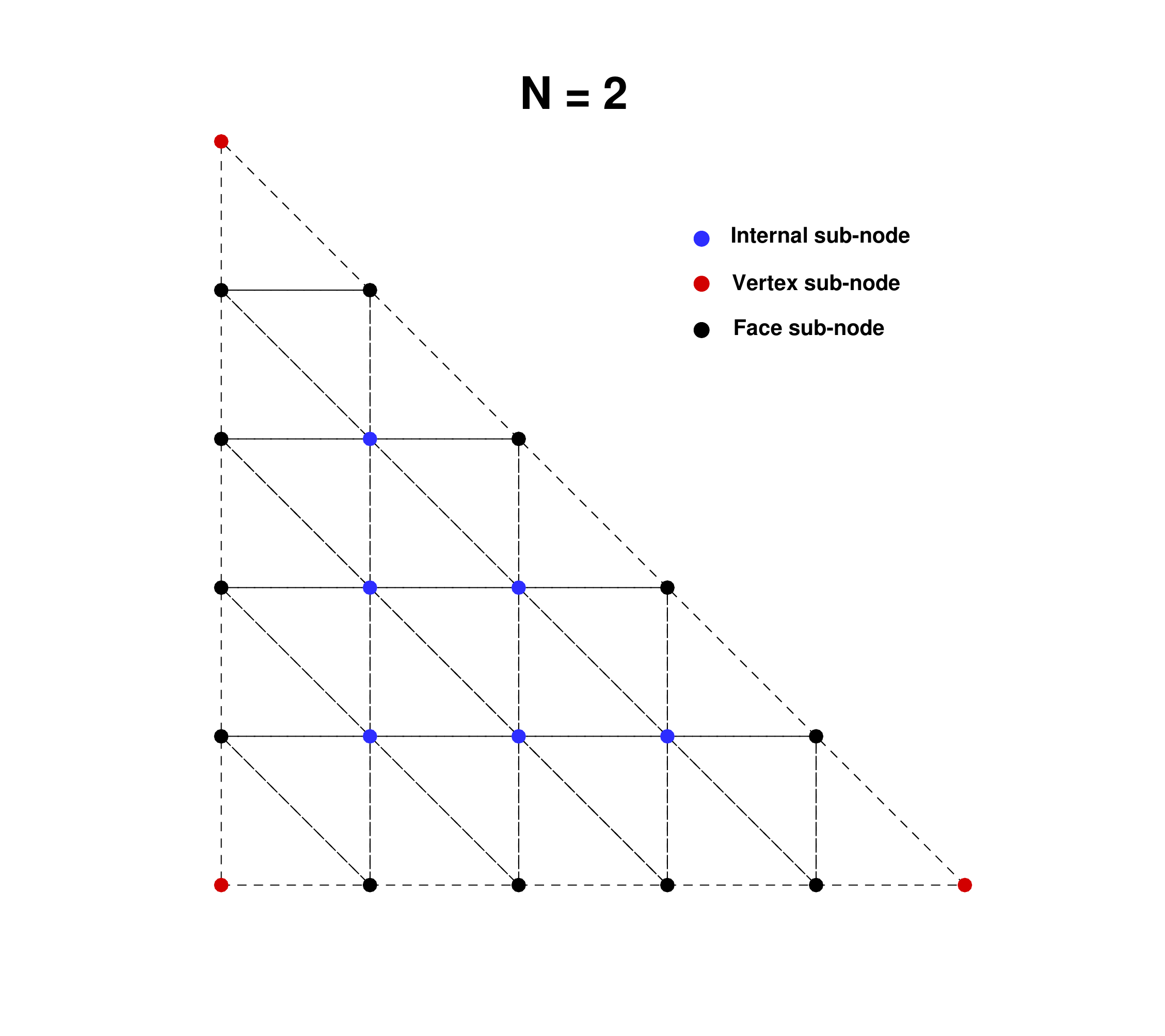} &
			\includegraphics[width=0.4\textwidth]{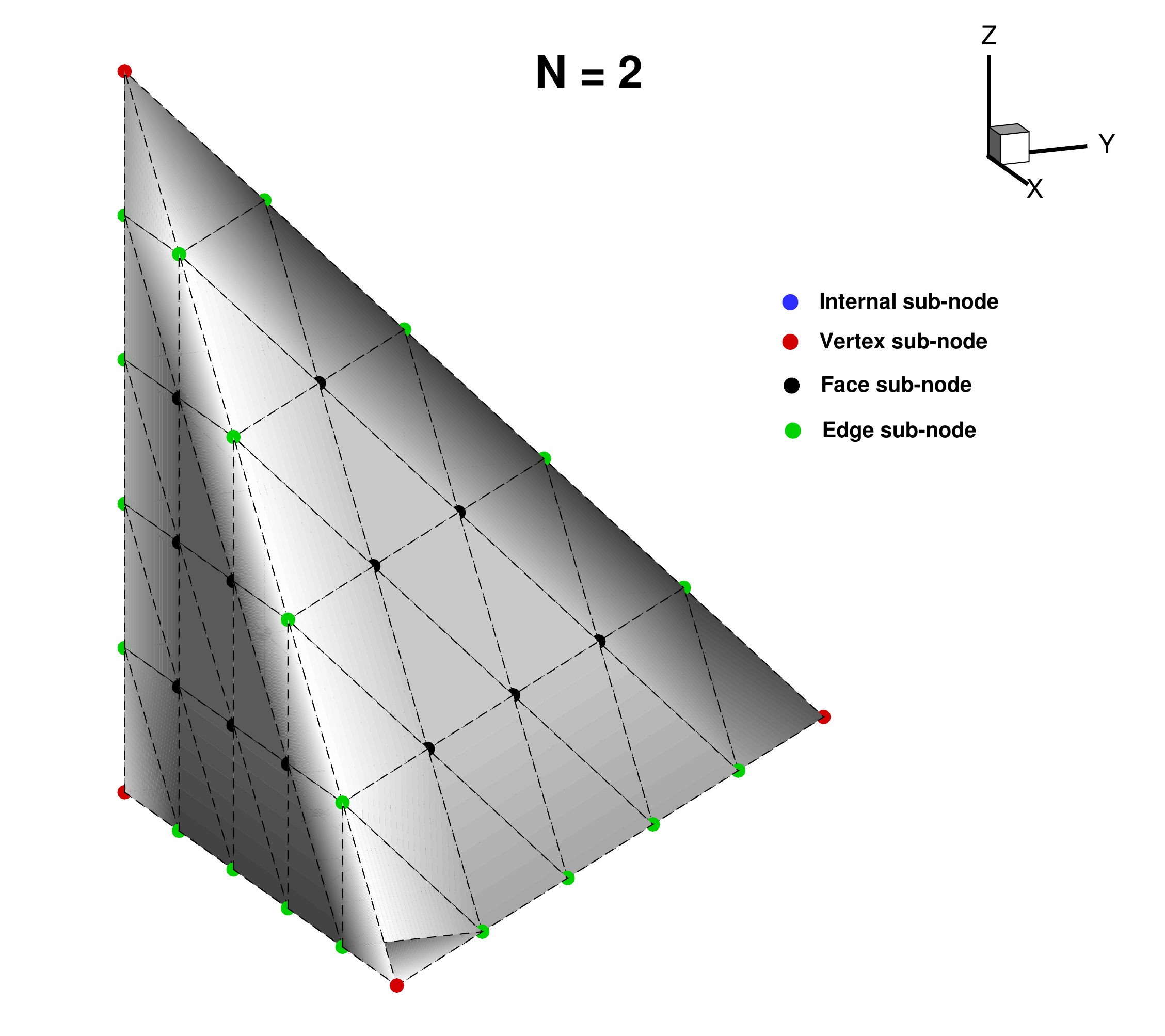} \\
			\includegraphics[width=0.4\textwidth]{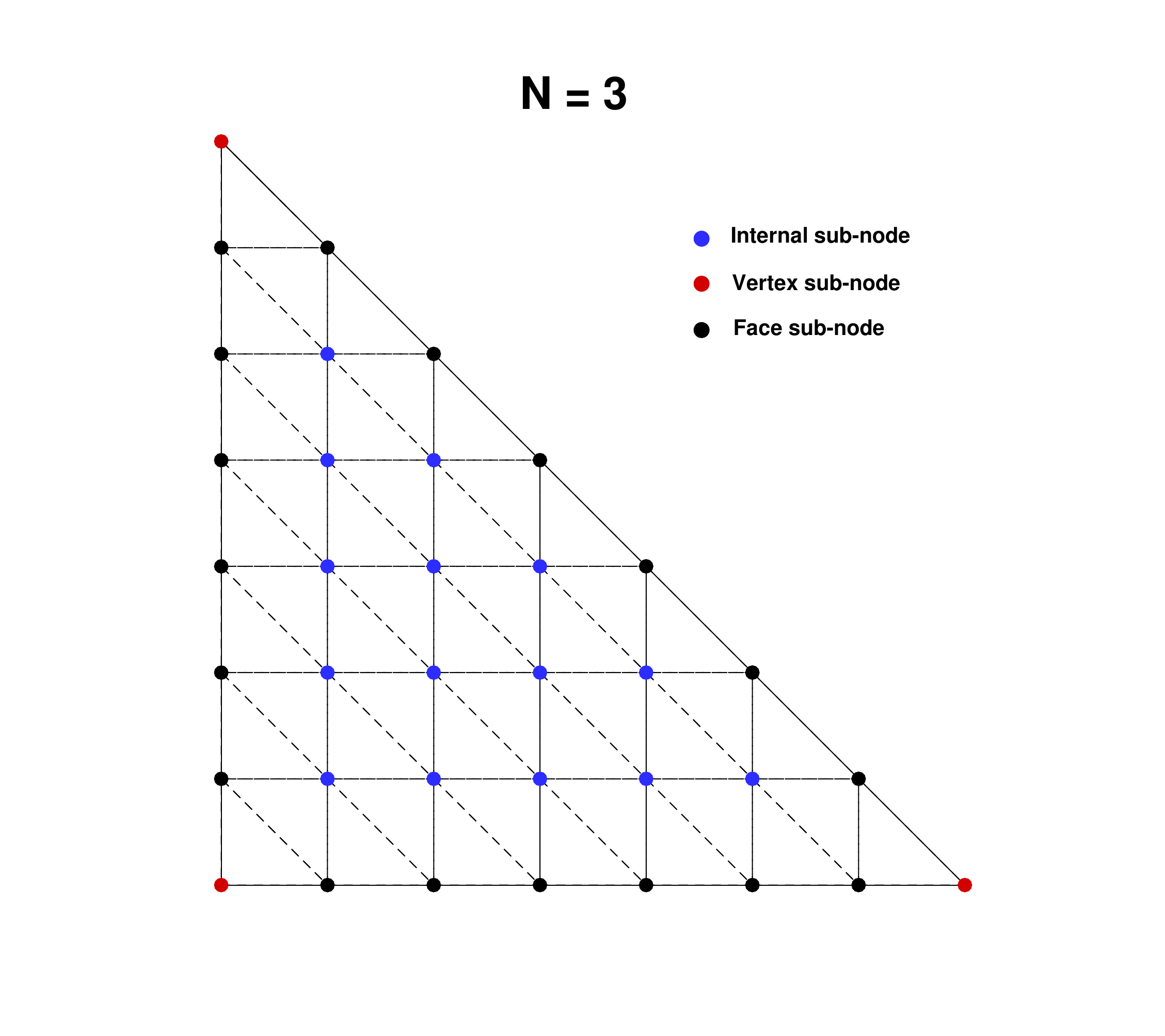} &
      \includegraphics[width=0.4\textwidth]{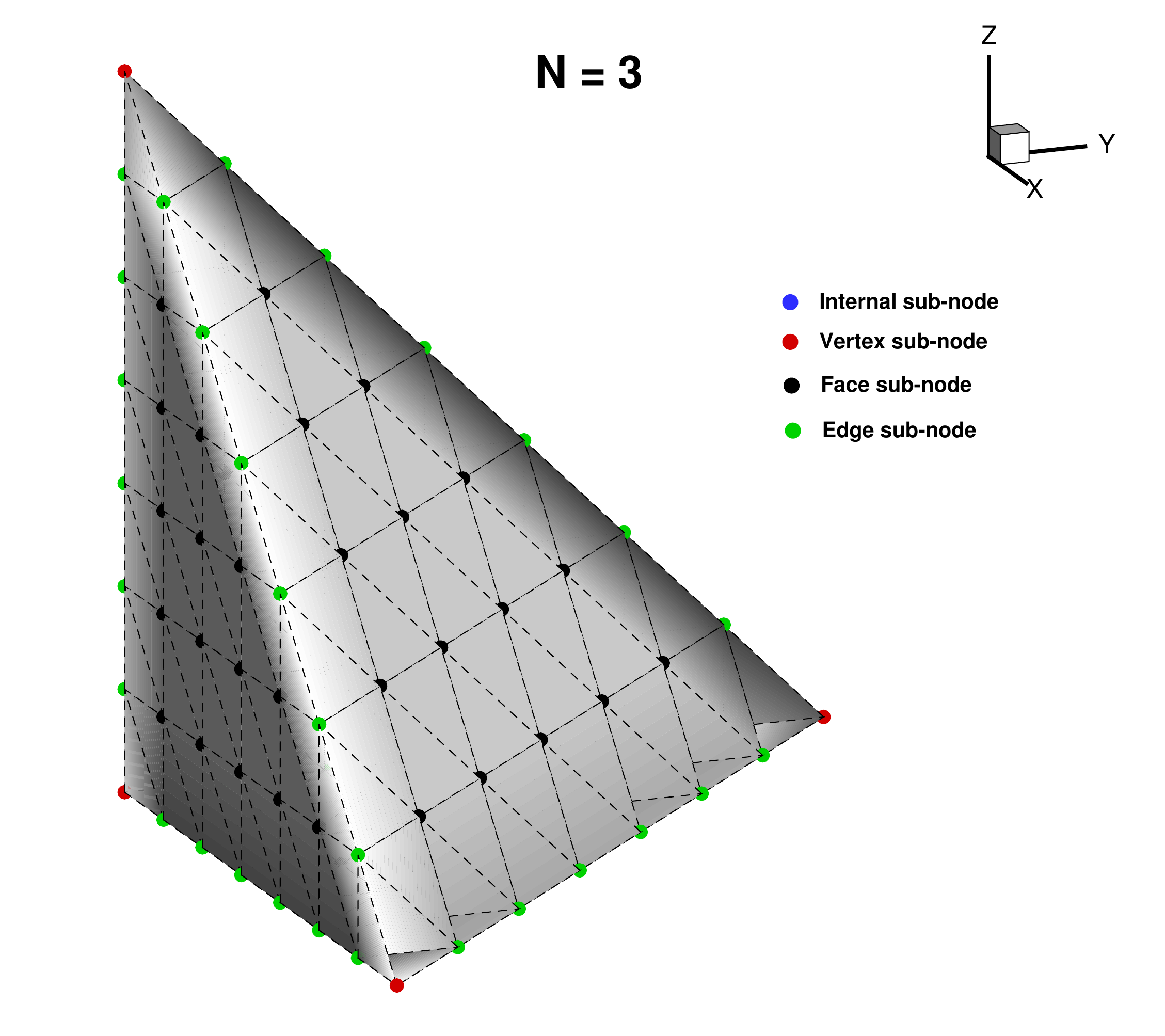} \\
			\includegraphics[width=0.4\textwidth]{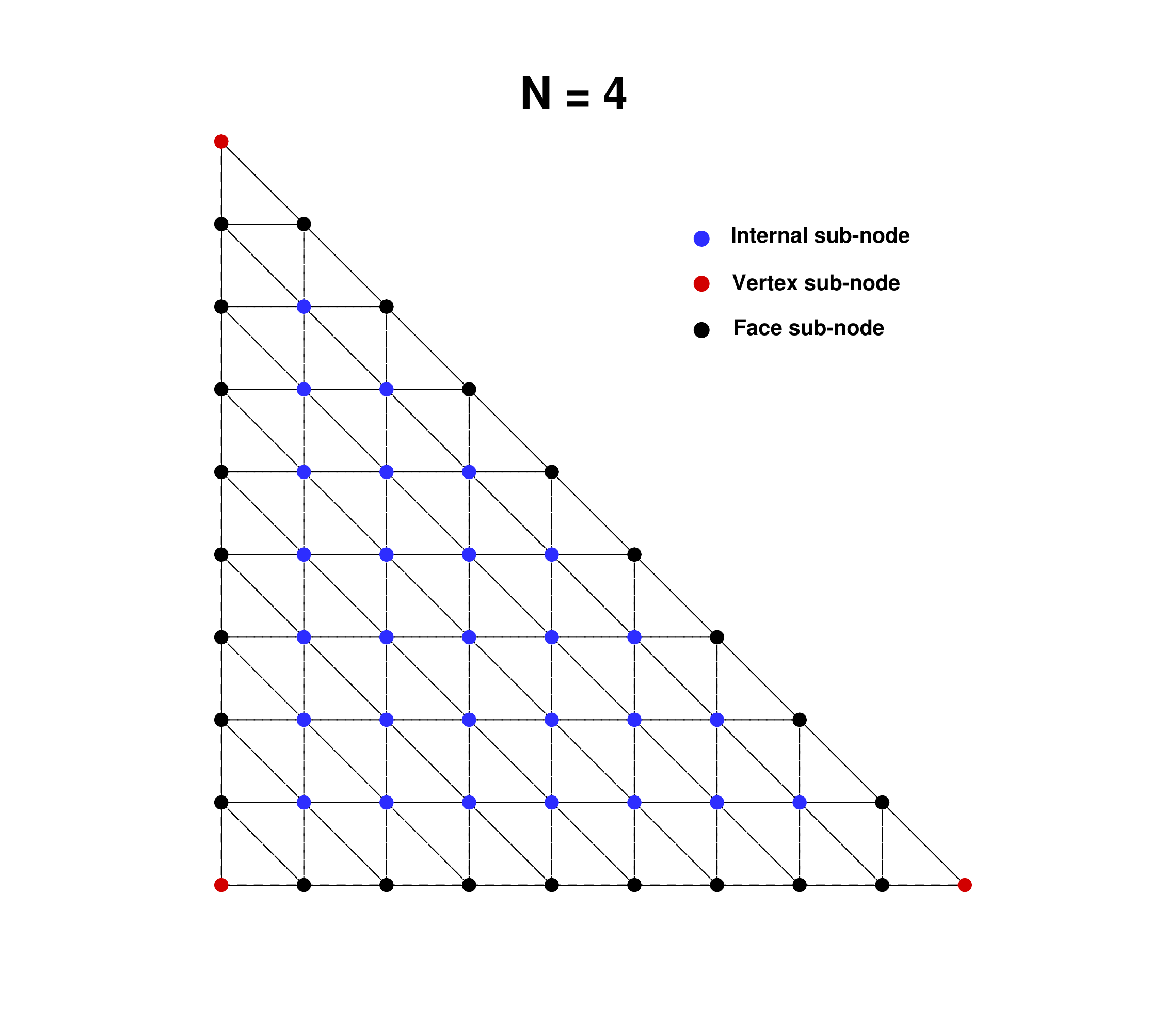} &
      \includegraphics[width=0.4\textwidth]{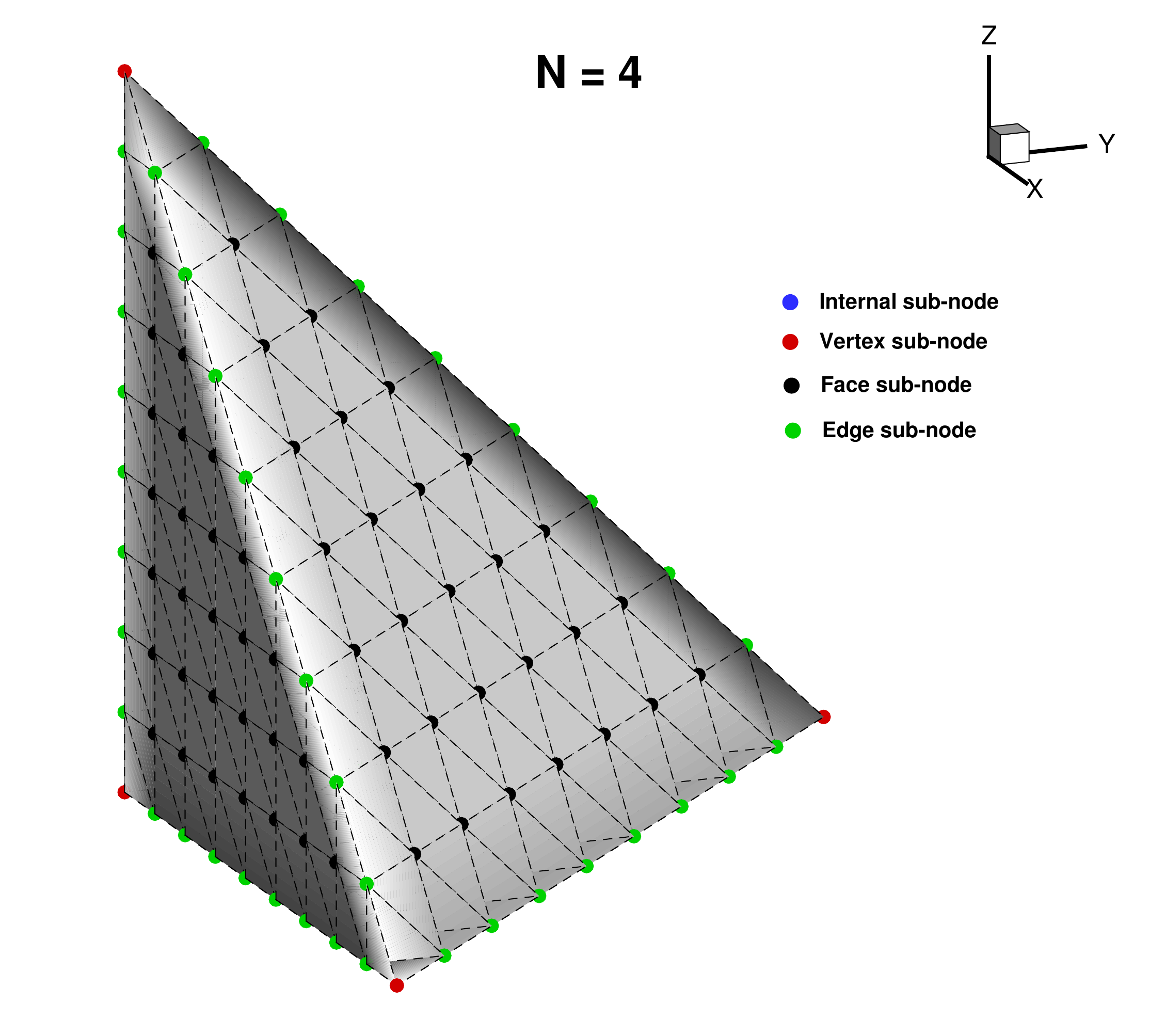} \\
		\end{tabular} 
    \caption{Sub-grid configuration on the reference element $T_E$ for $N=2,3,4$ (from top to bottom row) in two (left column) and three (right column) space dimensions. The types of sub-node (internal, vertex, face and edge sub-node) are highlighted with different colors.}
    \label{fig.SubGrid}
	\end{center}
\end{figure}

For so-called \textit{troubled cells}, i.e. for those cells which need limiting, the element shape is described by means of a set of \textit{sub-cells} arising from the splitting of each element edge into $N_s=2N+1$ sub-edges, as done in \cite{DGLimiter3}. The sub-grid is built in the reference element $T_E$, as shown in Figure \ref{fig.SubGrid}, and is composed by a total number of $\mathcal{S}=(N_s)^d$ sub-cells which are defined by $\mathcal{K}$ sub-nodes, whose coordinates $\kk$ are provided by the standard nodes of classical high order conforming finite elements on triangular and tetrahedral meshes, therefore
\begin{equation}
 \kk_{k,p}^{2D} = \left(\frac{k}{N_s},\frac{p}{N_s}\right) \quad \textnormal{ and } \quad \kk_{m,k,p}^{3D} = \left(\frac{m}{N_s},\frac{k}{N_s},\frac{p}{N_s}\right)
\label{eqn.xiSN}
\end{equation}
with
\begin{equation}
0 \leq p \leq N_s, \qquad 0 \leq k \leq (N_s-p), \qquad 0 \leq m \leq (N_s-p-k).
\label{eqn.idxSN}
\end{equation}
The total number of sub-nodes is given by \eqref{eqn.nDOF}, hence $\mathcal{K}^{2D}=\mathcal{N}(N_s,2)=(N_s+1)(N_s+2)/2$ and $\mathcal{K}^{3D}=\mathcal{N}(N_s,3)=(N_s+1)(N_s+2)(N_s+3)/6$. Each sub-cell $S_{k,p}^{2D}$ and $S_{m,k,p}^{3D}$ is assigned a local connectivity specified in \cite{DGLimiter3} and, to ease the notation, we will refer to sub-cell $S_{k,p}$ in 2D or $S_{m,k,p}$ in 3D of element $T_i^n$ simply with $S_{ij}^n$. The same applies to the vertexes, thus $\kk_{k,p}$ becomes $\kk_{ij}$. 

We emphasize that the subgrid description of the geometry at the aid of a \textit{piecewise linear} simplex subgrid essentially corresponds to the \textit{agglomeration approach} recently forwarded by  Bassi et al. in \cite{BassiAgglomeration}. 
The use of the sub-grid also allows us to introduce an \textit{alternative data representation} $\mathbf{v}_h(\x,t^n)$ given by a set of piecewise constant sub-cell averages $\mathbf{v}_{ij}^n$ that are computed according to \cite{DGLimiter3} as
\begin{equation}
\mathbf{v}_{ij}(\x,t^n) = \frac{1}{|S_{ij}^n|} \int \limits_{S_{ij}^n} \mathbf{u}_{h}(\x,t^n) \, d\x = \frac{1}{|S_{ij}^n|} \int \limits_{S_{ij}^n} \phi_l(\x) d\x \, \hat{\mathbf{u}}^{n}_{l} \qquad \forall j \in [1,\mathcal{S}],
\label{eqn.vh}
\end{equation}
where $|S_{ij}^n|$ denotes the volume of sub-cell $S_{ij}$ of element $T_i^n$. The $L_2$ projection operator \eqref{eqn.vh} is defined by $\mathbf{v}_{h}^n:=\mathcal{P}(\mathbf{u}_h^n)$ and it can be computed and stored once and for all in the pre-processing step because the sub-grid connectivity as well as the spatial basis functions $\phi_l(\x)$ are defined on the reference element $T_E$ in $\xxi$ which does not change in time.

The \textit{reconstruction operator} is the inverse of the projection \eqref{eqn.vh} and it permits to recover the piecewise polynomial solution $\mathbf{u}_{h}(\x,t^n)$ of the DG scheme on the main grid. This is done solving a classical reconstruction problem, where one requires the following condition to be satisfied:
\begin{equation}
\int \limits_{S_{ij}^n} \mathbf{u}_{h}(\x,t^n) \, d\x = \int \limits_{S_{ij}^n} \mathbf{v}_{h}(\x,t^n) \, d\x  \qquad \forall j \in [1,\mathcal{S}].
\label{eqn.intRec}
\end{equation}
Due to the choice of taking $N_s=2N+1 \geq N+1$, equation \eqref{eqn.intRec} in general leads to an overdetermined linear system that is solved using a constrained least-squares technique \cite{DumbserKaeser06b} in which the reconstruction is imposed to be \textit{conservative} on the main cell $T_i^n$, hence yielding the additional linear constraint
\begin{equation}
\int \limits_{T_i^n} \mathbf{u}_{h}(\x,t^n) \, d\x = \int \limits_{T_i^n} \mathbf{v}_{h}(\x,t^n) \, d\x.
\label{eqn.LSQ}
\end{equation}
The reconstruction operator is shortened by $\mathbf{u}_{h}^n:=\mathcal{R}(\mathbf{v}_h^n)$ and it is also defined on the reference element $T_E$, so that system \eqref{eqn.intRec} can be written as
\begin{equation}
\frac{1}{|S_{ij}^n|} \int \limits_{S_{ij}^n} \phi_l(\x) d\x \, \hat{\mathbf{u}}^{n}_{l} = \mathbf{v}_{ij}(\x,t^n),
\label{eqn.intRec2}
\end{equation}
and the reconstruction matrix given by the integrals on the left hand side is conveniently computed and stored once at the beginning of the simulation.

\subsection{Local space-time predictor} 
\label{sec.lst} 
The ADER approach is based on the solution of the generalized Riemann problem, which requires the time derivatives, that are needed to evolve the solution in time, to be computed from the governing PDE \eqref{PDE} in terms of spatial derivatives. Here, the local space-time predictor aims at providing an element-local predictor solution of the PDE without needing any neighbor information. This strategy has been successfully developed and applied in the Eulerian framework on fixed grids in \cite{DumbserEnauxToro,Dumbser20088209,USFORCE2,HidalgoDumbser} and subsequently extended to moving meshes in the ALE context \cite{Lagrange2D,LagrangeNC,Lagrange3D,LagrangeMHD,LagrangeMDRS,ALELTS2D,ALEMOOD1,ALEMOOD2,LagrangeISO}. The starting point of the local-space time strategy is given by the polynomials which represent the numerical solution at the current time $t^n$ that will be then evolved \textit{locally} up to the next time level $t^{n+1}$ within the space-time control volume $\tilde{T}_i = T_i(t) \times \left[t^n,t^{n+1}\right]$. In the finite volume framework, the starting polynomials are obtained via reconstruction from the known cell averages of the conserved quantities, while in the Discontinuous Galerkin approach they are directly available from \eqref{eqn.uh}. As a result we obtain piecewise space-time polynomials $\q_h(\x,t)$ of degree $N$, which will then be employed in the corrector step described in Section \ref{sec.SolAlg} for computing a numerical flux function (Riemann solver) that provides the coupling between neighbor elements.  

\begin{figure}[!htbp]
  \begin{center}
    \begin{tabular}{cc} 
      \includegraphics[width=0.3\textwidth]{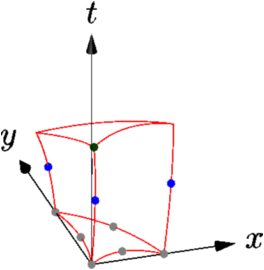} &
			\includegraphics[width=0.3\textwidth]{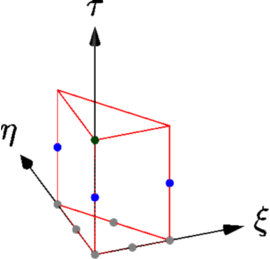} \\
		\end{tabular} 
    \caption{Physical (left) and reference (right) space-time element in 2D with the corresponding space-time nodes for $N=2$. The element configuration as well as the predictor solution are approximated using a set of isoparametric basis functions $\theta(\tilde{\xxi})$ of degree $N$.}
    \label{fig.Te-lcg}
	\end{center}
\end{figure}

The local predictor strategy is based on an element-local \textit{weak formulation} of the governing PDE \eqref{PDE} in space and time, which reads
\begin{eqnarray}
\int \limits_{t^n}^{t^{n+1}} \int \limits_{T_i(t)} \theta_k \, \frac{\partial \q_h}{\partial t}\, d\x \, dt + \int \limits_{t^n}^{t^{n+1}} \int \limits_{T_i(t)} \theta_k \, \nabla \cdot \F(\q_h,\nabla \q_h)\, d\x \, dt = \mathbf{0}, 
\label{eqn.PDEweak}
\end{eqnarray}
where the time step $\Delta t=t^{n+1}-t^n$ is given by \eqref{eq:timestep}. In the above expression, $\theta_k=\theta_k(\x,t)$ are a set of space-time test functions of degree $N$ that are also used to approximate the predictor solution $\q_h(\x,t)$, hence
\begin{equation}
  \q_h(\x,t) = \sum \limits_{l=1}^\mathcal{L} \theta_l(\tilde{\boldsymbol{\xi}}) \hat \q^{n}_{l,i}:= \theta_l \hat \q^{n}_{l,i}.
\label{eqn.qh}
\end{equation}
According to \cite{Lagrange2D,Lagrange3D} the basis functions $\theta_l(\tilde{\boldsymbol{\xi}})$ are defined by the Lagrange interpolation polynomials passing through a set of space-time nodes $\hat{\boldsymbol{\xi}}_l$ specified in \cite{Dumbser2008}, yielding a \textit{nodal} basis. $\mathcal{L}$ represents the total number of degrees of freedom and it is given by \eqref{eqn.nDOF} with $d+1$ dimensions, since also time is now considered. The symbol tilde ($\tilde{}$) is used for space-time quantities and the space-time basis and test functions as well as the integrals appearing in \eqref{eqn.PDEweak} are conveniently defined on the space-time reference element $T_E \times [0,1]$ with $\tilde{\boldsymbol{\xi}}=(\xi,\eta,\zeta,\tau)$, shown in Figure \ref{fig.Te-lcg}.

In order to take into account the initial condition, which is given by the known polynomials $\mathbf{u}_h(\x,t^n)$, the first term of \eqref{eqn.PDEweak} is integrated by parts in time leading to
\begin{equation}
\int \limits_{T_E} \theta_k(\boldsymbol{\xi},1) \q_h \, d\boldsymbol{\xi} - \int \limits_{T_E} \theta_k(\boldsymbol{\xi},0) \mathbf{u}_h \, d\boldsymbol{\xi} - \int \limits_{0}^{1} \int \limits_{T_E} \frac{\partial \theta_k}{\partial \tau} \q_h \, d\boldsymbol{\xi} d\tau + \int \limits_{0}^{1} \int \limits_{T_E} \theta_k \nabla_{\xi} \cdot \F(\q_h, \nabla \q_h) \, d\boldsymbol{\xi} d\tau = \mathbf{0}, 
\label{eqn.PDElst}
\end{equation}
where the integrals are defined in the space-time reference element and are evaluated using multidimensional Gaussian quadrature rules of suitable order of accuracy, see \cite{stroud} for details.

In the ALE context the mesh is moving in time, thus changing the geometry of the space-time control volume $\tilde{T}_i=T_i(t) \times \Delta t$. The mesh motion is governed by the trajectory equation  
\begin{equation}
\frac{d\x}{dt} = \mathbf{V}(\x,t),
\label{ODEmesh}
\end{equation}
where $\mathbf{V}(\x,t)=(U,V,W)$ is the local mesh velocity. We adopt an \textit{isoparametric} approach, where the \textit{same} space-time basis functions $\theta_l$, used for the approximation of the predictor solution $\q_h$, are also employed to discretize the element geometry configuration as well as the mesh velocity, therefore
\begin{eqnarray}
  \tilde{\x}_h(\tilde{\boldsymbol{\xi}}) &=& \sum \limits_{l=1}^\mathcal{L} \theta_l(\tilde{\boldsymbol{\xi}}) \tilde{\hat{\x}}^{n}_{l,i}:= \theta_l \tilde{\hat{\x}}^{n}_{l,i}, \label{eqn.stgeom1} \\
	\mathbf{V}_h(\tilde{\boldsymbol{\xi}}) &=& \sum \limits_{l=1}^\mathcal{L} \theta_l(\tilde{\boldsymbol{\xi}}) \hat{\mathbf{V}}^{n}_{l,i}:= \theta_l \hat{\mathbf{V}}^{n}_{l,i}, \label{eqn.stgeom2}
\end{eqnarray}
where $\tilde{\hat{\x}}^{n}_{l,i}$ and $\hat{\mathbf{V}}^{n}_{l,i}$ are the space-time coordinates and the corresponding velocities which, in the nodal approach, also provide the degrees of freedom of the expansions \eqref{eqn.stgeom1}-\eqref{eqn.stgeom2}. The trajectory equation \eqref{ODEmesh}, i.e. the time evolution of the element configuration, must be computed \textit{together} with the space-time predictor solution $\q_h$ given by the nonlinear equation \eqref{eqn.PDElst}. Such coupled system is solved by using an iterative procedure which stops when the residuals of the two systems are less than a prescribed tolerance $tol$ (typically $tol\approx 10^{-12}$). 

Once the above procedure is performed for all cells, an element-local predictor for the numerical solution $\q_h$, for the mesh velocity $\mathbf{V}_h$ as well as for the element configuration $\tilde{\x}_h$ is available.

\subsection{Mesh motion}
\label{sec.meshMot}
In the ALE framework the computational mesh changes its configuration $\mathcal{T}_{\Omega}$ at each time step, hence requiring a procedure to determine how the control volumes move in time. The local predictor strategy described in the previous section provides a high order predictor solution $\q_h$ as well as a high order isoparametric description of the element configuration $\tilde{\x}_h$, which has been computed \textit{locally}. As a consequence, the mesh configuration at the new time level $\mathcal{T}_{\Omega}^{n+1}$ might be discontinuous, due to the different local evolution of each space-time control volume $\tilde{T}_i$. In order to recover mesh continuity at time $t^{n+1}$, we rely on a \textit{nodal solver} algorithm. It is a widespread technique used in Lagrangian numerical schemes \cite{chengshu1,chengshu2,chengshu3,Mai07,PH09,phm109,maire_loubere_vachal10,MaireRezoning,DepresMazeran2003,Des05} which aims at evaluating a \textit{unique} velocity vector $\overline{\mathbf{V}}$ for each geometrical degree of freedom. If elements are bounded by straight edges, i.e. a piecewise linear description is adopted, such degrees of freedom are simply given by the vertexes of each cell \cite{Lagrange2D,Lagrange3D,LagrangeMHD}, while if the control volumes are defined by a high order geometry involving \textit{curvilinear} boundaries, as done in \cite{LagrangeISO}, one has to fix a velocity vector also for all the other corresponding degrees of freedom. In any case the velocity vectors $\overline{\mathbf{V}}$ allow  the \textit{Lagrangian}  mesh configuration $\mathcal{T}_{\Omega}^{Lag}$ to be determined, that is the geometry of the computational domain at the next time level obtained solving \textit{locally} the trajectory equation \eqref{ODEmesh} and applying \textit{globally} a nodal solver algorithm. Such a configuration might lead to highly compressed, twisted or even tangled control volumes if the fluid or the grid motion involves very complex flow patterns as vortexes, shock waves or other discontinuities. This is why the Lagrangian phase is typically followed by a \textit{rezoning} strategy which improves the local and global mesh quality, generating the rezoned mesh configuration $\mathcal{T}_{\Omega}^{Rez}$. Finally, the \textit{new} triangulation or tetrahedrization $\mathcal{T}_{\Omega}^{n+1}$ is given as a linear combination between the Lagrangian and the rezoned position of the mesh degrees of freedom, where the blending factor is evaluated according to the \textit{relaxation} algorithm proposed in \cite{MaireRezoning}.

In the following, we present separately the three steps needed for obtaining the final new mesh configuration $\mathcal{T}_{\Omega}^{n+1}$, namely the Lagrangian phase, the rezoning phase and the relaxation phase. If the local mesh velocity $\mathbf{V}(\x,t)=(U,V,W)$ is prescribed and known \textit{a priori}, then we do not need any of the aforementioned strategies and the new mesh configuration is simply obtained by
\begin{equation}
\mathbf{X}_{k}^{n+1} = \mathbf{X}_{k}^{n} + \bar{\mathbf{V}}_k \cdot \Delta t.
\label{eqn.ALEmesh}
\end{equation}
Next, we will discuss the more interesting case in which a Lagrangian-like mesh motion is solved by the trajectory equation \eqref{ODEmesh}, where the local mesh velocity is chosen to be equal to the local fluid velocity, i.e $\mathbf{V}=\mathbf{v}$.
 
\subsubsection{The Lagrangian step}
\begin{figure}[!htbp]
  \begin{center}
    \begin{tabular}{c} 
      \includegraphics[width=0.9\textwidth]{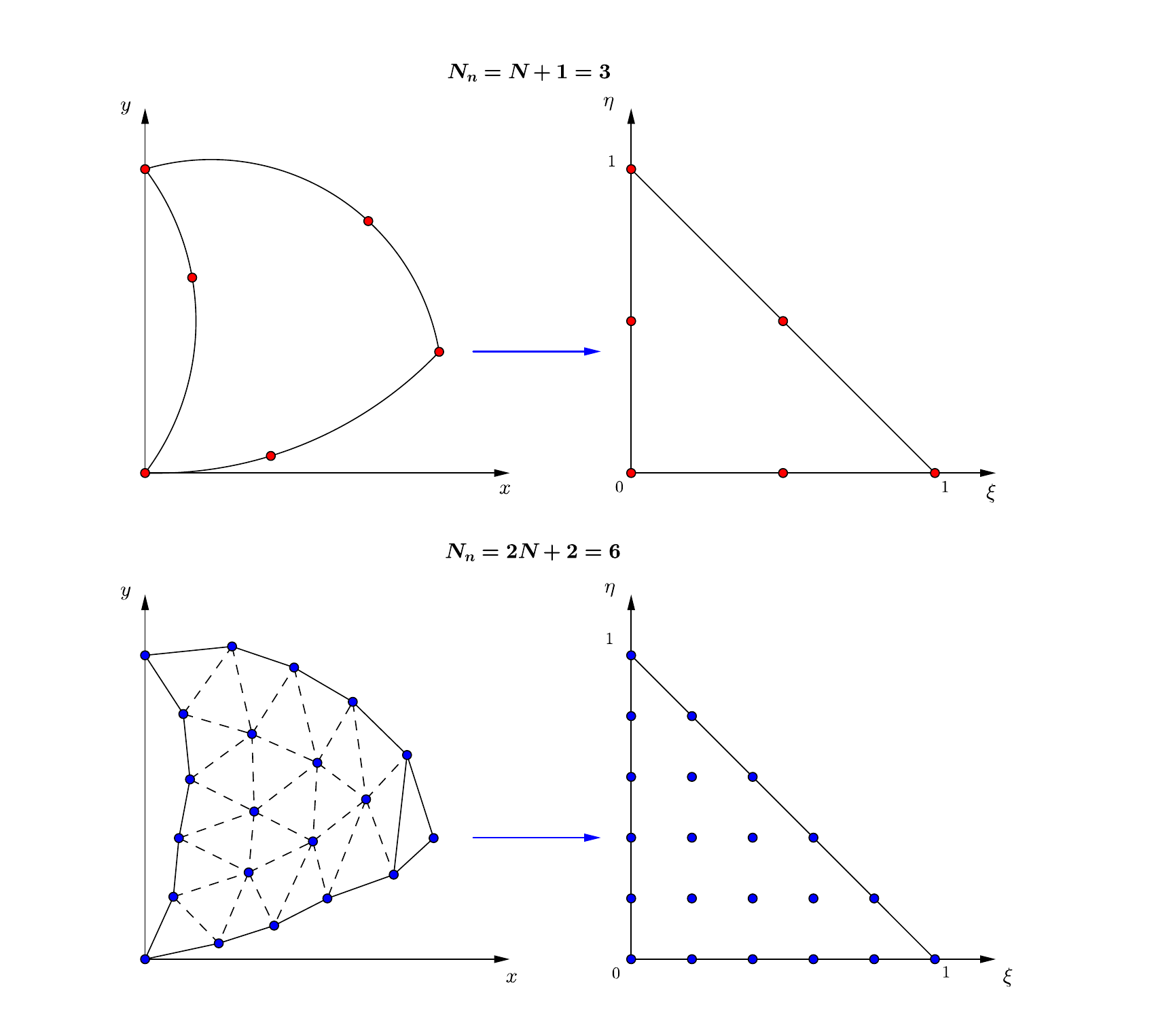}
		\end{tabular} 
    \caption{Element geometry description with the isoparametric approach (top row) and the sub-grid setting (bottom row) with the corresponding degrees of freedom highlighted in red and blue, respectively.}
    \label{fig.ISOvsSG}
	\end{center}
\end{figure}

We rely on two different settings for moving the computational mesh, namely either the \textit{isoparametric} approach and the \textit{sub-grid} approach. Equation \eqref{eqn.xiSN} provides the reference  coordinates $\kk_k=(\xi_k,\eta_k,\zeta_k)$ of the degrees of freedom $k$ needed to describe the element geometry. For the \textit{piecewise linear} subgrid approach, which is very similar to the \textit{agglomeration approach} of Bassi et al. \cite{BassiAgglomeration}, we use $N_s = 2N+1$ \textit{subelements} and thus $N_n=2N+2$ \textit{subnodes} along each element edge in order to describe the 
piecewise linear geometry on the sub-grid level. 
The basis functions $\phi_l(\x)$ that approximate the numerical solution \eqref{eqn.uh} are still defined on the reference element $T_E$ (see Figure \ref{fig.Te}), which is split according to the sub-grid definition provided in Section \ref{sec.subgrid}, and we apply a \textit{piecewise linear sub-mapping} to \textit{each} simplex sub-element $S_{ij}^n$ for the transformation from the reference coordinates in $\boldsymbol{\xi}=(\xi,\eta,\zeta)$ to the physical space in $\x=(x,y,z)$, that is
\begin{equation} 
 \mathbf{x} = \mathbf{X}^n_{1,i} + 
\left( \mathbf{X}^n_{2,i} - \mathbf{X}^n_{1,i} \right) \xi + 
\left( \mathbf{X}^n_{3,i} - \mathbf{X}^n_{1,i} \right) \eta + 
\left( \mathbf{X}^n_{4,i} - \mathbf{X}^n_{1,i} \right) \zeta,
 \label{xietaTransf} 
\end{equation} 
with $\mathbf{X}^n_{k,i} = (X^n_{k,i},Y^n_{k,i},Z^n_{k,i})$ denoting the vector of physical spatial coordinates of the $k$-th vertex of sub-element $S_{ij}^n$, according the local sub-grid connectivity given in \cite{DGLimiter3}. In the case of an \textit{isoparametric} description of curved spatial control volumes we use $N_s=N+1$ and, for the element configuration, a fully isoparametric mapping is adopted that can be retrieved by applying Eqn. \eqref{eqn.stgeom1} at the reference time $\tau=0$. Figure \ref{fig.ISOvsSG} shows the two-dimensional degrees of freedom, also called \textit{sub-nodes}, for both settings in the case of $N=2$, i.e. leading to $N_n=N+1=3$ in the isoparametric case and $N_n=2N+2=6$ in the case of piecewise linear subgrid elements.

The physical coordinates $\mathbf{X}_k^n(\kk_k):=\mathbf{X}_k^n$ of each sub-node $k$ at the current time level $t^n$ can be conveniently computed using the expansion \eqref{eqn.stgeom1} with $\tilde{\boldsymbol{\xi}}=(\kk_k,0)$. Please note that at time $t^n$ the mesh is continuous by definition, thus we can use either the local isoparametric description \eqref{eqn.stgeom1} or the piecewise linear subgrid mapping to evaluate the spatial coordinates of the degrees of freedom for the entire computational mesh, ensuring its continuity. This is the starting point for computing the corresponding  Lagrangian positions $\mathbf{X}_k^{Lag}$.

For each sub-node let $\mathcal{V}_k$ and $\mathcal{W}_k$ represent its associated main grid and sub-grid Voronoi neighborhood, respectively, composed by all corresponding neighbor cells $T_{j}$ and sub-cells $S_{j}$ that share the common sub-node $k$. Let furthermore $b_j$ denote the \textit{MOOD indicator} associated with each main cell $T_j$: it is allowed to assume only two values, either $b_j^n=0$ or $b_j^n=1$. As discussed later in Section \ref{sec.limiter}, if $b_j^n=0$ the element does not need any limiting procedure, while $b_j^n=1$ is used to mark the so-called \quotew{problematic cells} that are affected by the limiter. Therefore, we have to consider the \textit{effective neighborhood} $\mathcal{G}_k$ of sub-node $k$ which is built by adding either the main neighbor element $T_{j}$ in $\mathcal{V}_k$ if $b_j^n=0$ or the Voronoi sub-cells $S_j$ of $T_j$ in $\mathcal{W}_k$ when $b_j^n=1$, that is
\begin{equation}
\mathcal{G}_k^n = \left\{ \begin{array}{ccc}  T_{j} \in \mathcal{V}_k & \textnormal{ if } & b_j^n=0 \\
                                            S_{j} \in \mathcal{W}_k & \textnormal{ if } & b_j^n=1 
                \end{array}  \right. . 
\label{eqn.VorSN}
\end{equation}
Neighborhood $\mathcal{G}_k^n$ is composed by a total number $N_g^n$ neighbor elements $T_g$ and it is time-dependent, since the MOOD indicator $b_j^n$ may change in principle for all cells at each time step. The local velocity contribution $\mathbf{V}_{k,g}$ to sub-node $k$ from the effective neighbor $T_g$ is extracted from the corresponding state $\Q_{k,g}$, that is given by
\begin{equation}
\Q_{k,g} = \left\{ \begin{array}{ccc}  \left( \int \limits_{0}^{1} \theta_l(\kk_{m(k)}, \tau) d \tau \right) \hat \q^{n}_{l,g} & \textnormal{ if } & b_g^n=0 \\
                                       \mathbf{v}_{S_{j}}(\x,t^n) & \textnormal{ if } & b_g^n=1 
                \end{array}  \right. , 
\label{eqn.stateSN}
\end{equation}
where $m(k)$ denotes a mapping from the global sub-node number $k$ defined in $\mathcal{T}^n_{\Omega}$ to the local sub-node number in element $T_j$. In other words, we take the time integral of the high order extrapolated state at sub-node $k$ if the neighbor cell is not marked as problematic, otherwise we rely on the projector operator \eqref{eqn.vh} applied to element $T_g$ for getting a robust low order state. In this way, we guarantee the sub-node state $\Q_{k,g}$ to be always valid, i.e. physically and numerically acceptable. The corresponding velocity vector $\mathbf{V}_{k,g}$ is extracted from the state $\Q_{k,g}$ according to the governing equations \eqref{PDE}.

\begin{figure}[!htbp]
  \begin{center}
    \begin{tabular}{c} 
      \includegraphics[width=0.8\textwidth]{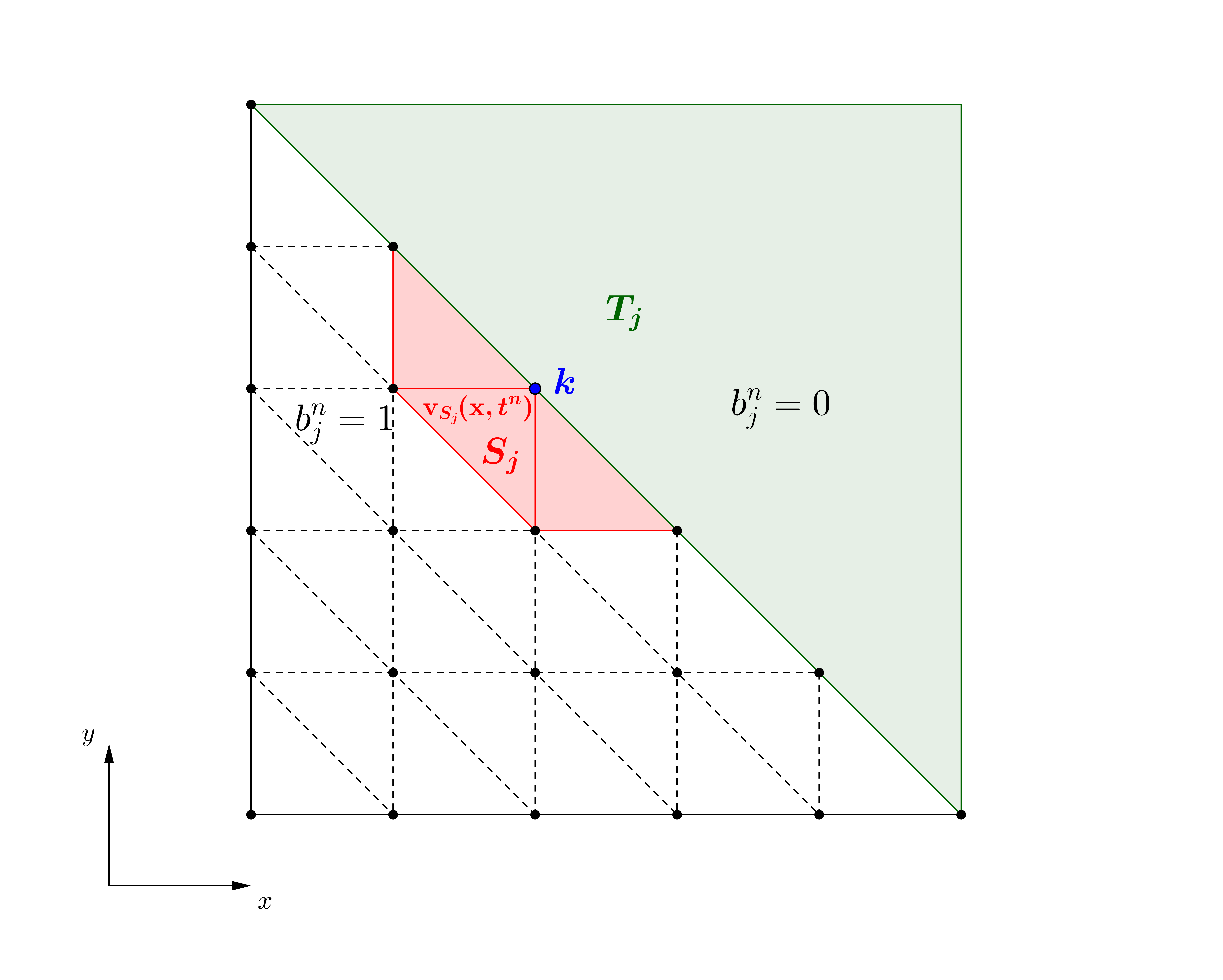}
		\end{tabular} 
    \caption{Neighborhood $\mathcal{G}_k^n$ of sub-node $k$: the left element is assigned with $b_j^n=1$, hence the states $\Q_{k,g}$ are given by the sub-cell finite volume solution $\mathbf{v}_{S_{j}}(\x,t^n)$, while the right cell does not need any limiting procedure and the corresponding velocity vector $\mathbf{V}_{k,g}$ is extracted from the high order extrapolated state at sub-node $k$.}
    \label{fig.NS_k}
	\end{center}
\end{figure}

For each sub-node a unique velocity vector $\overline{\mathbf{V}}_k$ must be computed. The sub-nodes can be classified into four different types, as depicted in Figure \ref{fig.SubGrid}:
\begin{enumerate}
	\item \textit{vertex} sub-nodes coincide with the vertexes of cell $T_i$ on the main grid. The associated velocity vector $\overline{\mathbf{V}}_k$ is extracted from the corresponding state $\Q_k$ which is simply computed as the arithmetic average among the local contributions coming from the neighborhood of sub-node $k$, hence
	\begin{equation}
    \Q_k = \frac{1}{N_g^n} \sum \limits_{T_g \in \mathcal{G}_k^n}{\Q_{k,g}};
    \label{eqn.vrtxVel}
\end{equation}
	\item \textit{edge} sub-nodes appear only for the three-dimensional case and they are aligned along each edge of the tetrahedron. For these sub-nodes the velocity is computed in the same manner used for the vertex sub-nodes, thus relying on \eqref{eqn.vrtxVel};
	\item \textit{face} sub-nodes belong to the faces of the main element. Here, we propose to evaluate the HLL state \cite{hll} at sub-node $k$ in order to obtain the associated velocity, since in the neighborhood $\mathcal{G}_k^n$ only two main elements are involved, i.e. the right $T_R$ and the left $T_L$ neighbor. Let $\mathbf{n}^n=(n_x,n_y,n_z)^n$ denote the outward pointing unit normal vector of the face of cell $T_L$ where sub-node $k$ is lying and let $\Q_{k,L}$ and $\Q_{k,R}$ be the left and right state computed with \eqref{eqn.stateSN}. The sub-node HLL state $\Q_k:=\Q_{k,HLL}$ is evaluated according to \cite{ToroBook} as
	\begin{equation}
	 \Q_{k,HLL} = \frac{s_R \Q_{k,R} - s_L \Q_{k,L} + \left( \F(\Q_{k,L}) - \F(\Q_{k,R}) \right) \cdot \mathbf{n}^n}{s_R-s_L}. 
	\label{eqn.hllVel}
	\end{equation}
		The signal speeds $s_L$ and $s_R$ are defined as usual,  
		\begin{equation}
		s_L = \min \left(0,\Lambda_{L},\Lambda_{R}\right) \qquad s_R = \max \left(0,\Lambda_{L},\Lambda_{R}\right),
		\label{eqn.sHLL}
		\end{equation}
		with $\Lambda_L$ and $\Lambda_R$ denoting the diagonal matrix of eigenvalues of the Jacobian matrix of the flux in normal direction $\mathbf{A}=\partial \mathbf{F} / \mathbf{Q} \cdot \mathbf{n}$, computed from the corresponding left and right states, respectively;
\item \textit{internal} sub-nodes are located in the inner part of $T_i$ and their velocity vectors are determined by solving a \textit{local Laplace equation} within each element, that is
	\begin{equation}
	 \Delta \mathbf{V} = \mathbf{0},
	\label{eqn.Laplace}
	\end{equation}
	with Dirichlet-type boundary conditions given by the velocities previously computed for vertex, edge and face sub-nodes. Equation \eqref{eqn.Laplace} is solved by a classical second order finite element method on the local sub-grid level.
\end{enumerate}

Once the sub-node velocity vector $\overline{\mathbf{V}}_k$ is known, the evaluation of the Lagrangian position is straightforward and reads
\begin{equation}
 \mathbf{X}_k^{Lag} = \mathbf{X}_k^{n} + \Delta t \cdot \overline{\mathbf{V}}_k \qquad \forall k \in \mathcal{T}_{\Omega}^{n}.
\label{eqn.LagX}
\end{equation}

\subsubsection{The rezoning step}
After the Lagrangian phase, coordinates $\mathbf{X}_k^{Lag}$ might yield a complex mesh configuration, with highly compressed or twisted control volumes which could degenerate even to tangled elements. As a consequence, the time step would become very small according to \eqref{eq:timestep}, or the computation would blow up due to the presence of invalid cells, i.e. computational elements with negative volume. To improve the mesh quality, a so-called rezoning strategy is usually applied \cite{KnuppRezoning,MaireRezoning} in order to improve the mesh quality. Rezoning algorithms do not take into account any physical aspect, but they are based and developed on geometrical considerations. Here we use the \textit{same} strategy described in \cite{LagrangeMHD,Lagrange3D} for triangular and tetrahedral elements: it consists in optimizing a goal function which is defined locally for each control volume. The crucial point in our approach is that the entire rezoning procedure is carried out \textit{on the sub-grid level}, which is composed by simplex elements defined by straight boundaries that perfectly match the requirements needed to perform the rezoning algorithm detailed in \cite{LagrangeMHD,Lagrange3D}. Therefore, once the goal function has been optimized, the rezoned coordinates $\mathbf{X}_k^{Rez}$ are available for each sub-node of the computational mesh.

\subsubsection{The relaxation step}
The final mesh configuration $\mathcal{T}_{\Omega}^{n+1}$ is then given by a weighted linear combination between the Lagrangian coordinates $\mathbf{X}_k^{Lag}$ and the rezoned position vectors $\mathbf{X}_k^{Rez}$ of each sub-node $k$, hence
\begin{equation}
\mathbf{X}_k^{n+1} = \mathbf{X}_k^{Lag} + \omega_k \left( \mathbf{X}_k^{Rez} - \mathbf{X}_k^{Lag} \right),
\label{eqn.relaxation}
\end{equation}
with $\omega_k$ representing a sub-node coefficient bounded in the interval $[0,1]$. According to \cite{MaireRezoning}, $\omega_k$ is associated to the deformation of the Lagrangian grid over the time step $\Delta t$, that is
\begin{equation}
\mathbf{F} = \frac{\partial \mathbf{X}^{Lag}}{\partial \mathbf{X}^{n}},
\label{eqn.F}
\end{equation}   
where $\mathbf{F}$ denotes the deformation gradient tensor. Since the sub-elements, on which the mesh motion procedure is carried out, are simplex elements, one can rely either on the original technique given in \cite{MaireRezoning} or on the variant recently proposed in \cite{LagrangeISO} to compute the tensor $F$ and subsequently to extract the blending factor $\omega_k$. All the details can be found in the aforementioned references.

After completion of the mesh motion algorithm, the mesh configuration at the new time level $t^{n+1}$ is \textit{known} and \textit{continuous}, hence allowing the space-time control volumes to be uniquely defined within the time step $\Delta t = t^{n+1} - t^{n}$. To maintain algorithm simplicity, the old mesh configuration is connected to the new mesh configuration by \textit{straight lines}, thus obtaining a linear description in the time evolution of the control volumes. The mapping in time is linear and simply reads
\begin{equation}
t = t^n + \tau \Delta t.
\label{eqn.tau}
\end{equation}
Now, a direct high order Arbitrary-Lagrangian-Eulerian DG scheme can be applied to solve the governing equations \eqref{PDE}.

\subsection{ADER-DG scheme on moving unstructured meshes}
\label{sec.SolAlg}

\begin{figure}[!htbp]
  \begin{center}
    \begin{tabular}{c} 
      \includegraphics[width=0.7\textwidth]{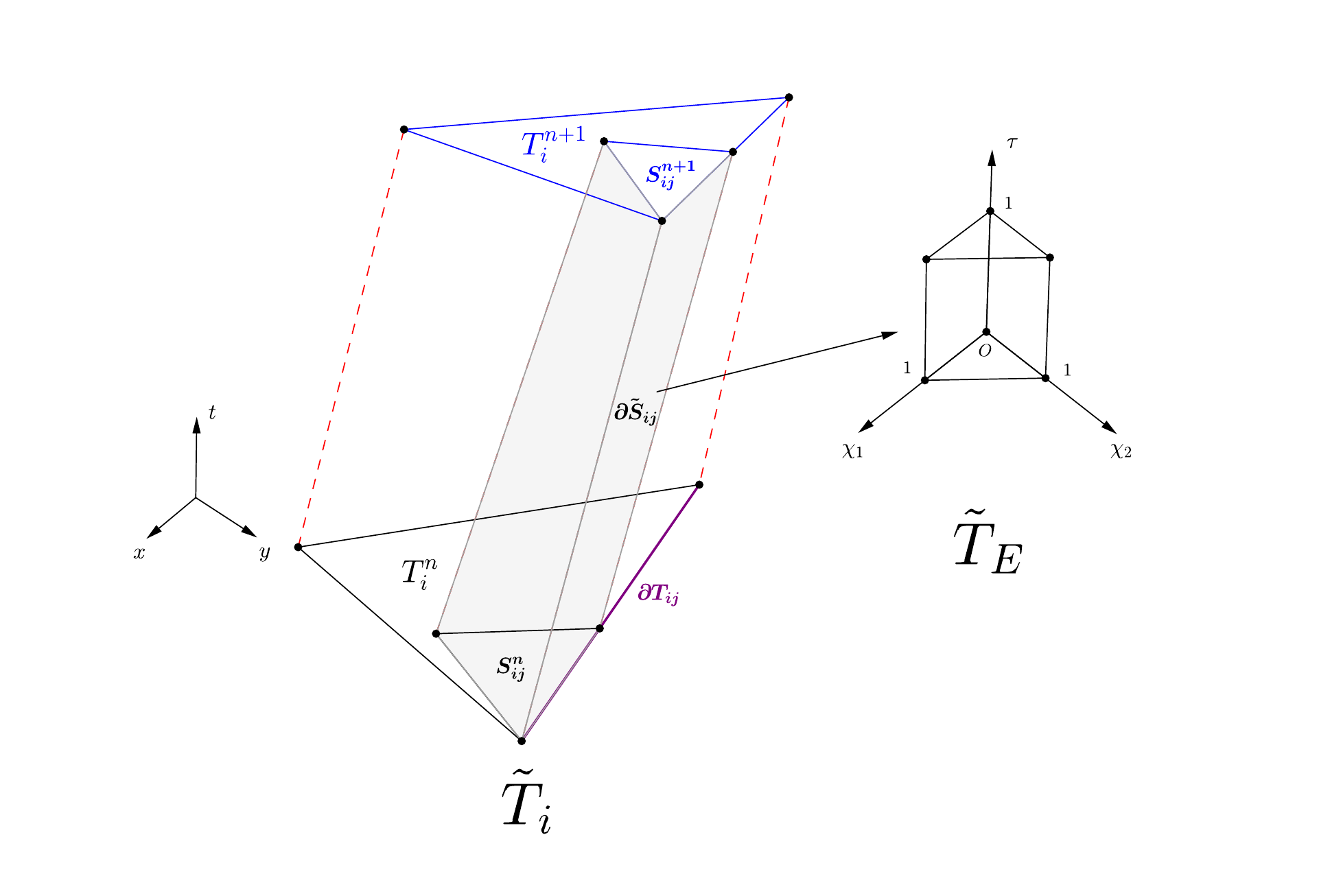}
		\end{tabular} 
    \caption{Space-time evolution of element $T_i^n$ within one time step $\Delta t$. The physical space-time sub-volumes $\tilde{S}_{ij}$ (left) are mapped onto a reference element (right) defined in $\tilde{\boldsymbol{\chi}}$ at the aid of a set of linear basis functions.}
    \label{fig.DGscheme}
	\end{center}
\end{figure}

A fully discrete one-step ADER-DG scheme is derived starting from the predictor solution $\q_h(\x,t)$, available from the local predictor strategy described in Section \ref{sec.lst}, and the space-time control volumes $\tilde{T}_i = T_i(t) \times \left[t^n,t^{n+1}\right]$, which ensure a continuous mesh configuration in space and time thanks to the mesh motion procedure illustrated in Section \ref{sec.meshMot}. The PDE system \eqref{PDE} is written in a more compact space-time divergence form as
\begin{equation}
\tilde \nabla \cdot \tilde{\F} = \mathbf{0} 
  \qquad \tilde \nabla  = \left( \frac{\partial}{\partial x}, \, \frac{\partial}{\partial y}, \, \frac{\partial}{\partial z}, \, \frac{\partial}{\partial t} \right)^T,
\label{eqn.st.pde}
\end{equation}
with  $\tilde \nabla$ representing a space-time divergence operator and $\tilde{\F} = \left( \mathbf{f}, \, \mathbf{g}, \, \mathbf{h}, \, \Q \right)$ the corresponding space-time flux tensor. Multiplication of \eqref{eqn.st.pde} by test functions $\phi_k$, which are taken to be identical with the spatial basis functions used in \eqref{eqn.uh}, and subsequent integration over the four-dimensional space-time control volume $\tilde{T}_i$ yields
\begin{equation}
\int \limits_{\tilde{T}_i} \phi_k \tilde \nabla \cdot \tilde{\F} \, d\mathbf{x} dt = \mathbf{0}. 
\label{eqn.intPDE}
\end{equation}
Application of Gauss' theorem allows the above expression to be reformulated as
\begin{equation}
\int \limits_{\partial \tilde{T}_i} \phi_k  \tilde{\F} \cdot \ \mathbf{\tilde n} \, dS dt - \int \limits_{\tilde{T}_i} \tilde \nabla \phi_k \cdot \tilde{\F} \, d\mathbf{x} dt = \mathbf{0}, 
\label{eqn.GaussPDE}
\end{equation}
where $\mathbf{\tilde n} = (\tilde n_x,\tilde n_y,\tilde n_z,\tilde n_t)$ is the outward pointing space-time unit normal vector on the space-time face $\partial \tilde{T}_i$, that is given by the evolution of each face of element ${T}_i$ within the timestep $\Delta t$. Specifically, a total number of five or six space-time faces are needed to bound the space-time volume $\tilde{T}_i$ for $d=2$ or $d=3$, respectively. Figure \ref{fig.DGscheme} shows the two-dimensional case: the lateral space-time faces of element $T_i$ involve the \textit{Neumann neighborhood} $\mathcal{N}_i$, which is the set of directly adjacent neighbors $T_j$ that share a common face $\partial {T}_{ij}$ with element $T_i$, then volume $\tilde{T}_i$ is closed by the cell configuration at the old and at the new time level, that is
\begin{equation}
\partial \tilde{T}_i = \left( \bigcup \limits_{\tilde{T}_j \in \mathcal{N}_i} \partial \tilde{T}_{ij} \right) 
\,\, \cup \,\, T_i^{n} \,\, \cup \,\, T_i^{n+1}.
\label{dCi}
\end{equation}  

The boundary integral appearing in \eqref{eqn.GaussPDE} is replaced by a numerical flux function that provides the coupling between neighbor elements, which was not considered in the predictor step presented in Section \ref{sec.lst}. The numerical flux, also known as Riemann solver, is written in space-time normal direction as $\mathcal{G}((\q_h^-,\nabla \q_h^-),(\q_h^+,\nabla \q_h^+)) \cdot \mathbf{\tilde n}$ and it involves the left $(\q_h^-,\nabla \q_h^-)$ and right $(\q_h^+,\nabla \q_h^+)$ high order boundary-extrapolated data and gradients. Using approximation \eqref{eqn.uh} and the predictor solution $\q_h$, the arbitrary high order direct ALE one-step ADER-DG scheme reads
\begin{equation}
\left( \int \limits_{T_i^{n+1}} \phi_k \phi_l d\x \right) \hat{\mathbf{u}}^{n+1}_{l} = \left( \int \limits_{T_i^{n}} \phi_k \phi_l d\x \right) \hat{\mathbf{u}}^{n}_{l} - \int \limits_{\partial \tilde{T}_i} \phi_k \mathcal{G}\left((\q_h^-,\nabla \q_h^-),(\q_h^+,\nabla \q_h^+) \right) \cdot \mathbf{\tilde n} \, dS dt  
 + \int \limits_{\tilde{T}_i} \tilde \nabla \phi_k \cdot \tilde{\F}(\q_h,\nabla \q_h) \, d\mathbf{x} dt. 
\label{eqn.ALE-DG}
\end{equation}
In this work we rely on a simple and very robust Rusanov flux \cite{Rusanov:1961a} to evaluate the term $\mathcal{G}$. It has already been applied to the ALE context \cite{Lagrange2D,Lagrange3D} and, following \cite{DumbserNSE}, it includes both the convective and the viscous terms, hence
\begin{equation}
  \mathcal{G}\left((\q_h^-,\nabla \q_h^-),(\q_h^+,\nabla \q_h^+) \right) \cdot \mathbf{\tilde n}  =  
  \frac{1}{2} \left( \tilde{\F}(\q_h^+,\nabla \q_h^+) + \tilde{\F}(\q_h^-,\nabla \q_h^-)  \right) \cdot \mathbf{\tilde n}  - 
  \frac{1}{2} \left(|s_{\max}| + 2\eta |s_{\max}^{\nu}| \right) \left( \q_h^+ - \q_h^- \right).  
  \label{eqn.rusanov} 
\end{equation} 
Here, $|s_{\max}|$ represents the maximum eigenvalue of the ALE Jacobian matrix in space normal direction, which is
\begin{equation} 
\mathbf{A}^{\!\! \mathbf{V}}_{\mathbf{n}}(\Q,\nabla \Q):=\left(\sqrt{\tilde n_x^2 + \tilde n_y^2 + \tilde n_z^2 }\right)\left[\mathbf{A} \cdot \mathbf{n}  - 
(\mathbf{V} \cdot \mathbf{n}) \,  \mathbf{I}\right], \qquad    
\mathbf{n} = \frac{(\tilde n_x, \tilde n_y, \tilde n_z)^T}{\sqrt{\tilde n_x^2 + \tilde n_y^2 + \tilde n_z^2 }},  
\label{eqn.An}
\end{equation} 
where $\mathbf{I}$ is the identity matrix, $\mathbf{V} \cdot \mathbf{n}$ denotes the local normal mesh velocity and $\mathbf{A}=\frac{\partial \mathbf{F}(\Q,\nabla \Q)}{\partial \Q}$. Then, $|s_{\max}^{\nu}|$ is the maximum eigenvalue of the Jacobian matrix of the \textit{viscous} operator given by $\mathbf{D}=\frac{\partial \mathbf{F}(\Q,\nabla \Q)}{\partial (\nabla \Q \cdot \mathbf{n})} \cdot \mathbf{n}$.
Finally, the factor $\eta$ is estimated according to \cite{MunzDiffusionFlux,DumbserNSE} from the solution of the generalized diffusive Riemann problem as
\begin{equation}
\eta = \frac{2N+1}{h_\nu},
\label{eqn.eta}
\end{equation}
where the characteristic size $h_\nu$ is given by the sum of the distances between the barycenter of the adjacent elements ($T_i$ and $T_j$)  and the barycenter of the face $\partial T_{ij}$ along which the numerical flux is computed. 

At this point we have two different approaches to carry on with the one-step ALE ADER-DG scheme \eqref{eqn.ALE-DG}, depending on the strategy adopted for the mesh motion, namely the \textit{piecewise linear} or the \textit{isoparametric} geometry approximation. \\ 
In the first case, the mesh motion procedure yields an element configuration in which the cell is bounded by a set of linear space-time surfaces, hence leading to a general polyhedral element. As a consequence, the space-time volume $\tilde{T}_i$ is \textit{decomposed} into a set of $\mathcal{S}$ corresponding space-time sub-grid volumes $\tilde{S}_{ij}$, as shown in Figure \ref{fig.DGscheme}. Each space-time sub-volume $\tilde{S}_{ij}$ is parametrized using a set of \textit{linear} basis functions $\alpha$ defined on a local reference system $\tilde{\boldsymbol{\chi}}=(\chi_1,\chi_2,\chi_3,\tau)$, in which the reference time coordinate $\tau$ is orthogonal to the reference space coordinates that lie on element $T_i^n$. Such parametrization reads
\begin{equation}
\tilde{S}_{ij} = \sum\limits_{a=1}^{N_{\alpha}}{\alpha_a(\tilde{\xxi}_a) \, \mathbf{\tilde{X}}_{ij,a} },
\label{eqn.alpha}
\end{equation} 
where the degrees of freedom $\mathbf{\tilde{X}}_{ij,a}$ are known and are given by the coordinates of the sub-cell vertexes at time $t^n$ and $t^{n+1}$. As fully detailed in \cite{Lagrange2D,Lagrange3D}, $N_{\alpha}=2(d+1)$ and the positions $\mathbf{\tilde{X}}_{ij,a}$ are directly available from the corresponding sub-node $k$ with $\mathbf{X}_k^{n}$ and $\mathbf{X}_k^{n+1}$. Index $k$ is obtained relying on the local sub-grid connectivity of sub-cell $S_{ij}$, see \cite{DGLimiter3}. Thus, the direct ALE ADER-DG scheme \eqref{eqn.ALE-DG} with piecewise linear sub-cell representation of the geometry looks  very similar to the corresponding finite volume scheme presented in \cite{Lagrange2D,Lagrange3D} and it can be formulated as 
\begin{eqnarray}
\left( \sum \limits_{s=1}^{\mathcal{S}} \int \limits_{S_{ij}^{n+1}} \phi_k \phi_l d\x \right) \hat{\mathbf{u}}^{n+1}_{l} &=& \left( \sum \limits_{s=1}^{\mathcal{S}} \int \limits_{S_{ij}^{n}} \phi_k \phi_l d\x \right) \hat{\mathbf{u}}^{n}_{l} + \sum \limits_{s=1}^{\mathcal{S}} \int \limits_{\tilde{S}_{ij}} \tilde \nabla \phi_k \cdot \tilde{\F}(\q_h,\nabla \q_h) \, d\mathbf{x} dt \nonumber \\
&& - \sum \limits_{s=1}^{\mathcal{S}} \int \limits_{\partial \tilde{S}_{ij}} \phi_k \mathcal{G}\left((\q_h^-,\nabla \q_h^-),(\q_h^+,\nabla \q_h^+)\right) \cdot \mathbf{\tilde n} \, dS dt.
\label{eqn.ALE-DG-subgrid}
\end{eqnarray}
Expression \eqref{eqn.alpha} allows the evaluation of the space-time normal vectors $\mathbf{\tilde n}$ as well as the Jacobian of the transformation, as done in \cite{Lagrange2D,Lagrange3D}. 

If the \textit{isoparametric} approach is employed to approximate the element geometry, the whole space-time volume $\tilde{T}_i$ is parametrized to the space-time reference element $\tilde{T}_E$ using a set of \textit{high order} basis functions $\beta$, therefore
\begin{equation}
\tilde{T}_{i} = \sum\limits_{b=1}^{N_{\beta}}{\beta_b(\tilde{\xxi}_b) \, \mathbf{\tilde{X}}_{ij,b} }.
\label{eqn.beta}
\end{equation} 
In this case, the total number of degrees of freedom $N_{\beta}$ is given by \eqref{eqn.nDOF} with $d+1$ dimensions and for both approximations \eqref{eqn.alpha} and \eqref{eqn.beta} the basis functions $\alpha_a$ and $\beta_b$ are defined by the Lagrange interpolation polynomials passing through the space-time nodes $\tilde{\xxi}_a$ and $\tilde{\xxi}_b$, respectively. The degrees of freedom $\mathbf{\tilde{X}}_{ij,b}$ are known: those ones defined at time $t^n$ are given by the current mesh configuration, the ones at time $t^{n+1}$ are provided by the mesh motion algorithm, then the time linear mapping \eqref{eqn.tau} gives the element configuration at all intermediate time levels needed for evaluating the missing degrees of freedom of order $N+1$. Here, we are dealing with \textit{curvilinear} elements that are approximated by a set of high order basis functions and the integrals appearing in \eqref{eqn.ALE-DG} are evaluated following the approach recently proposed in \cite{LagrangeISO}, where high order finite volume schemes have been applied to curvilinear simplex elements.

For the sake of clarity all the integrals which are present in \eqref{eqn.ALE-DG} are computed on the space-time reference element $T_E \times [0,1]$ employing Gaussian quadrature rules of sufficient precision, see \cite{stroud} for details.\\
Finally, we point out that even for the direct ALE ADER-DG algorithm presented in this paper, the scheme provided by \eqref{eqn.ALE-DG-subgrid} automatically satisfies the geometric conservation law (GCL) for \textit{all} test functions $\phi_k$, since according to Gauss' theorem it follows
\begin{equation}
\int \limits_{\partial \tilde{T}_i} \phi_k \mathbf{\tilde n} \, dS dt - \int \limits_{\tilde{T}_i} \tilde \nabla \phi_k \, d\mathbf{x} dt = 0.
\label{eqn.GCL}
\end{equation}

\subsection{\textit{A posteriori} sub-cell finite volume limiter on moving unstructured meshes}
\label{sec.limiter}
The numerical scheme presented in the previous section needs a nonlinear limiting procedure to avoid the Gibbs phenomenon at shock waves or other discontinuities, which typically occur while solving nonlinear hyperbolic systems of the form \eqref{PDE}. In our approach, we rely on the very recent technique developed in \cite{DGLimiter1,DGLimiter2,DGLimiter3} based on the MOOD paradigm \cite{CDL1,CDL2,CDL3}, where an \textit{a posteriori} limiter is applied in order to stabilize the numerical solution. All the details can be found in the aforementioned references, hence we limit us to briefly recall the main features of the limiter that makes use of a robust second order TVD finite volume scheme on the sub-grid level.

The unlimited ALE ADER-DG scheme \eqref{eqn.ALE-DG} generates a so-called \textit{candidate solution} $\mathbf{u}^*_{h}(\x,t^{n+1})$, which is checked against a set of \textit{detection criteria} that must be fulfilled in order to accept the discrete solution at the new time level. If the candidate solution does not satisfy all the requirements, the numerical solution is \textit{locally} recomputed using a second order direct ALE ADER finite volume scheme, based on a TVD reconstruction with Barth \& Jespersen slope limiter, as done in \cite{ALEMOOD1,ALEMOOD2}. \\
The \textit{a posteriori} sub-cell limiter procedure can be summarized as follows:
\begin{itemize}
	\item [\tiny$\bullet$] compute the candidate solution $\mathbf{u}^*_{h}(\x,t^{n+1})$ for each cell $T_i^{n}$ by means of \eqref{eqn.ALE-DG};
	\item [\tiny$\bullet$] use the projection operator \eqref{eqn.vh} to obtain the candidate solution $\mathbf{v}^*_h(\x,t^{n+1})$ on the sub-grid level for each sub-cell $S_{ij}$ of element $T_i^{n}$;
	\item [\tiny$\bullet$] check the candidate solution $\mathbf{v}^*_h(\x,t^{n+1})$ against the detection criteria: according to \cite{DGLimiter3}, the first criterion is given by requiring physical positivity for some quantities related to the governing system \eqref{PDE}, such as density and pressure, if the compressible Euler equations for gas dynamics are considered. Then, a relaxed discrete maximum principle (RDMP) is applied in the sense of polynomials, hence verifying
	\begin{equation}
	\min_{m \in \mathcal{V}_i} (\mathbf{v}_h(\x_m,t^{n})) - \delta \leq \mathbf{v}^*_h(\x,t^{n+1}) \leq \max_{m \in \mathcal{V}_i} (\mathbf{v}_h(\x_m,t^{n})) + \delta \qquad \forall \x \in T_i^{n},
	\label{eqn.RDMP}
	\end{equation}
	where $\mathcal{V}_i$ represents the Voronoi neighborhood of cell $T_i^{n}$ and $\delta$ is a parameter which, according to \cite{DGLimiter1,DGLimiter2}, reads
	\begin{equation}
	\delta = \max \left[ \delta_0, \epsilon \cdot \left( \max_{m \in \mathcal{V}_i} (\mathbf{v}_h(\x_m,t^{n})) - \min_{m \in \mathcal{V}_i} (\mathbf{v}_h(\x_m,t^{n})) \right) \right],
	\label{eqn.deltaRDMP}
	\end{equation}
	with $\delta_0=10^{-4}$ and $\epsilon=10^{-3}$. If a cell passes the detection criteria in \textit{all} its sub-cells, then the cell is marked as \quotew{good} using the MOOD indicator $b_i^n=0$, otherwise the cell is \quotew{problematic} or troubled with $b_i^n=1$. Such indicator is also employed for determining the local velocity contribution \eqref{eqn.VorSN} in the Lagrangian phase of the mesh motion procedure;
	\item [\tiny$\bullet$] at this point the numerical solution at the new time level $\mathbf{u}_{h}(\x,t^{n+1})$ must be determined: if $b_i^n=0$ then we simply have $\mathbf{u}_{h}(\x,t^{n+1}) = \mathbf{u}^*_{h}(\x,t^{n+1})$. In the case of a troubled cell, i.e. $b_i^n=1$, the new numerical solution is first computed on the sub-grid level for each sub-cell $S_{ij}$, hence obtaining $\mathbf{v}_h(\x,t^{n+1})$. To this purpose, we propose to use a second order direct ALE ADER finite volume scheme which exactly follows the algorithm fully detailed in \cite{Lagrange2D,Lagrange3D}: the only difference with the aforementioned references is that here the finite volume scheme is applied to each sub-cell and we do not use a WENO reconstruction, but a simple and robust TVD reconstruction with Barth \& Jespersen slope limiter \cite{BarthJespersen}. The piecewise polynomial solution of the DG scheme is now recovered from the robust and stable solution on the sub-grid level by applying the reconstruction operator \eqref{eqn.intRec}, thus $\mathbf{u}_{h}(\x,t^{n+1})=\mathcal{R}(\mathbf{v}_h(\x,t^{n+1}))$. 
\end{itemize} 

\paragraph*{Remark} In order to be strictly \textit{conservative}, in a good cell $T_j^n$ with $b_j^n=0$ which is a neighbor of a troubled cell $T_i^n$ with $b_i^n=1$, the numerical solution $\mathbf{u}_{h}(\x_j,t^{n+1})$ is also \textit{recomputed}.  Indeed, the numerical flux on the common boundary face $\partial \tilde{T}_{ij}$, shared by elements $T_j^n$ and $T_i^n$, has been evaluated on the sub-grid level with the second order TVD direct ALE finite volume scheme and it must be taken into account also in cell $T_j^n$.

\section{Test problems}
\label{sec.validation} 
\vspace{-2pt}

In this paper we focus on the \textit{Euler equations of compressible gas dynamics}, which can be cast into form \eqref{PDE} with
\begin{equation}
\label{eulerTerms}
\Q=\left( \begin{array}{c} \rho \\ \rho \mathbf{v} \\ \rho E \end{array} \right), \qquad \F(\Q) = 
\left( \begin{array}{c}  \rho \mathbf{v} \\ \rho \left(\mathbf{v}\otimes \mathbf{v}\right) + p \mathbf{I} \\ \mathbf{v} (\rho E + p)  \end{array} \right), 
\end{equation}
where $\rho$ and $p$ are the fluid density and pressure, respectively, while $\mathbf{v}=(u,v,w)$ denotes the velocity vector and $E$ represents the total energy density. The $d \times d$ identity matrix is addressed with $\mathbf{I}$ and the system is closed using the equation of state of a perfect gas with adiabatic index $\gamma$, i.e.
\begin{equation}
\label{eqn.eos} 
p = (\gamma-1)\left(\rho E - \frac{1}{2} \rho \mathbf{v}^2 \right).  
\end{equation}

If viscous flows with heat conduction are considered, the flux term in \eqref{eulerTerms} becomes
\begin{equation}
\F(\Q,\nabla \Q) = 
\left( \begin{array}{c}  \rho \mathbf{v} \\ \rho \left(\mathbf{v}\otimes \mathbf{v}\right) + \boldsymbol{\sigma}(\Q,\nabla \Q) \\ \mathbf{v} \cdot (\rho E \mathbf{I} + \boldsymbol{\sigma}(\Q,\nabla \Q)) - \kappa \nabla T  \end{array} \right),
\label{eqn.NSflux}
\end{equation}
hence obtaining the \textit{compressible Navier-Stokes equations}. The stress tensor $\boldsymbol{\sigma}(\Q,\nabla \Q)$ is computed under Stokes' hypothesis as
\begin{equation}
\boldsymbol{\sigma}(\Q,\nabla \Q) = \left(p + \frac{2}{3} \mu \nabla \cdot \mathbf{v} \right) \mathbf{I} - \mu \left( \nabla \mathbf{v} + \nabla \mathbf{v}^T \right),
\label{eqn.sigma}
\end{equation}
with $\mu$ denoting the viscosity. $T$ represents the temperature and the heat conduction coefficient $\kappa$ is linked to the viscosity through the Prandtl number $Pr$, thus
\begin{equation}
\kappa = \frac{\mu \gamma c_v}{Pr}.
\label{eqn.kappa}
\end{equation}
The specific heat at constant volume is given by $c_v=R/(\gamma-1)$ with $R$ being the gas constant which is assumed to be $R=1$, if not stated differently. In the case of viscous phenomena, the time step restriction is more severe and equation \eqref{eq:timestep} is modified into
\begin{equation}
\Delta t < \frac{\textnormal{CFL}}{2N+1} \, \min \limits_{T_i^n} \frac{h_i}{|\lambda_{\max,i}| + 2|\lambda_{\max,i}^\nu|\frac{2N+1}{h_i}}, \qquad \forall T_i^n \in \Omega^n,
\label{eq:timestep2}
\end{equation}
where, according to \cite{DumbserNSE}, the maximum viscous eigenvalue is $|\lambda_{\max,i}^\nu|=\max \left(\frac{4}{3}\frac{\mu}{\rho},\frac{\gamma \mu}{Pr \rho}\right)$.

In our ALE framework we choose to set the local mesh velocity equal to the local fluid velocity for each of the test cases shown in this paper, hence   
\begin{equation}
 \mathbf{V} = \mathbf{v}. 
\end{equation} 
Furthermore, we employ by default the piecewise linear mesh motion and the corresponding DG scheme, therefore the usage of the isoparametric approach will be explicitly declared when adopted. Finally, the initial condition might also be given in primitive variables $\U=(\rho,u,v,w,p)$.

\subsection{Numerical convergence studies}
\label{sec.conv}
\begin{figure}[!htbp]
\begin{center}
\begin{tabular}{ccc} 
\includegraphics[width=0.47\textwidth]{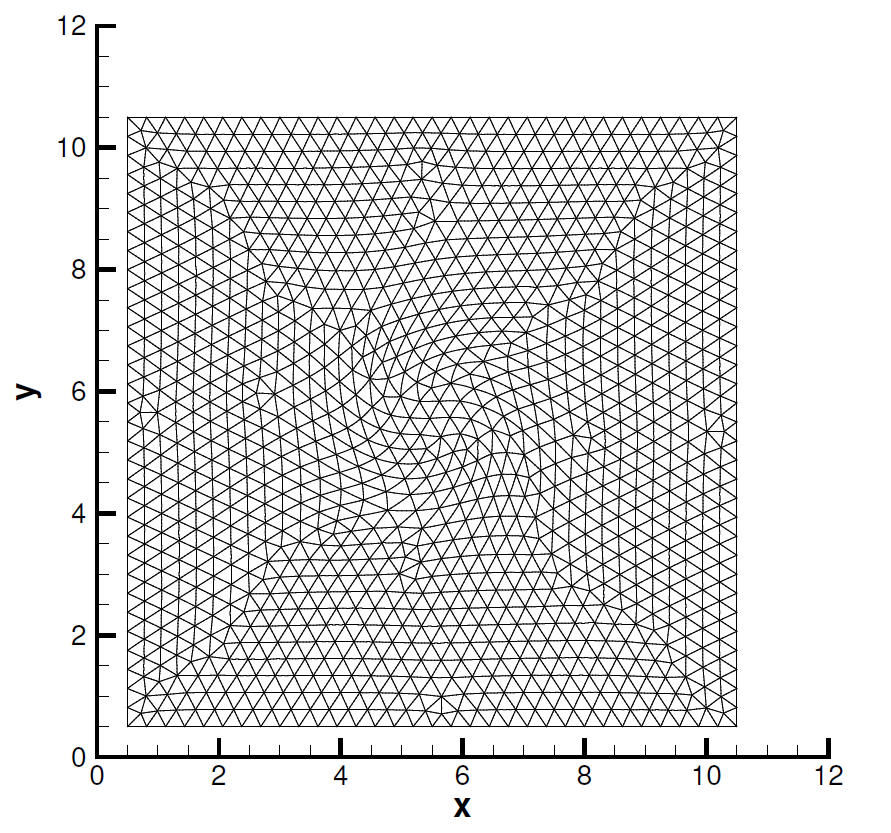}  &           
\includegraphics[width=0.47\textwidth]{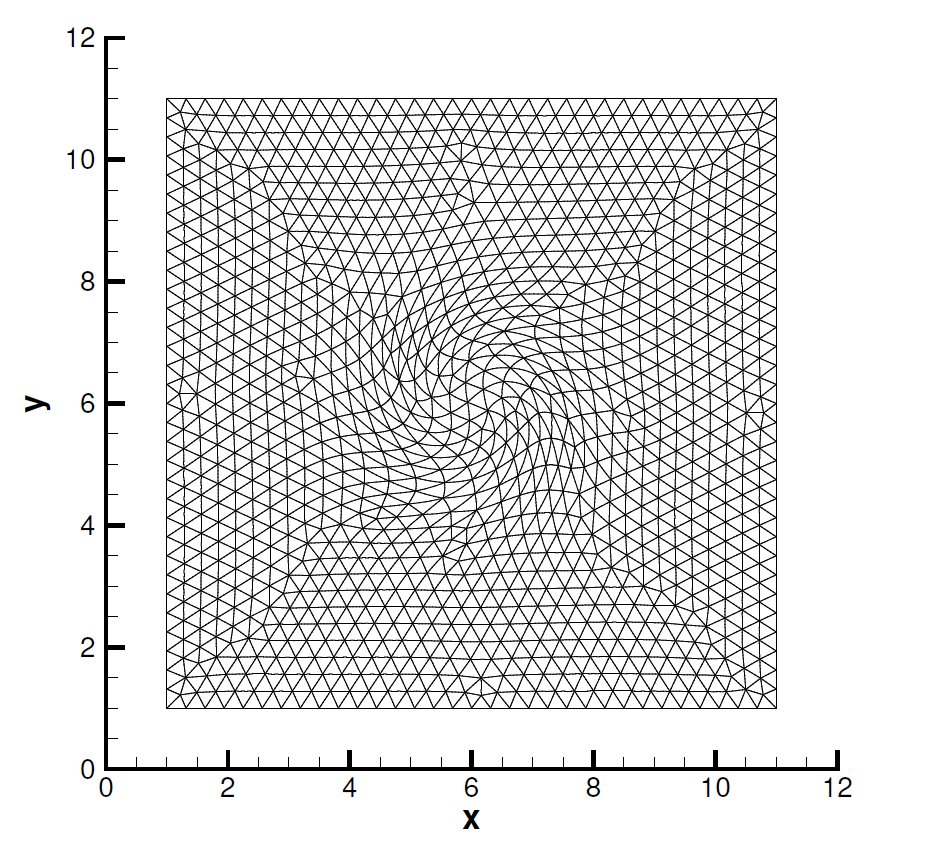} \\    
\includegraphics[width=0.47\textwidth]{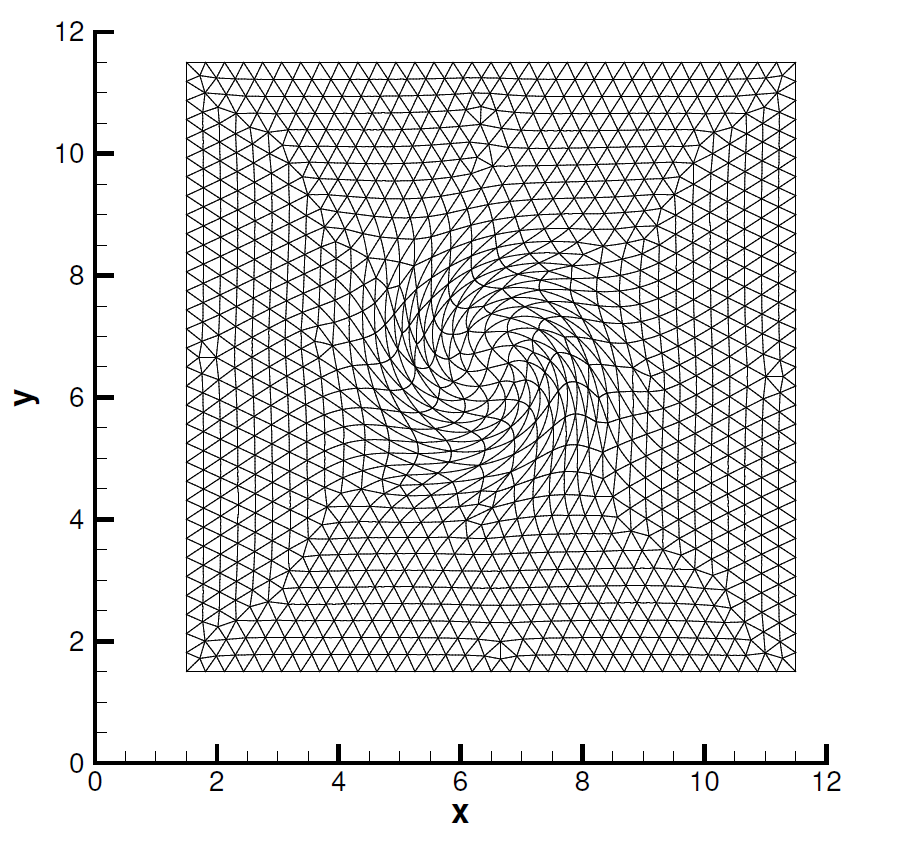}  &           
\includegraphics[width=0.47\textwidth]{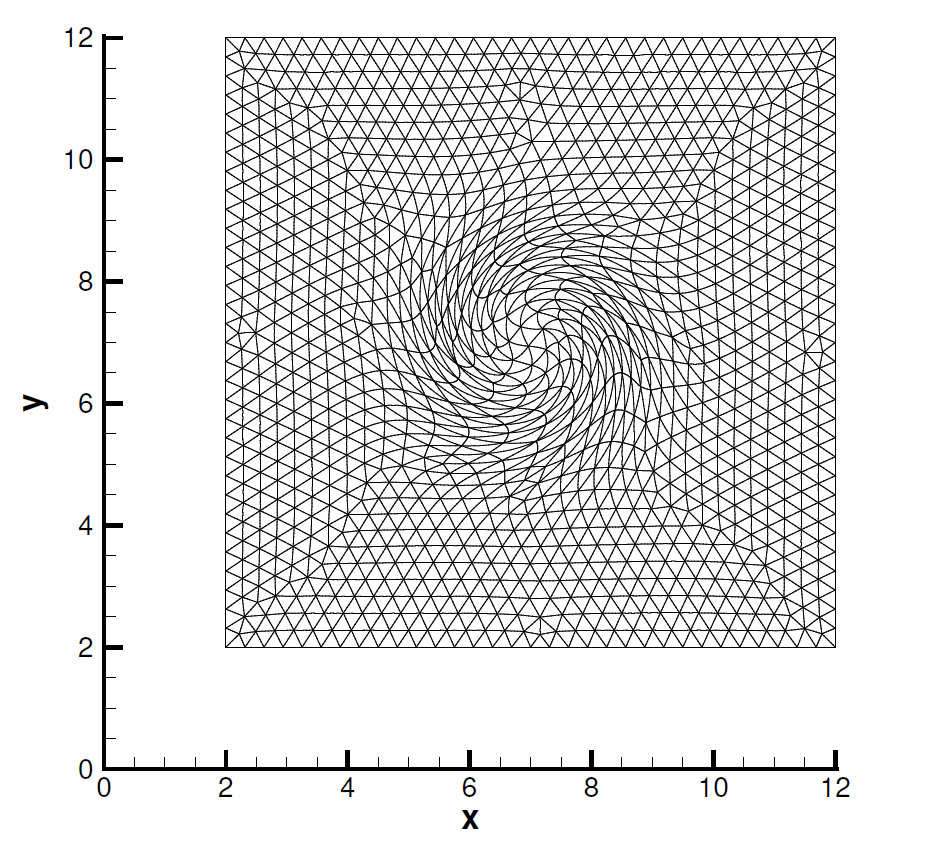} \\       
\end{tabular} 
\caption{Two-dimensional isentropic vortex test problem. Mesh configuration at output times $t=0.5$, $t=1.0$, $t=1.5$ and $t=2.0$ from top left to bottom right. A fourth order of accuracy is used to approximate the element geometry.} 
\label{fig.SV-mesh}
\end{center}
\end{figure}

The numerical convergence of the direct ALE ADER-DG schemes presented in this article is studied considering a test case proposed in \cite{HuShuTri}, which involves a smooth isentropic vortex flow that is furthermore convected on the horizontal plane $x-y$ with velocity $\mathbf{v}_c=(1,1,0)$. The initial computational domain is the square $\Omega(0)^{2D}=[0;10]\times[0;10]$ in 2D and the box $\Omega(0)^{3D}=[0;10]\times[0;10]\times[0;5]$ in 3D with periodic boundary conditions set everywhere. The initial condition is given by some perturbations $\delta$ that are superimposed onto a homogeneous background field $\U_0=(\rho,u,v,w,p)=(1,1,1,0,1)$, assuming that the entropy perturbation is zero, i.e. $S=\frac{p}{\rho^\gamma}=0$. The perturbations for density and pressure read
\begin{equation}
\label{rhopressDelta}
\delta \rho = (1+\delta T)^{\frac{1}{\gamma-1}}-1, \quad \delta p = (1+\delta T)^{\frac{\gamma}{\gamma-1}}-1, 
\end{equation}
with the temperature fluctuation $\delta T = -\frac{(\gamma-1)\epsilon^2}{8\gamma\pi^2}e^{1-r^2}$. According to \cite{HuShuTri}, the vortex strength is $\epsilon=5$ and the adiabatic index is set to $\gamma=1.4$, while the velocity field is affected by the following perturbations:
\begin{equation}
\label{ShuVortDelta}
\left(\begin{array}{c} \delta u \\ \delta v \\ \delta w \end{array}\right) = \frac{\epsilon}{2\pi}e^{\frac{1-r^2}{2}} \left(\begin{array}{c} -(y-5) \\ \phantom{-}(x-5) \\ 0 \end{array}\right).
\end{equation}
The exact solution $\Q_e$ can be simply obtained as the time-shifted initial condition, e.g. $\Q_e(\x,t_f)=\Q(\x-\v_c t_f,0)$, and the error is expressed in the continuous $L_2$ norm as  
\begin{equation}
  \epsilon_{L_2} = \sqrt{ \int \limits_{\Omega(t_f)} \left( \Q_e(\x,t_f) - \mathbf{u}_h(\x,t_f) \right)^2 d\x }.  
	\label{eqnL2error}
\end{equation}
Table \ref{tab.convEul} reports the convergence rates from second up to fourth order of accuracy for the vortex test problem run on a sequence of successively refined unstructured meshes. $h(\Omega(t_f))$ is the mesh size at the final time of the simulation $t_f=0.1$, which is taken to be the maximum diameter of the circumcircles or the circumspheres among all the control volumes of the final grid configuration $\Omega(t_f)$. The optimal order of accuracy is achieved both in space and time for $d=2$ as well as for $d=3$. Figure \ref{fig.SV-mesh} plots the two-dimensional mesh configuration at output times $t=0.5$, $t=1.0$, $t=1.5$ and $t=2.0$ for $N=3$.

\begin{table}[!htbp]  
\caption{Numerical convergence results for the compressible Euler equations using the direct ALE ADER-DG schemes from second up to fourth order of accuracy. The error norms refer to the variable $\rho$ (density) at time $t=0.1$.}  
\begin{center} 
\begin{small}
\renewcommand{\arraystretch}{1.0}
\begin{tabular}{c|cccccc} 
  \textbf{2D} & \multicolumn{2}{c}{$\mathcal{O}2$ ($N=1$)} & \multicolumn{2}{c}{$\mathcal{O}3$ ($N=2$)}  & \multicolumn{2}{c}{$\mathcal{O}4$ ($N=3$)}  \\
\hline
  $h(\Omega(t_f))$ & $\epsilon_{L_2}$ & $\mathcal{O}(L_2)$ & $\epsilon_{L_2}$ & $\mathcal{O}(L_2)$ &  $\epsilon_{L_2}$ & $\mathcal{O}(L_2)$ \\ 
\hline
3.26E-01 & 1.0004E-02 & -   & 7.5703E-04 & -   & 8.2888E-05 & -    \\ 
2.48E-01 & 5.4550E-03 & 2.2 & 3.1513E-04 & 3.2 & 1.8413E-05 & 5.5  \\ 
1.63E-01 & 2.4121E-03 & 2.0 & 9.7362E-05 & 2.8 & 4.1320E-06 & 3.6  \\ 
1.28E-01 & 1.3399E-03 & 2.4 & 4.1703E-05 & 3.5 & 1.3910E-06 & 4.5  \\ 
\hline
\multicolumn{7}{c}{} \\
\textbf{3D} & \multicolumn{2}{c}{$\mathcal{O}2$ ($N=1$)} & \multicolumn{2}{c}{$\mathcal{O}3$ ($N=2$)}  & \multicolumn{2}{c}{$\mathcal{O}4$ ($N=3$)}  \\
\hline
  $h(\Omega(t_f))$ & $\epsilon_{L_2}$ & $\mathcal{O}(L_2)$ & $\epsilon_{L_2}$ & $\mathcal{O}(L_2)$ &  $\epsilon_{L_2}$ & $\mathcal{O}(L_2)$ \\ 
\hline
5.92E-01 & 6.5631E-02 & -   & 9.7555E-03 & -   & 1.5405E-03 & -    \\ 
3.62E-01 & 2.6576E-02 & 1.8 & 2.4926E-03 & 2.8 & 3.3902E-04 & 3.1  \\ 
2.31E-01 & 1.1667E-02 & 1.8 & 7.5848E-04 & 2.7 & 3.8998E-05 & 4.8  \\ 
1.81E-01 & 6.5522E-03 & 2.3 & 3.8457E-04 & 2.8 & 1.2356E-05 & 4.7  \\ 
\hline 
\end{tabular}
\end{small}
\end{center}
\label{tab.convEul}
\end{table}

\subsection{The explosion problem} 
\label{sec.EP}
\begin{figure}[!htbp]
\begin{center}
\begin{tabular}{ccc} 
\includegraphics[width=0.47\textwidth]{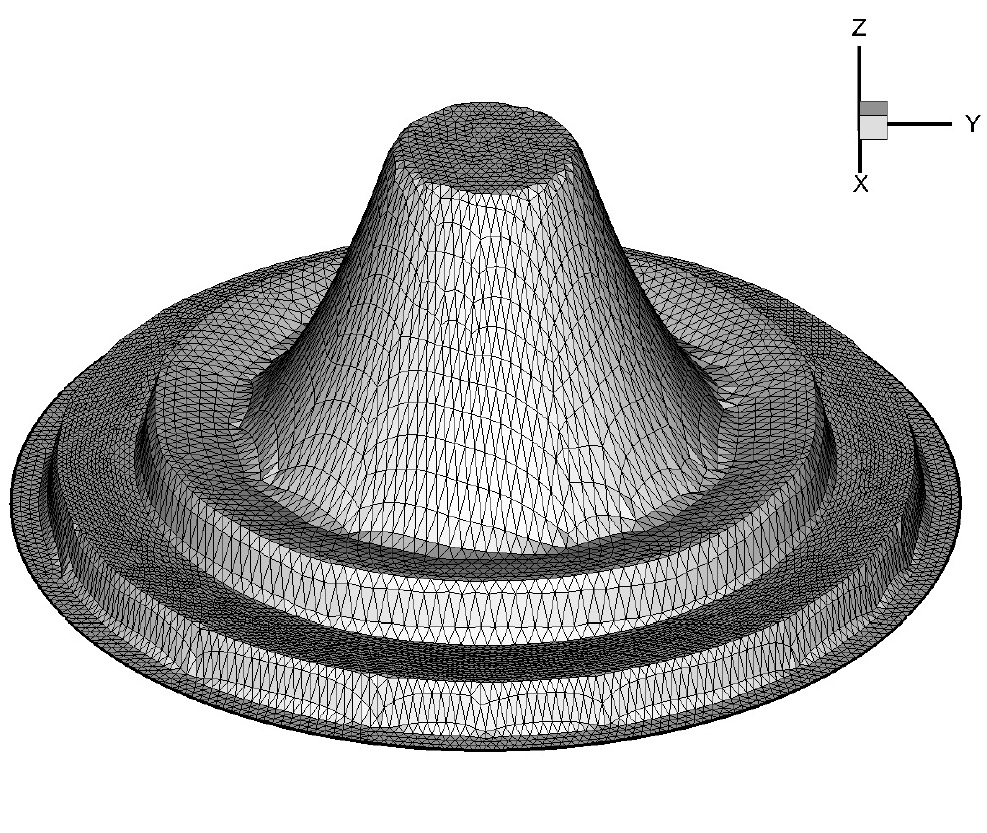}  &           
\includegraphics[width=0.47\textwidth]{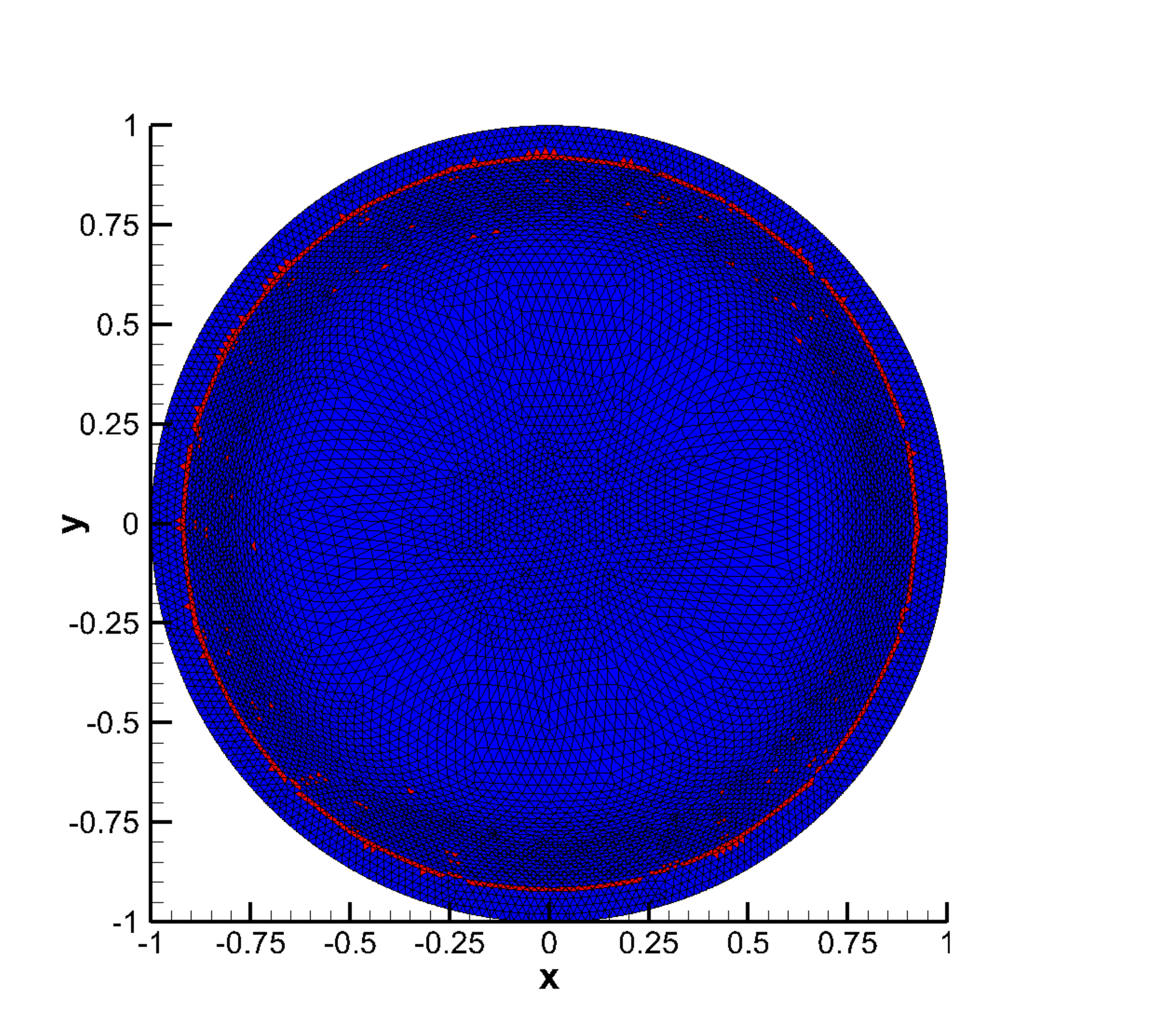} \\    
\includegraphics[width=0.47\textwidth]{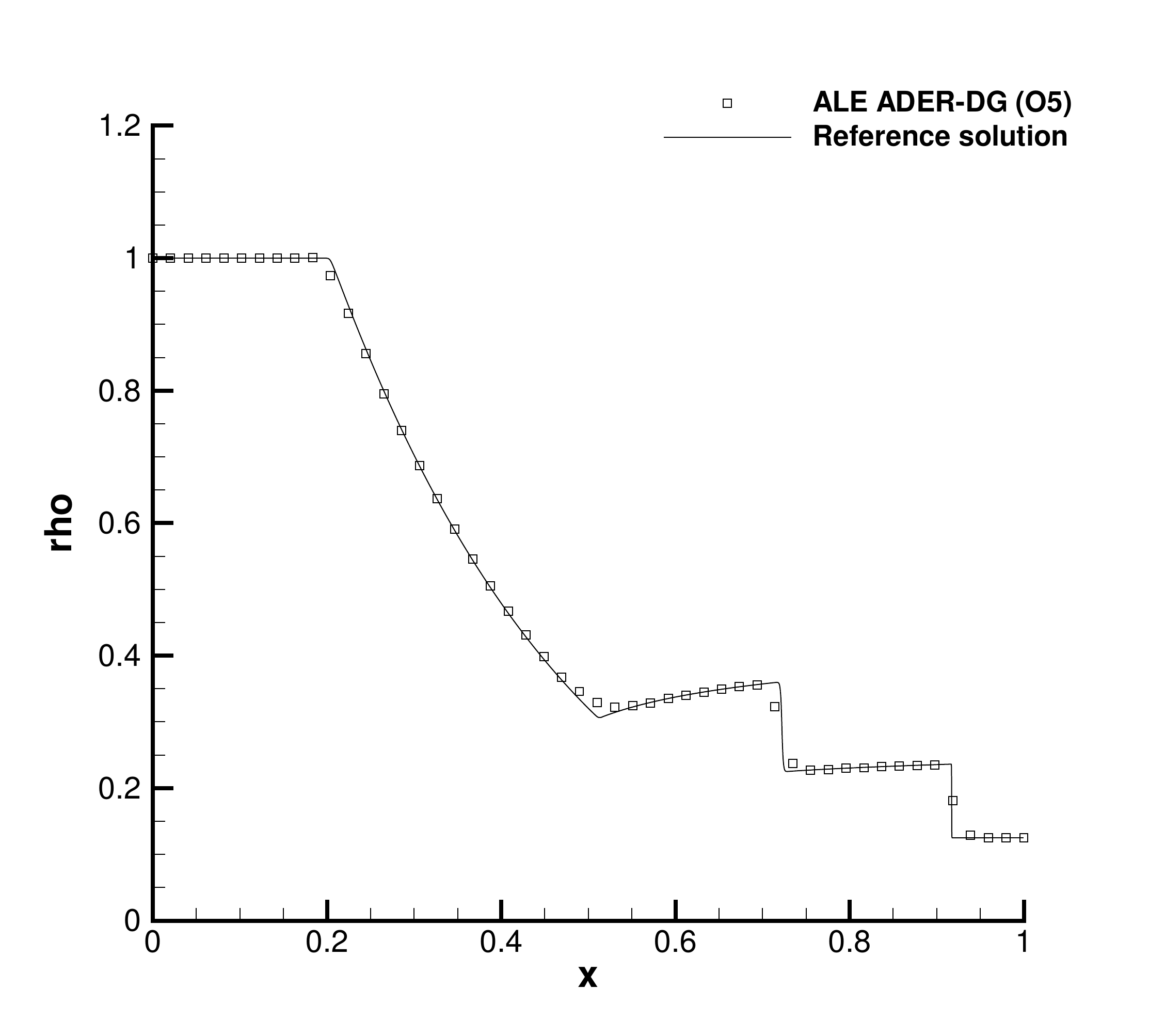}  &           
\includegraphics[width=0.47\textwidth]{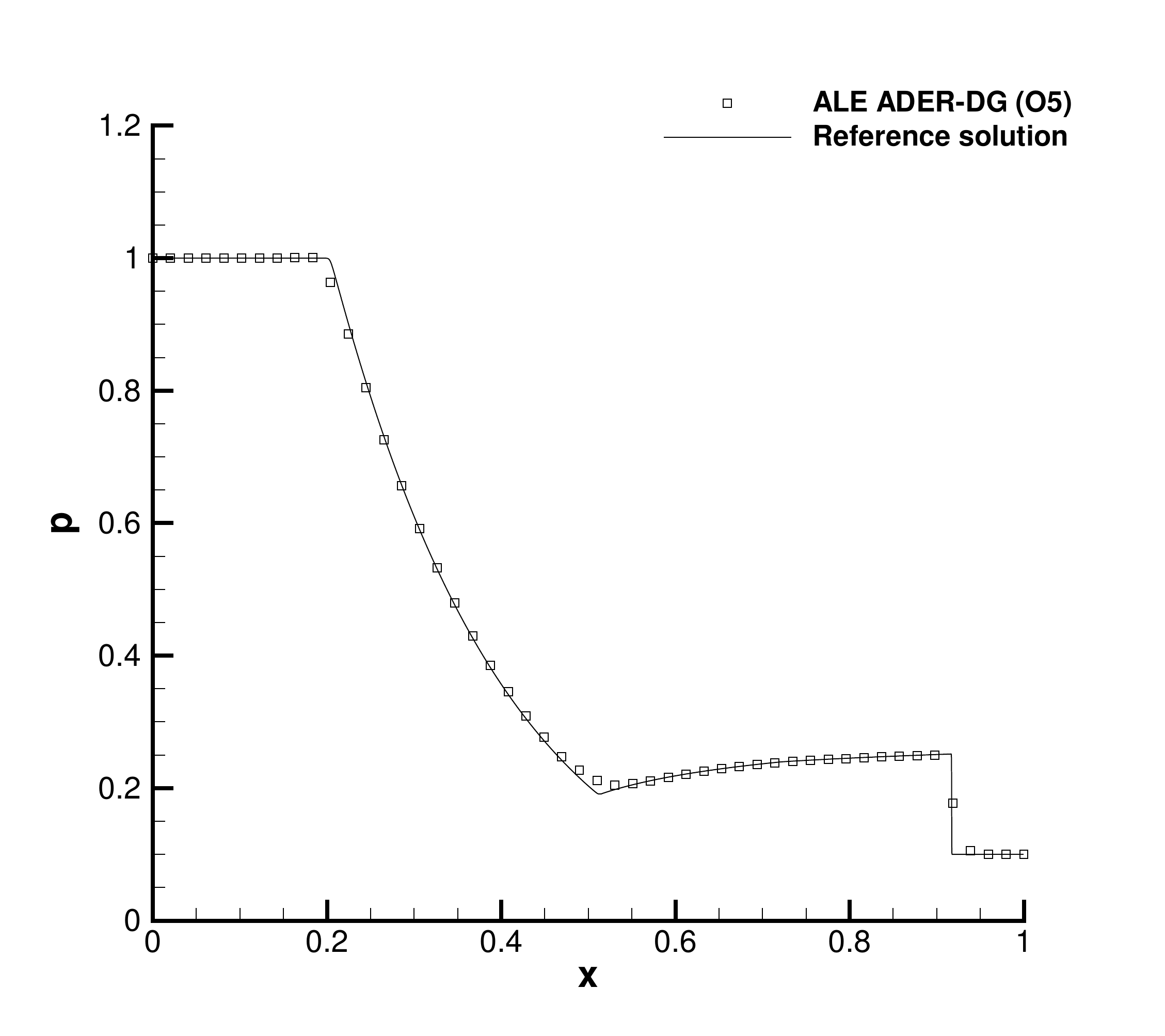} \\       
\end{tabular} 
\caption{Two-dimensional explosion problem at output time t=0.25 with $N=4$. Top: three-dimensional view of density distribution (left) and sub-cell limiter map (right). Bottom: 1D cut of density distribution (left) and pressure (right) compared against the reference solution. } 
\label{fig.EP}
\end{center}
\end{figure}

The first test case considered in this work is the multidimensional explosion problem. It represents a useful sanity check because it involves a rarefaction wave moving towards the center of the computational domain as well as a contact discontinuity and a shock wave that are traveling to the opposite direction. The initial computational domain is a circle or a sphere of radius $R=1$ and the initial condition is composed by two different states separated at radius $R_s=0.5$. The inner state $\U_i$ and the outer state $\U_o$ read
\begin{equation}
  \U(\x,0) = \left\{ \begin{array}{lcc} \U_i=(1,0,0,0,1),       & \textnormal{ if } & \left\| r \right\| \leq R_s \\ 
                                        \U_o=(0.125,0,0,0,0.1), & \textnormal{ if } & \left\| r \right\| > R_s        
                      \end{array}  \right. ,
\end{equation} 
with $r=\sqrt{\x^2}$ representing the generic radial position. Transmissive boundaries have been set on the external side and the domain is paved with $N_E=17340$ triangles. We set $\gamma=1.4$ and, at the final time of the simulation $t_f=0.25$, the fifth order accurate numerical solution is compared against the reference solution, whose derivation can be found in \cite{toro-book,Lagrange2D}. In Figure \ref{fig.EP} one can note an excellent resolution of the contact wave and the sub-cell limiter map shows that the limiter is active only across the shock wave.

\begin{figure}[!htbp]
\begin{center}
\begin{tabular}{ccc} 
\includegraphics[width=0.4\textwidth]{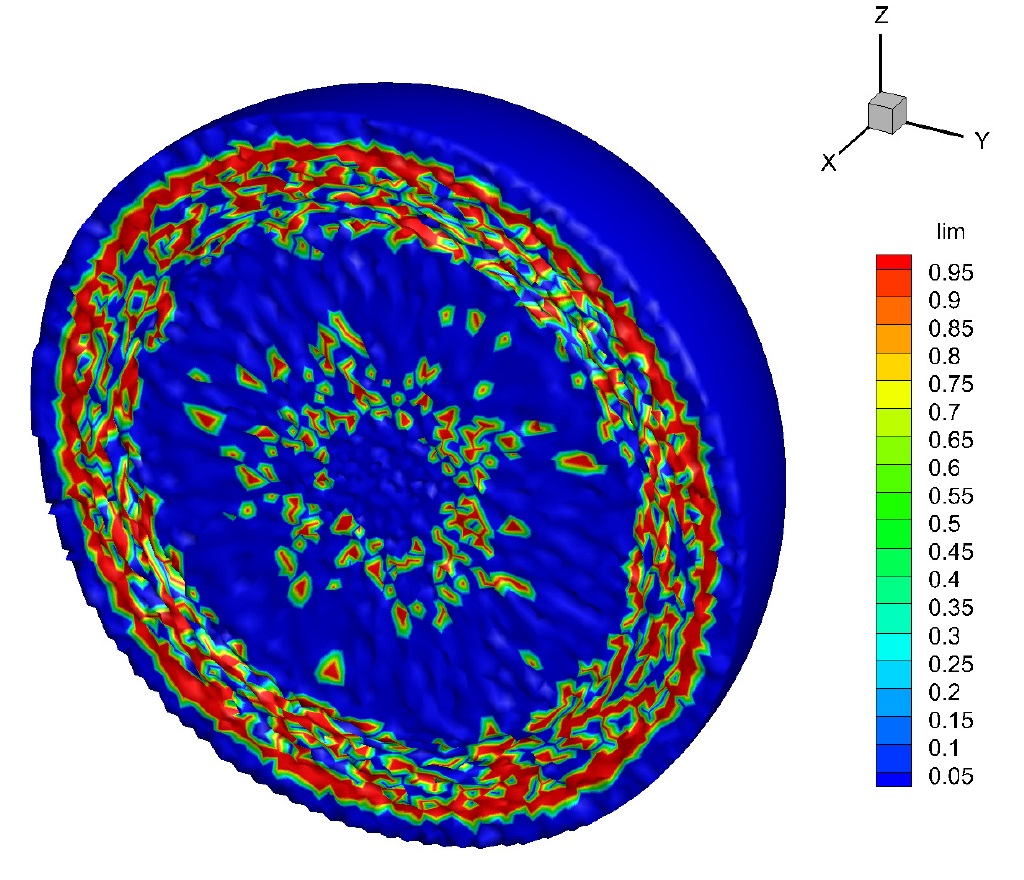}  &           
\includegraphics[width=0.4\textwidth]{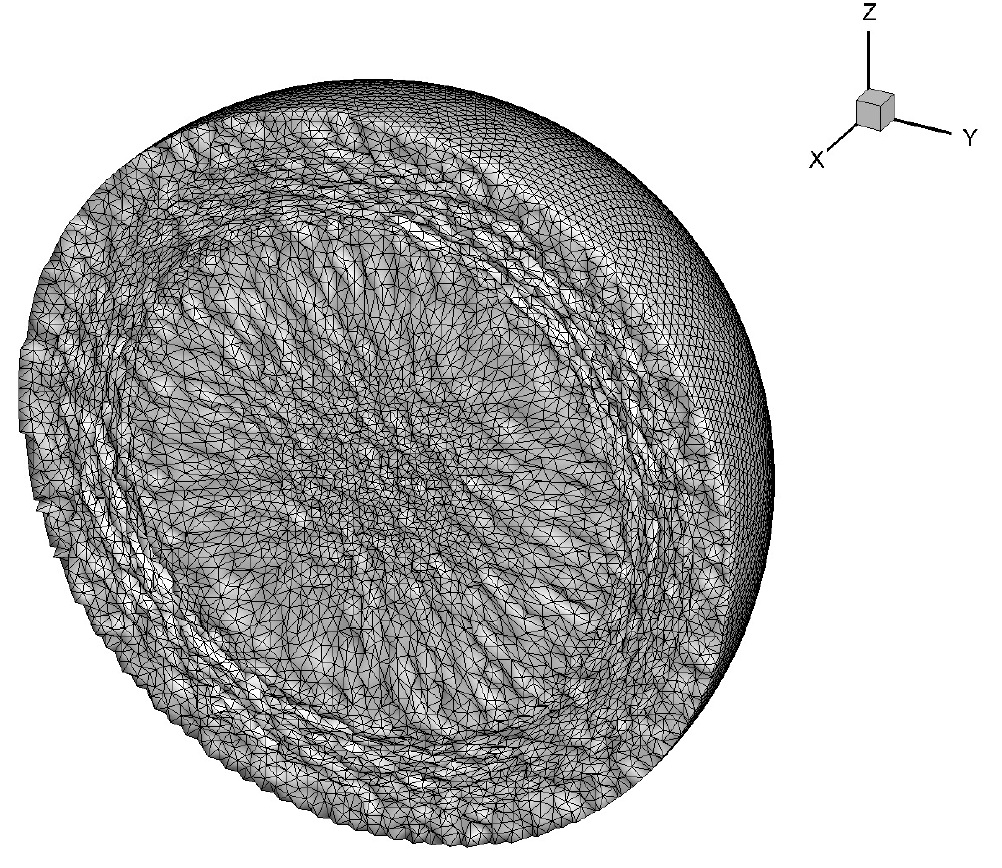} \\    
\includegraphics[width=0.4\textwidth]{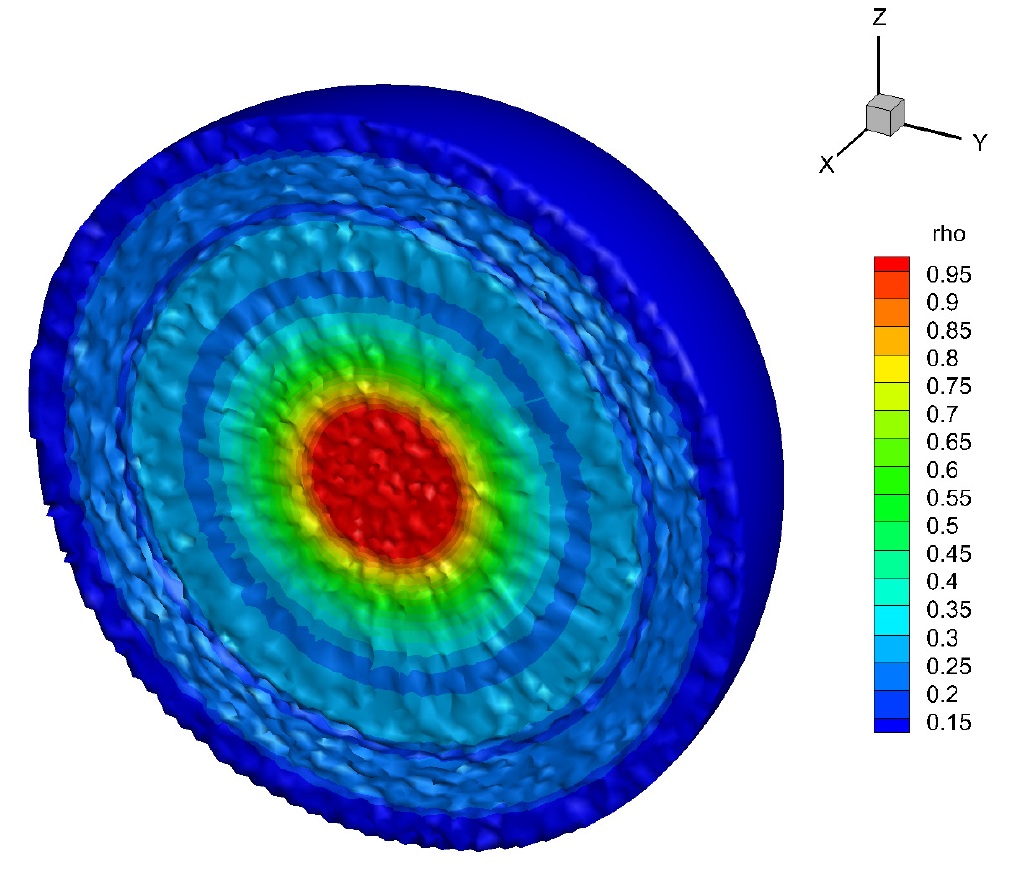}  &           
\includegraphics[width=0.4\textwidth]{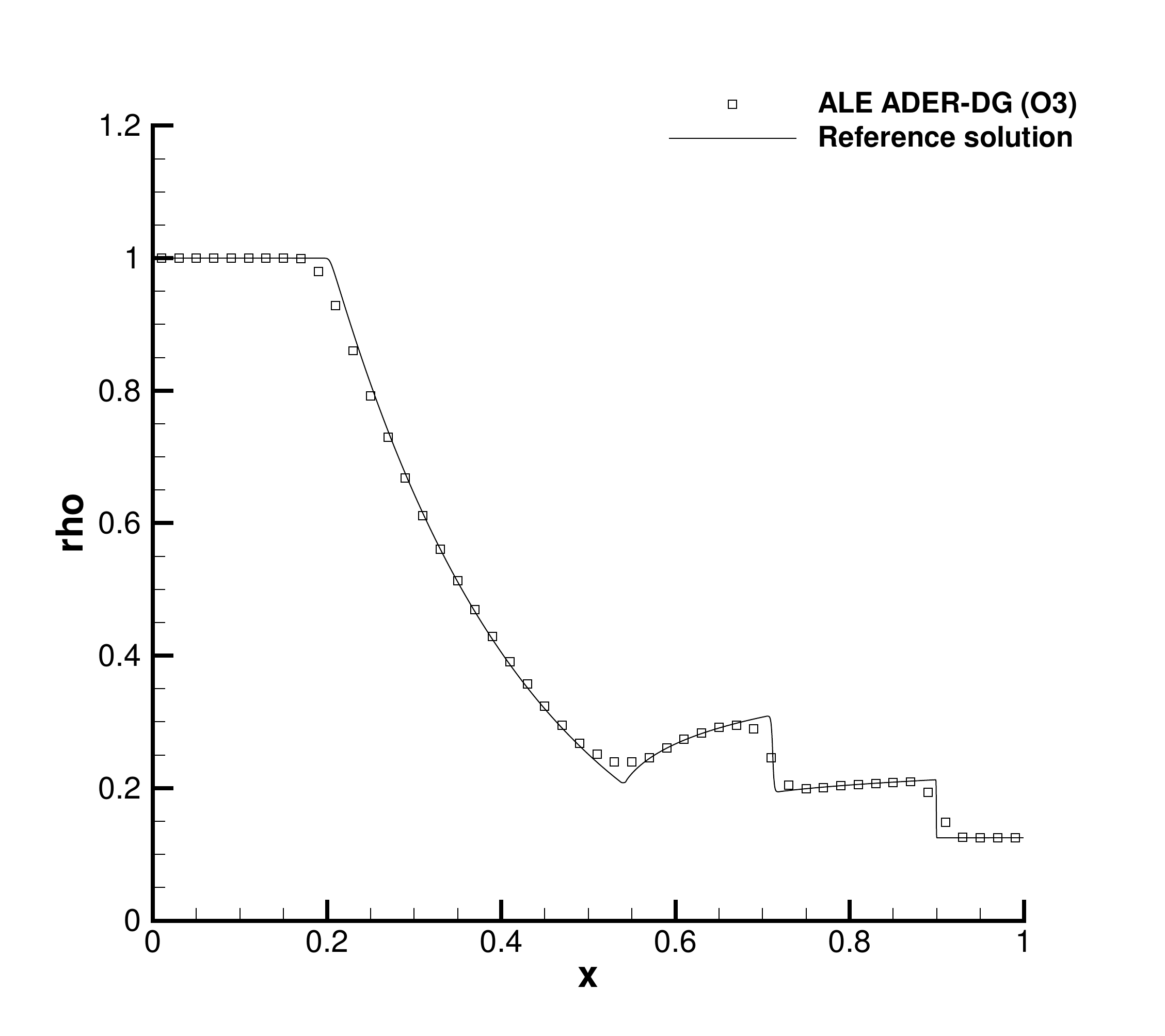} \\   
\includegraphics[width=0.4\textwidth]{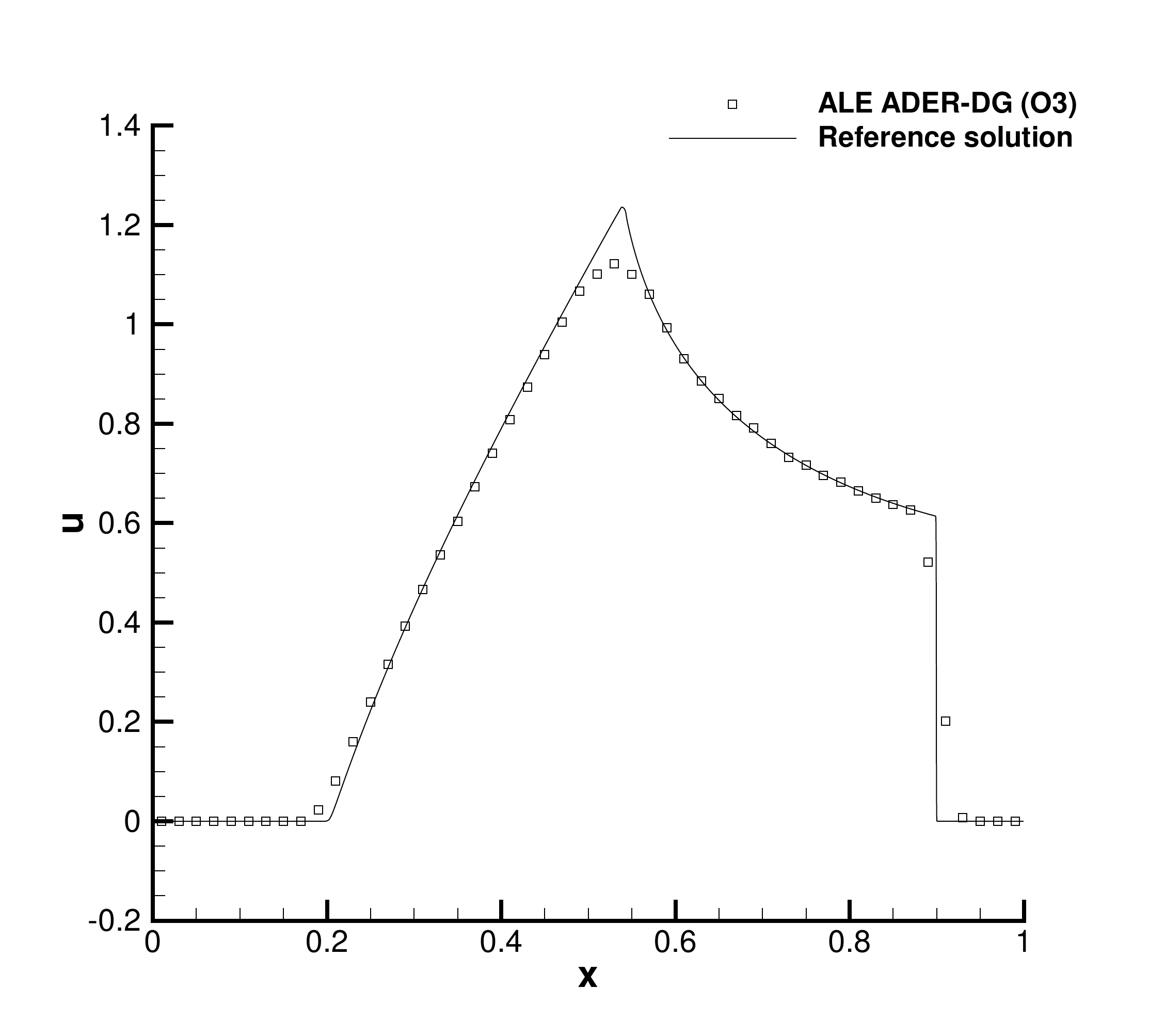}  &           
\includegraphics[width=0.4\textwidth]{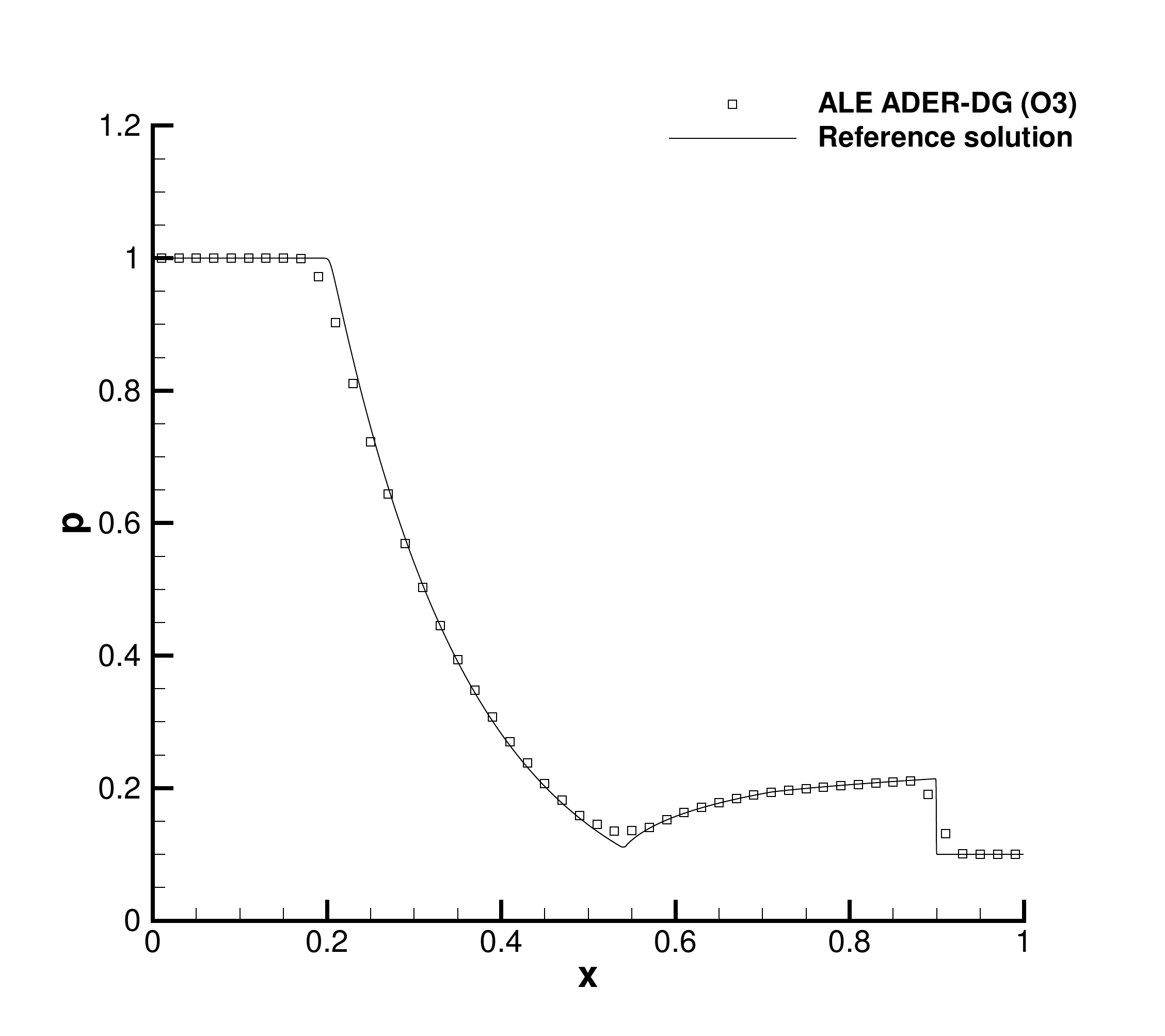} \\    
\end{tabular} 
\caption{Three-dimensional explosion problem at output time t=0.25 with $N=2$. Top row: sub-cell limiter map (left) and final mesh configuration (right). Middle row: three-dimensional view of density distribution (left) and 1D cut of density distribution (right). Bottom row: 1D cut of horizontal velocity (left) and pressure (right) distribution compared against the reference solution.} 
\label{fig.EP3D}
\end{center}
\end{figure}

Figure \ref{fig.EP3D} shows the results obtained running a third order three-dimensional simulation of the explosion problem. A total number of $N_E=1469472$ tetrahedra has been used to discretize the sphere. The discontinuity between internal and external state is located again at $R_s=0.5$. Also in this case the limiter well detects the region around the shock and a good agreement with the exact solution is achieved.

\subsection{The Saltzman problem} 
\label{sec.Saltzman}
\begin{figure}[!htbp]
\begin{center}
\begin{tabular}{cc} 
\includegraphics[width=0.47\textwidth]{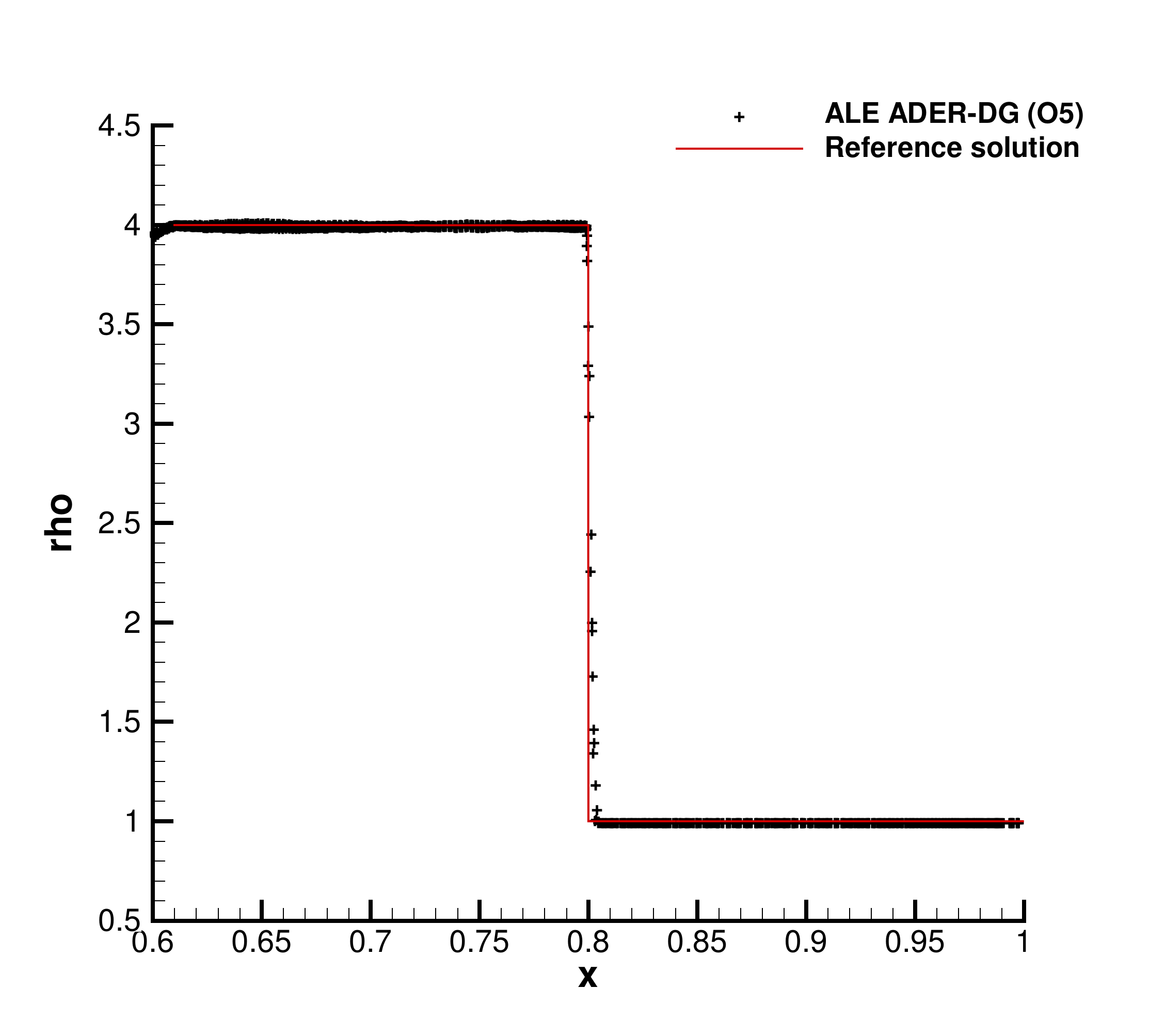}  &           
\includegraphics[width=0.47\textwidth]{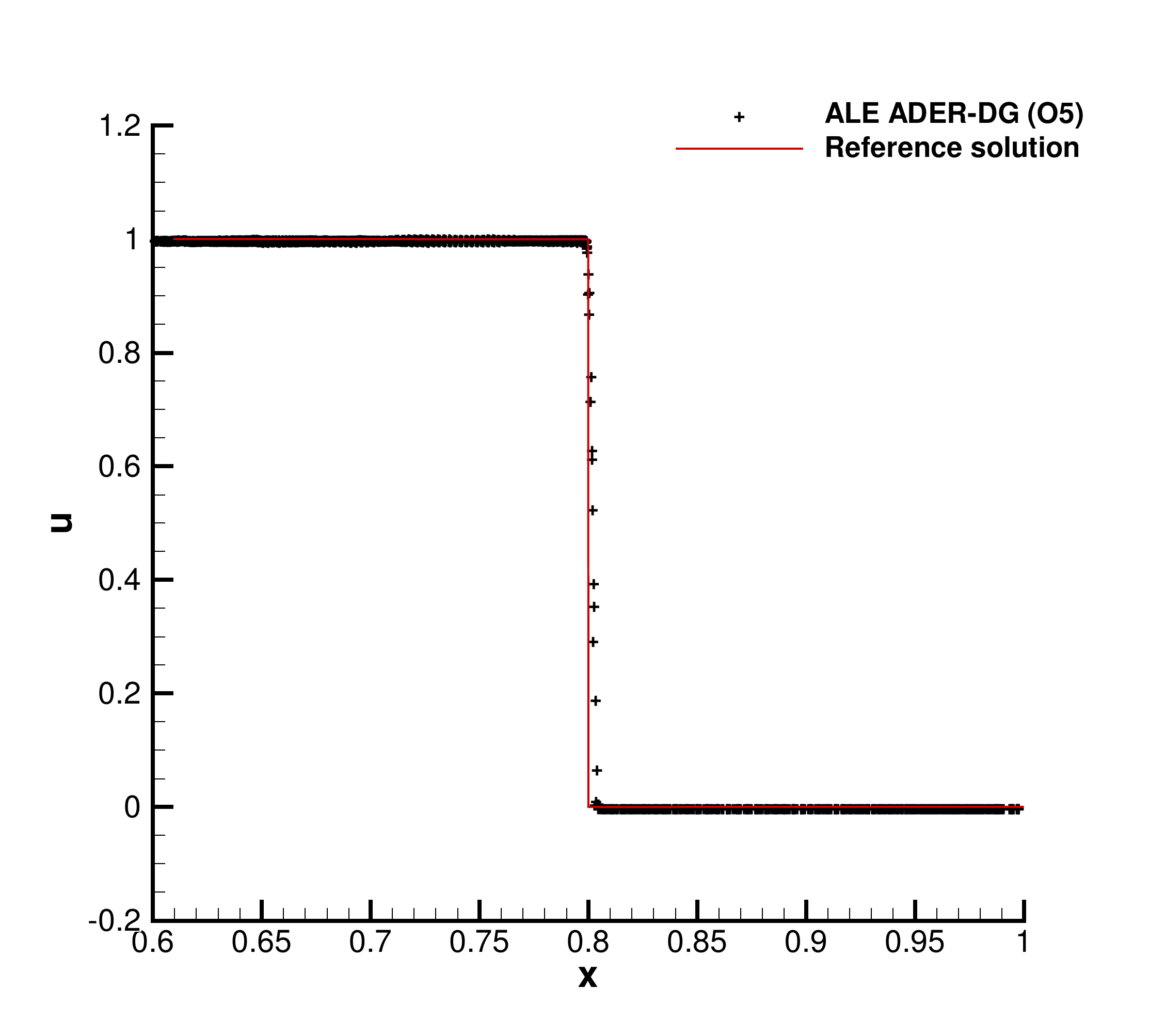} \\   
\includegraphics[width=0.47\textwidth]{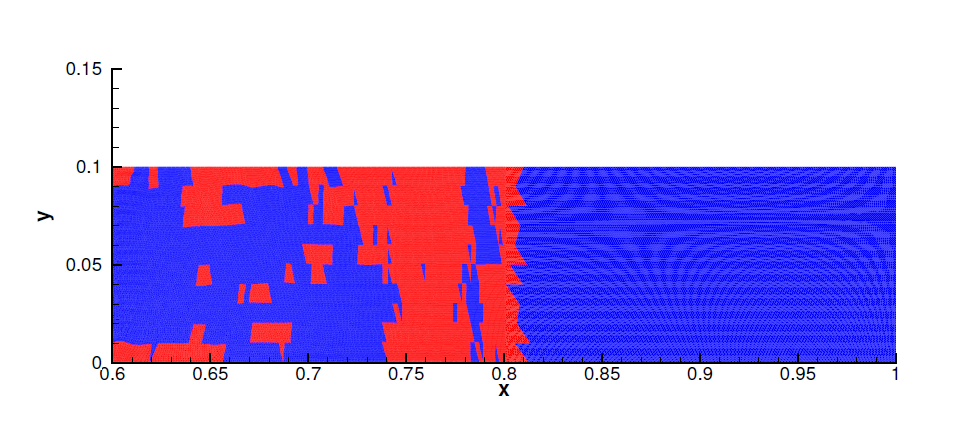}  &           
\includegraphics[width=0.47\textwidth]{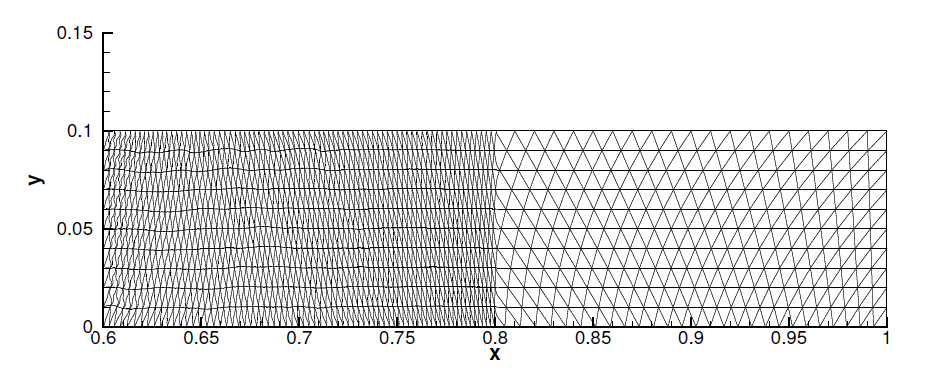} \\        
\end{tabular} 
\caption{Saltzman problem in 2D with $N=4$. Top: scatter plot of the cell density (left) and horizontal velocity (right) as a function of cell horizontal coordinate $x$ versus the exact solution. Bottom: sub-cell limiter map (left) and final mesh configuration (right).} 
\label{fig.Saltz2D}
\end{center}
\end{figure}

The Saltzman test case describes the motion of a piston which is impinging on a fluid at rest contained in the initial computational domain given by $\Omega(0)=[0;1]\times[0;0.1]$. This is a challenging test problem used in literature \cite{Maire2009,chengshu2,SaltzmanOrg} to assess the robustness of any Lagrangian algorithm. The piston moves with velocity $\v_p=(1,0)$ and generates a strong shock wave that is traveling along the main direction of the computational domain. The cells which lie near the piston are highly compressed during the simulation. Wall boundaries are considered everywhere except for the piston, on which we impose a moving slip wall boundary condition. The computational mesh is composed by $N_E=2000$ triangles and the grid is initially distorted applying a skewness that makes no sides of the mesh aligned with the main fluid flow, as fully explained in \cite{SaltzmanOrg}. At time $t=0$ the domain is filled with a perfect gas at rest with $\gamma=\frac{5}{3}$, uniform density $\rho=1$ and pressure $p=10^{-4}$, according to \cite{chengshu2}. The final time of the simulation is $t_f=0.6$ and the exact solution is given by a one-dimensional infinite strength shock wave with a post shock density of $\rho_e = 4.0$ and the shock front located at $x=0.8$, see \cite{Lagrange2D} for further details. Figure \ref{fig.Saltz2D} shows the two-dimensional results, highlighting the excellent agreement between the fifth order accurate numerical solution and the exact solution. Note that even the well-known wall heating effect close to the moving piston \cite{toro.anomalies.2002} almost disappears in this case, without any specific treatment. Furthermore, we also plot the sub-cell limiter map that marks in red those cells in which the limiter is active and in blue the unlimited elements.

\begin{figure}[!htbp]
\begin{center}
\begin{tabular}{cc} 
\includegraphics[width=0.47\textwidth]{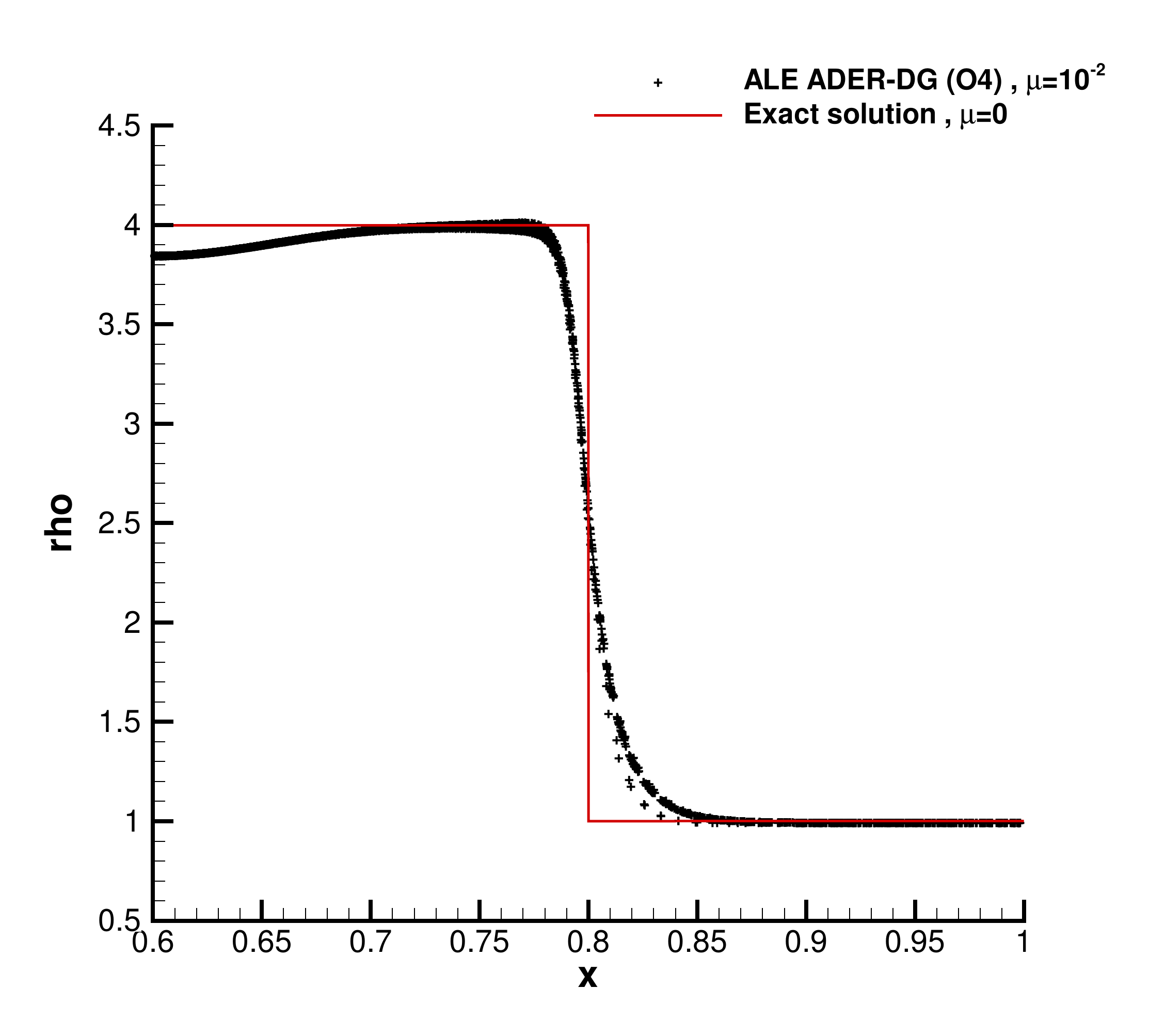}  &           
\includegraphics[width=0.47\textwidth]{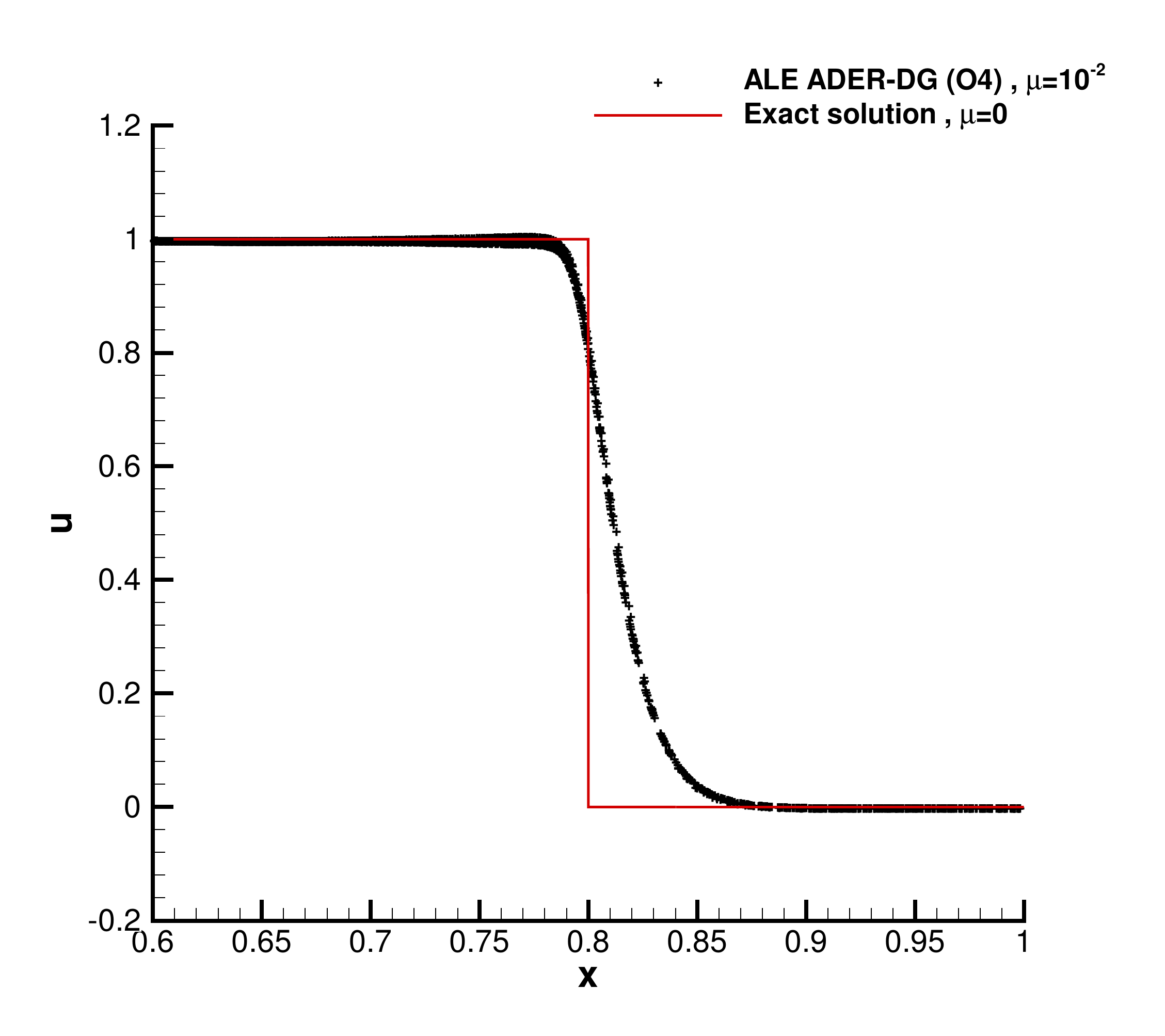} \\   
\includegraphics[width=0.47\textwidth]{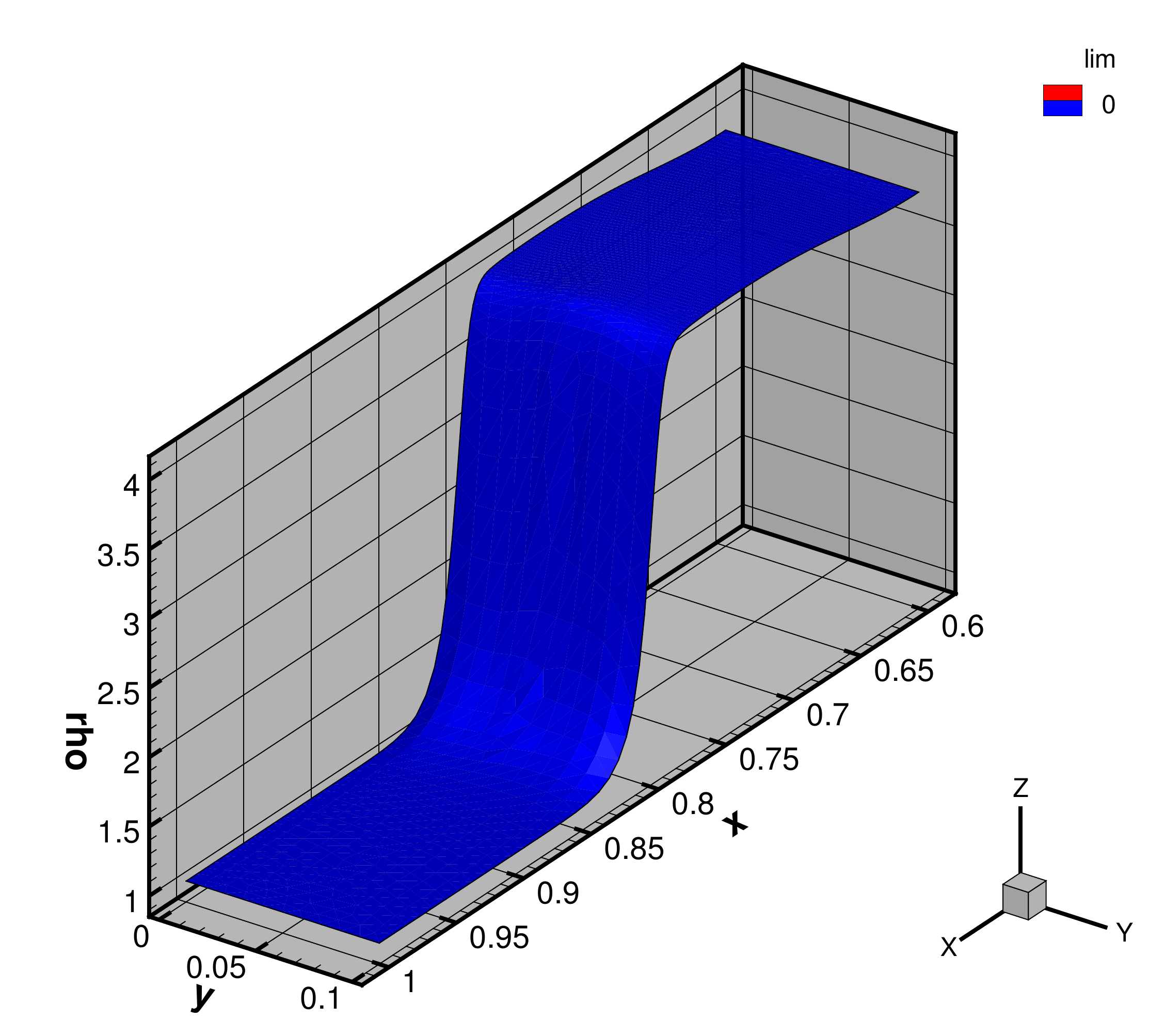}  &           
\includegraphics[width=0.47\textwidth]{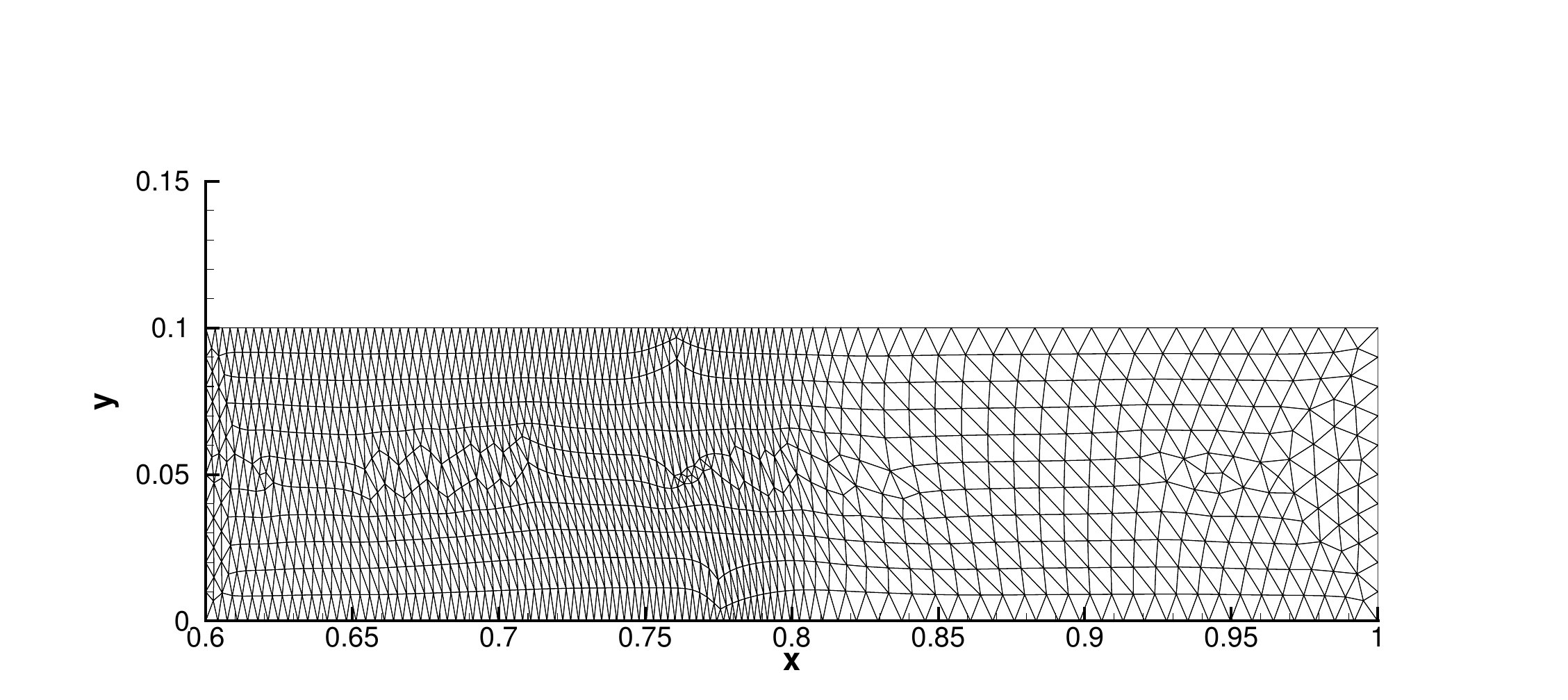} \\        
\end{tabular} 
\caption{Viscous Saltzman problem in 2D with $N=3$ and $\mu=10^{-2}$. Top: scatter plot of the cell density (left) and horizontal velocity (right) as a function of cell horizontal coordinate $x$ versus the exact solution for the inviscid flow. Bottom: sub-cell limiter map (left) and final mesh configuration (right).} 
\label{fig.Saltz2D-100}
\end{center}
\end{figure}

\begin{figure}[!htbp]
\begin{center}
\begin{tabular}{cc} 
\includegraphics[width=0.47\textwidth]{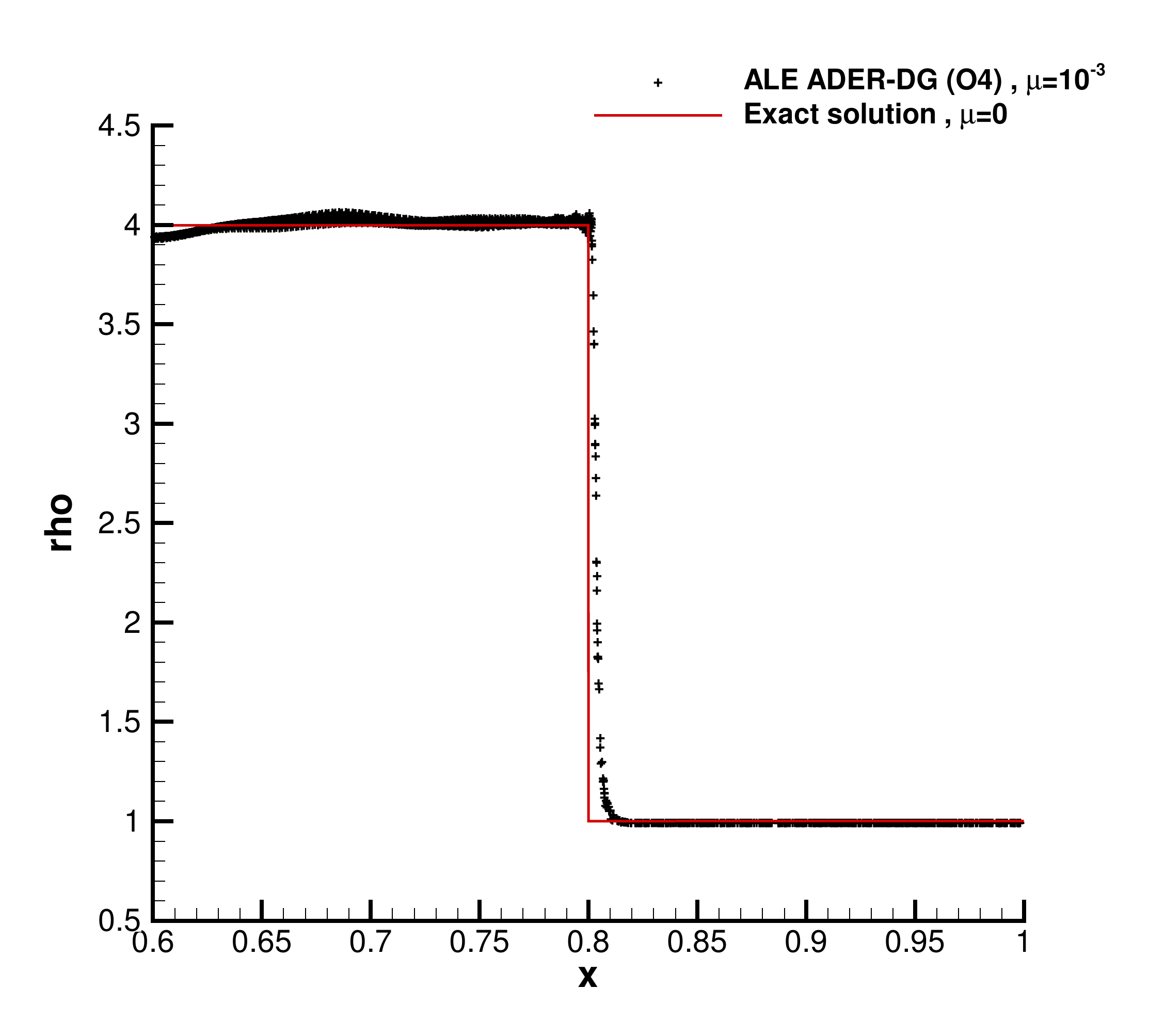}  &           
\includegraphics[width=0.47\textwidth]{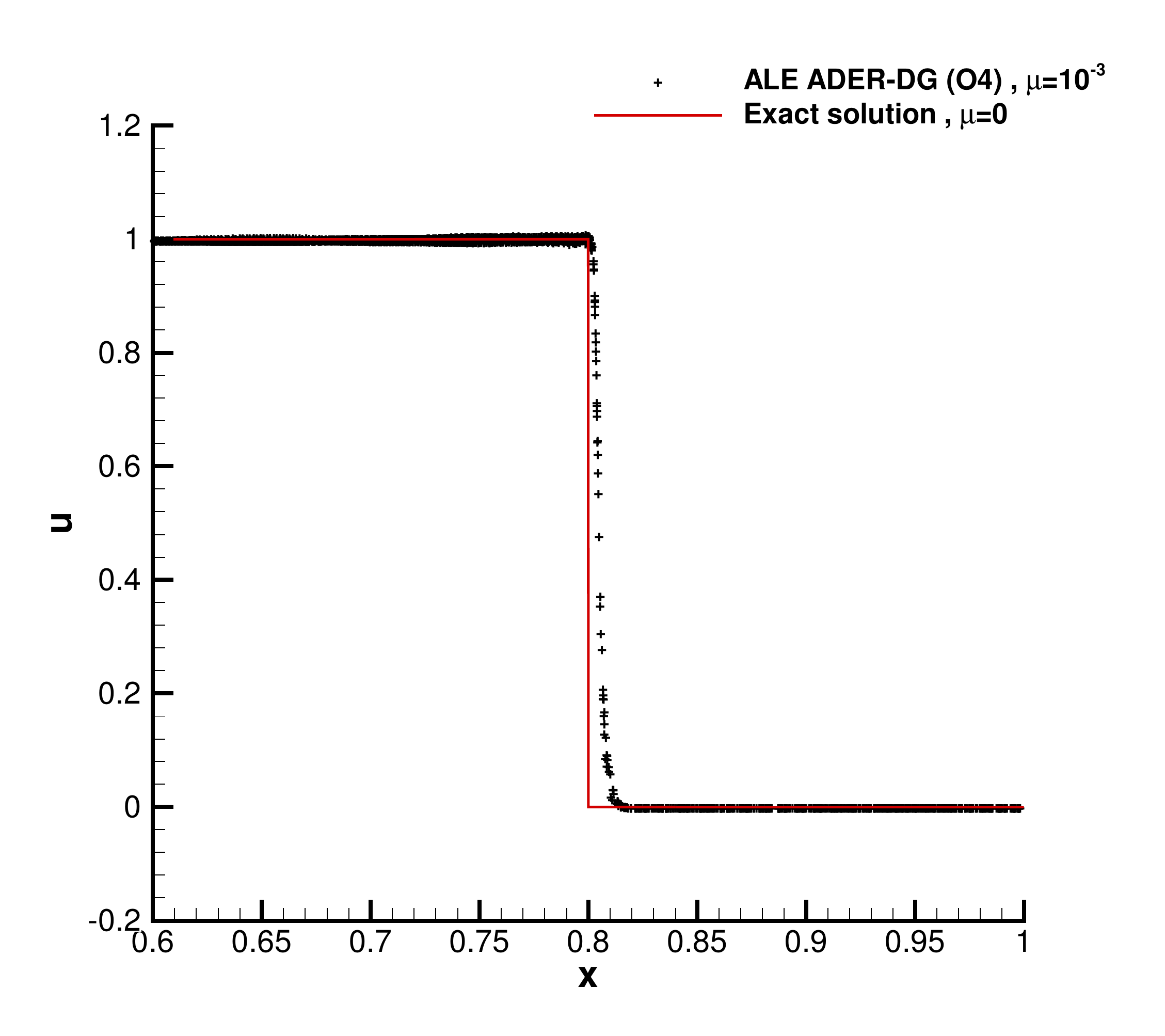} \\   
\includegraphics[width=0.47\textwidth]{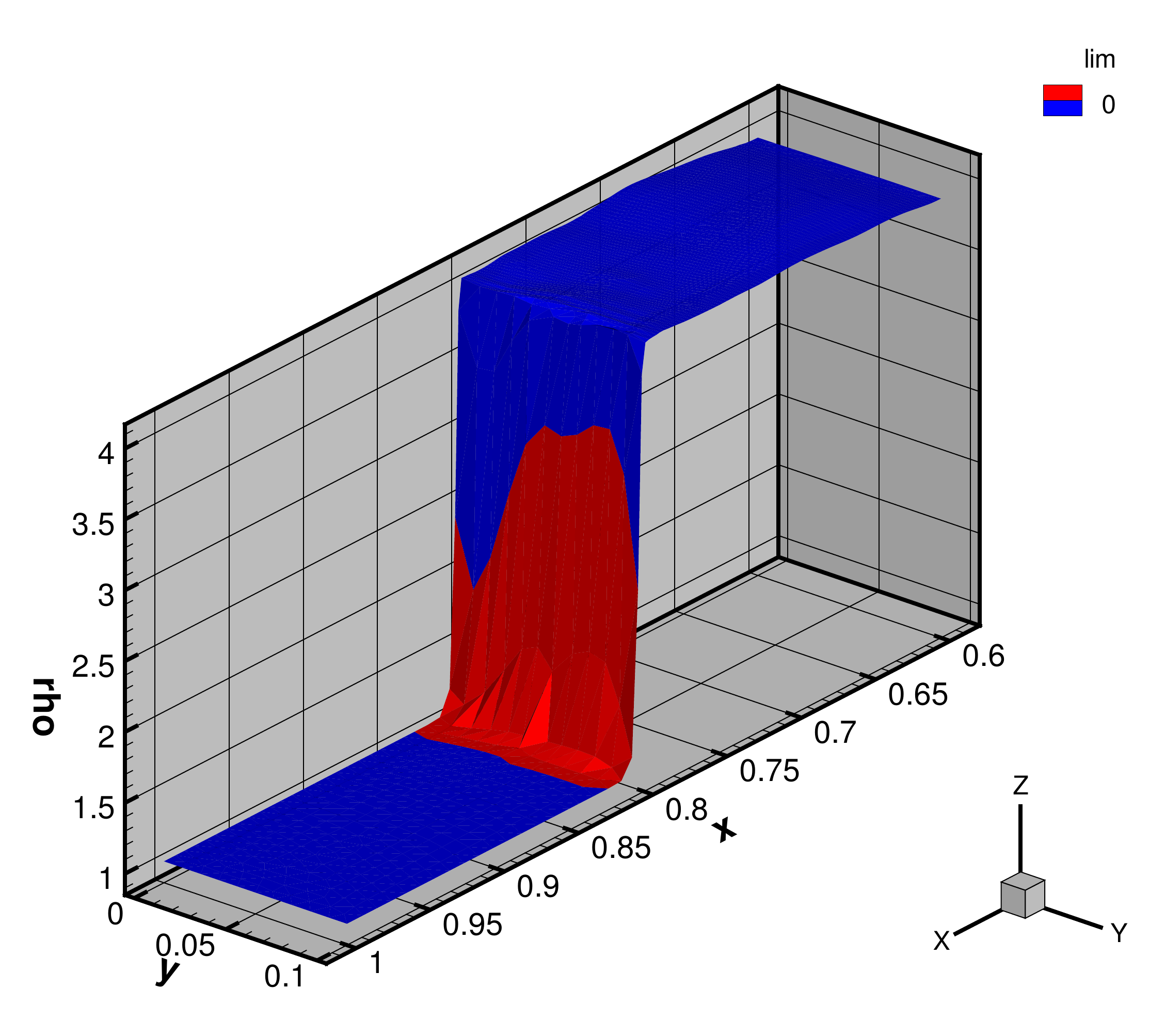}  &           
\includegraphics[width=0.47\textwidth]{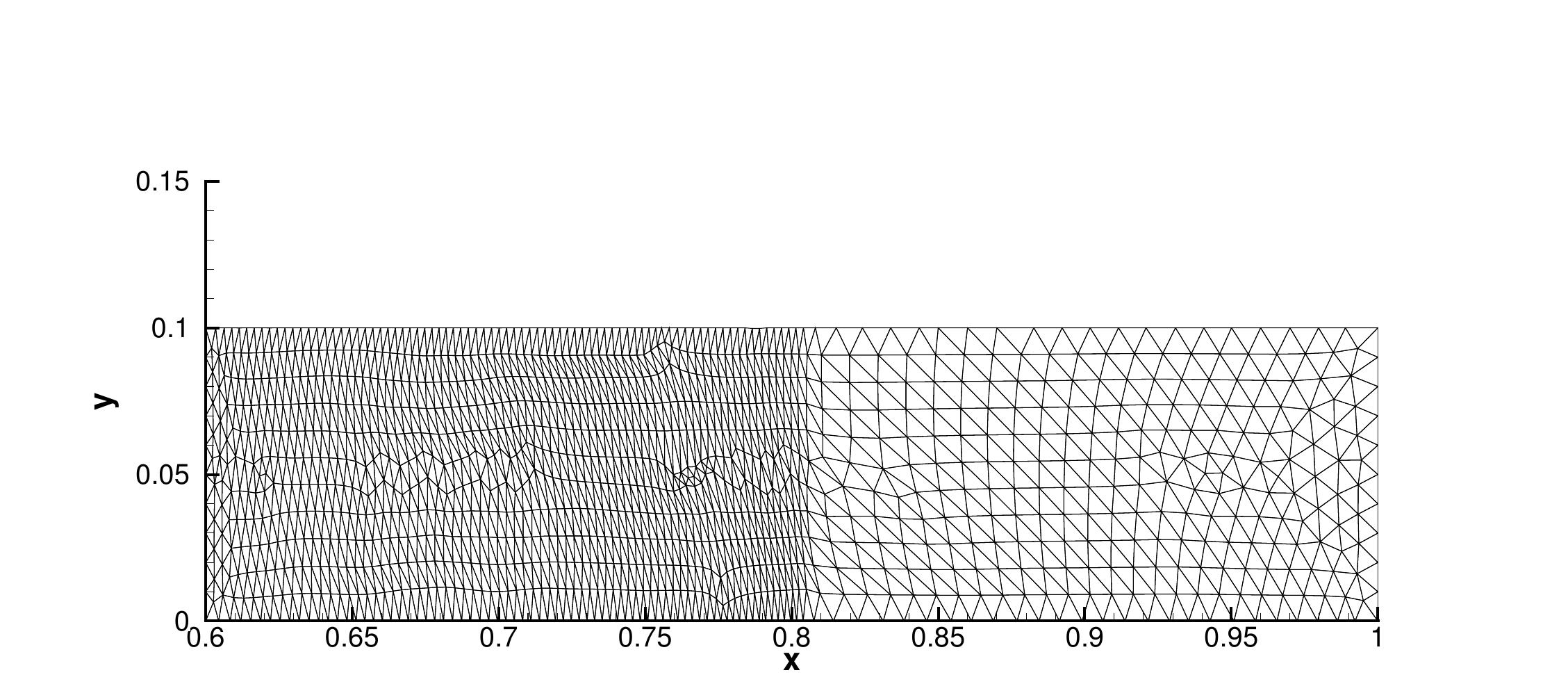} \\        
\end{tabular} 
\caption{Viscous Saltzman problem in 2D with $N=3$ and $\mu=10^{-3}$. Top: scatter plot of the cell density (left) and horizontal velocity (right) as a function of cell horizontal coordinate $x$ versus the exact solution for the inviscid flow. Bottom: sub-cell limiter map (left) and final mesh configuration (right).} 
\label{fig.Saltz2D-1000}
\end{center}
\end{figure}

Next, we propose to solve the Saltzman problem with physical viscosity, hence considering the compressible Navier-Stokes equations \eqref{eqn.NSflux}. The setting of the problem is the same one used for the inviscid case, but the initial computational mesh counts a total number of $N_E=2220$ fully unstructured triangles, that have been skewed according to the transformation explicitly given in \cite{SaltzmanOrg}. Slip-wall boundaries have been imposed on the lateral side of the domain in order to avoid the generation and the growth of the boundary layer. Figure \ref{fig.Saltz2D-100} shows the numerical results obtained with a viscosity coefficient of $\mu=10^{-2}$, hence leading to a Reynolds number of $Re=100$, while we plot a fourth order simulation of the Saltzman problem with $Re=1000$ in Figure \ref{fig.Saltz2D-1000}. The physical viscosity spreads the shock wave induced by the piston, so that the sub-cell limiter is not needed for $Re=100$, or it becomes active in very few cells with $Re=1000$. For both viscous simulations we have used the isoparametric version of our ALE ADER-DG schemes with a CFL number of $\textnormal{CFL}=0.1$. The results are compared against the exact solution obtained in the \textit{inviscid} case.

\subsection{The Kidder problem} 
\label{sec.Kidder}

This is a smooth test case proposed in \cite{Kidder1976} that considers the isentropic compression of a portion of a shell filled with an ideal gas. According to \cite{Maire2009,Despres2009}, the initial computational domain is bounded by $r_i(t) \leq r \leq r_e(t)$, where $r_i(t),r_e(t)$ represent the time-dependent internal and external radius, respectively, and $r=\sqrt{\x^2}$ denotes as usual the generic radial coordinate. Sliding wall boundaries are imposed everywhere apart from the internal and external frontiers, where we set a space-time dependent state computed according to the self-similar analytical solution $R(r,t)$, available in \cite{Kidder1976}. The gas is initially assigned a uniform entropy $s_0= \frac{p_0}{\rho_0^\gamma} = 1$ with the adiabatic index $\gamma=2$ and the initial condition
\begin{equation}
\left( \begin{array}{c} \rho_0(r) \\ \mathbf{v}_0(r) \\ p_0(r) \end{array}  \right) = \left( \begin{array}{c}  \left(\frac{r_{e,0}^2-r^2}{r_{e,0}^2-r_{i,0}^2}\rho_{i,0}^{\gamma-1}+\frac{r^2-r_{i,0}^2}{r_{e,0}^2-r_{e,0}^2}\rho_{e,0}^{\gamma-1}\right)^{\frac{1}{\gamma-1}} \\ 0 \\ s_0\rho_0(r)^\gamma \end{array}  \right), 
\label{eq:KidderIC}
\end{equation}
where $\rho_{i,0}=1$ and $\rho_{e,0}=2$ denote the initial values of density at the time-dependent internal and external frontier, respectively. The final time is taken to be $t_f=\frac{\sqrt{3}}{2}\tau$ with the focalisation time $\tau=0.217944947177$ computed according to \cite{Kidder1976,Lagrange2D,Lagrange3D}. The exact location of the shell at the final time is bounded by $0.45 \leq R \leq 0.5$, therefore the absolute error $|err|$ between analytical and numerical solution is easily computed and reported in Table \ref{tab.Kidder}. Since the Kidder problem does not involve any discontinuity, no cells are affected by the sub-cell limiter.

\begin{figure}[!htbp]
	\begin{center}
	\begin{tabular}{cc} 
	\includegraphics[width=0.47\textwidth]{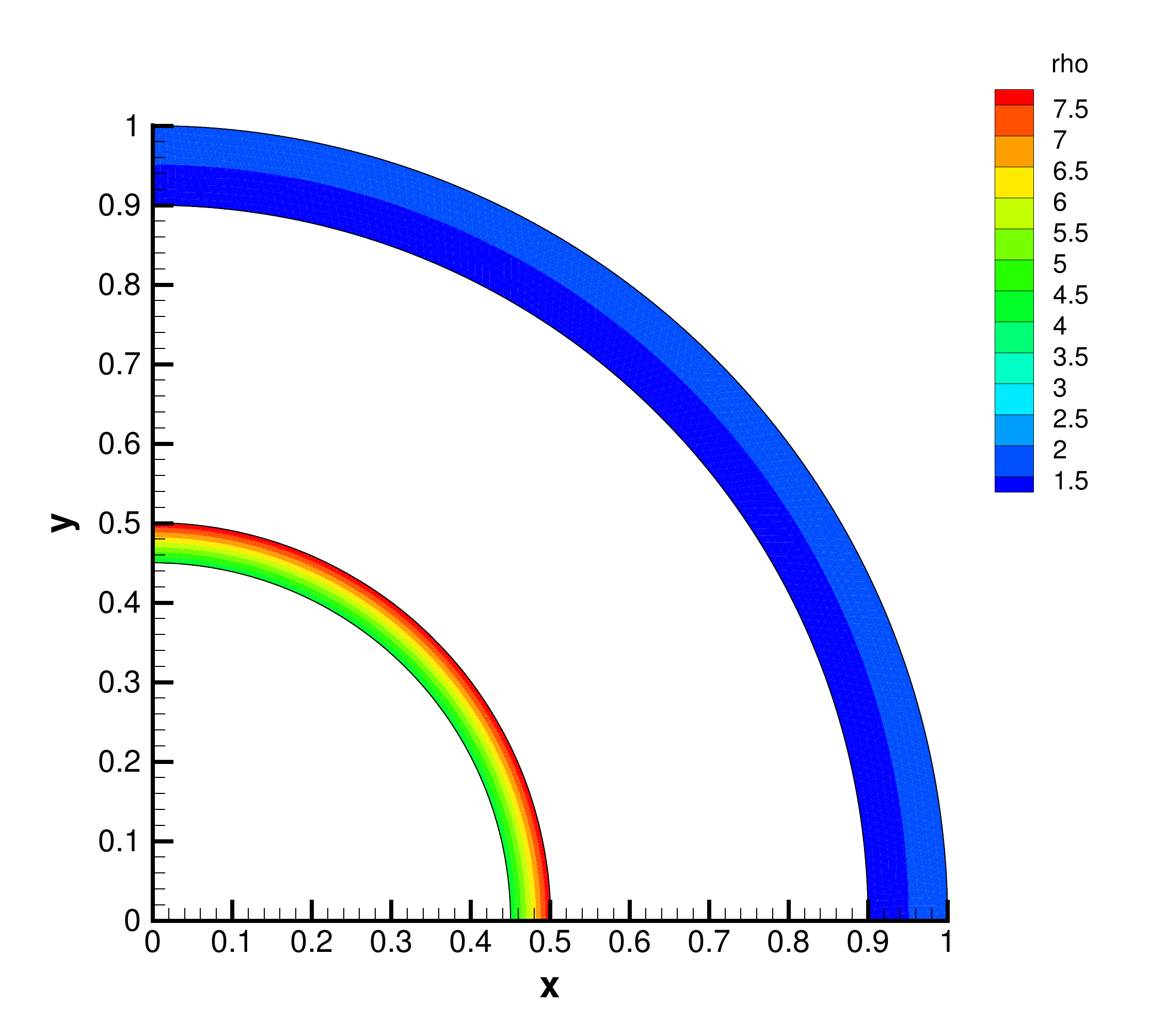}  &           
	\includegraphics[width=0.47\textwidth]{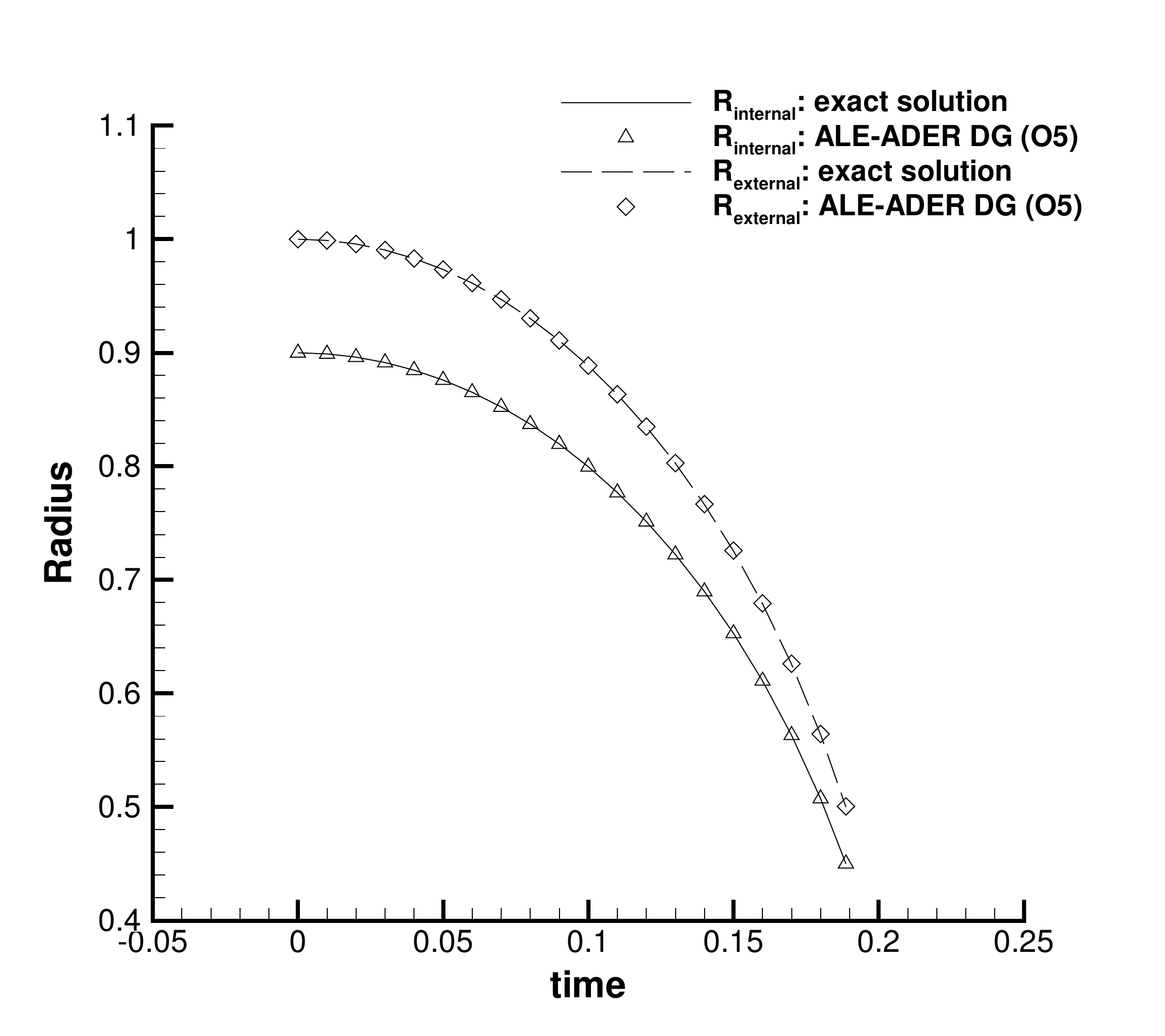} \\
	\end{tabular}
	\caption{Final computational domain (left) and evolution of the internal and external radius of the shell and comparison between analytical and numerical solution (right).}
	\label{fig.Kidder}
	\end{center}
\end{figure}

\begin{table}[!htbp]
 \centering
  \begin{center}
    \begin{tabular}{c|cc}
      $r_{ex}$ 	& $r_{num}$  & $|err|$   \\
      \hline
      0.450000 	& 0.45014   &  1.4e-4    \\
      \hline
      0.500000 	& 0.50041   &  4.1e-4    \\ 
      \hline
    \end{tabular}
  \end{center}
  \caption{Kidder problem. Absolute error for the internal and external radius location between exact ($r_{ex}$) and numerical ($r_{num}$) 
    solution.}
  \label{tab.Kidder}
\end{table}

\subsection{The Sedov problem} 
\label{sec.Sedov}
\begin{figure}[!htbp]
\begin{center}
\begin{tabular}{ccc} 
\includegraphics[width=0.47\textwidth]{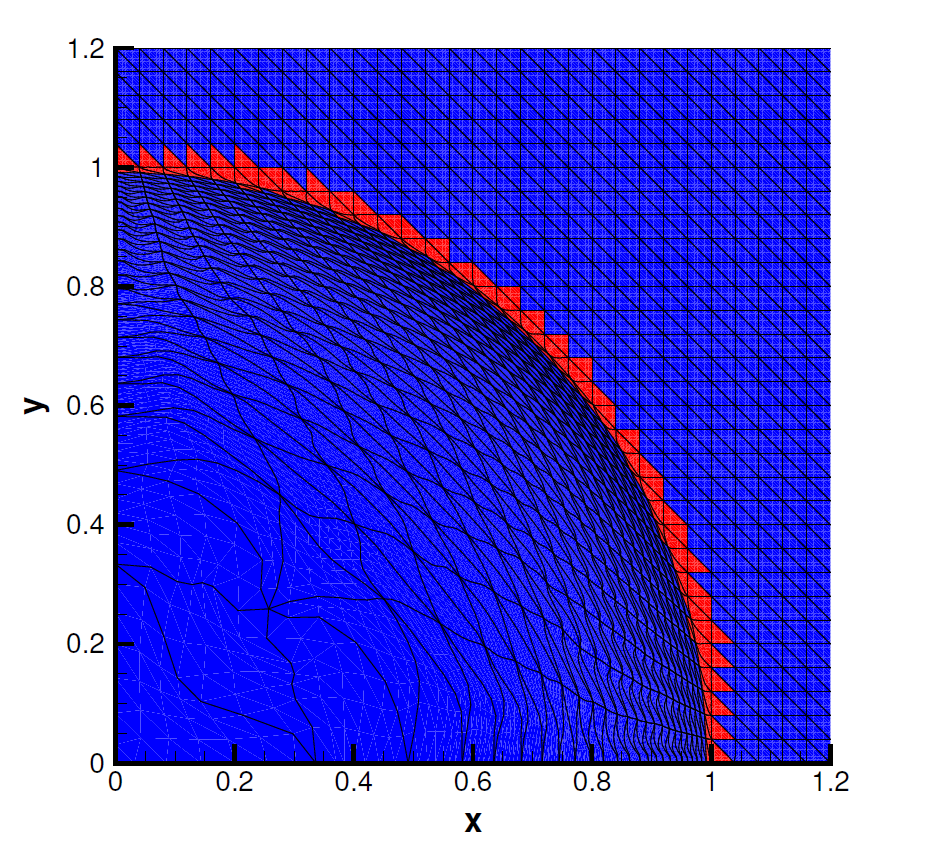}  &           
\includegraphics[width=0.47\textwidth]{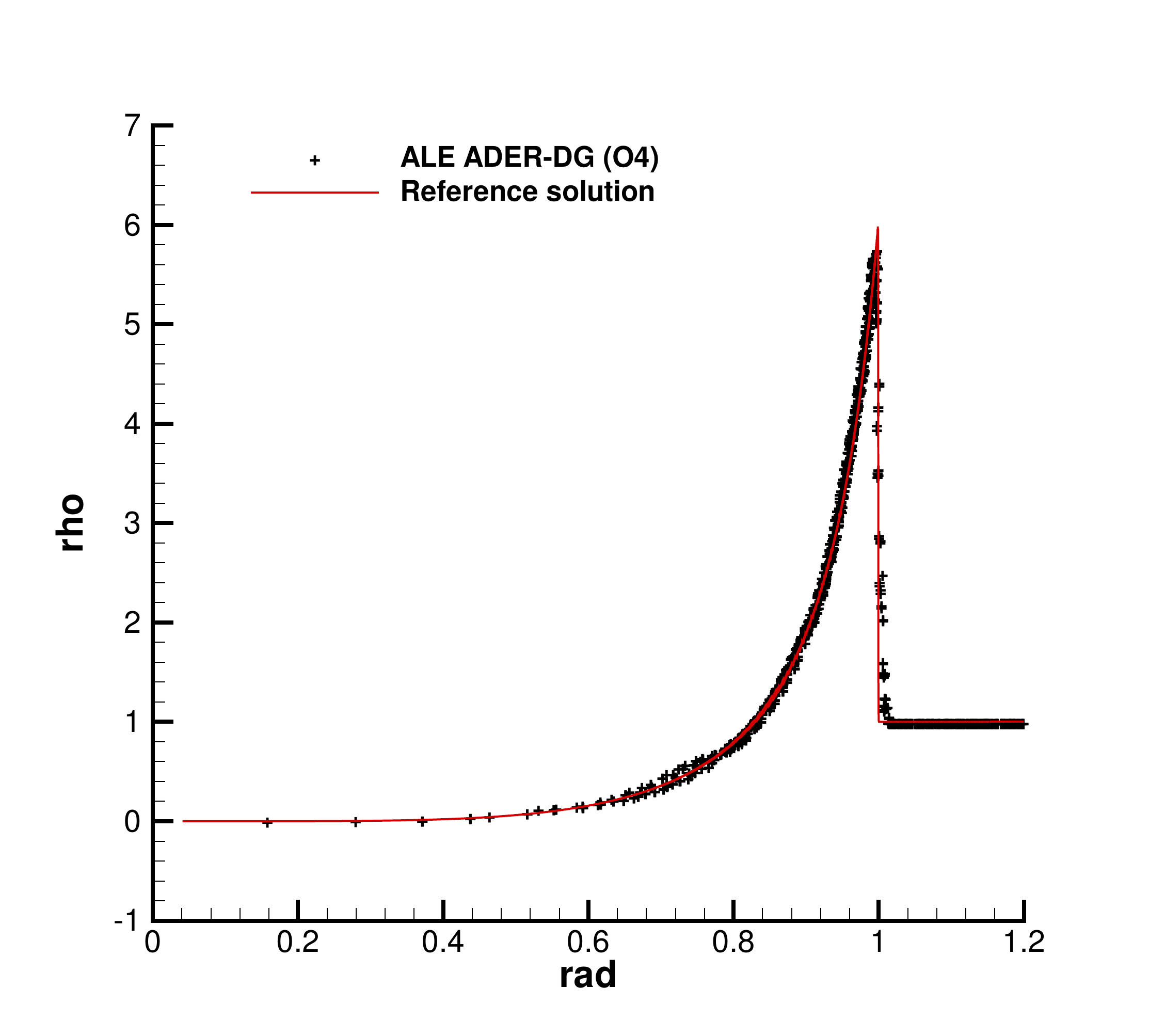} \\ 
\includegraphics[width=0.47\textwidth]{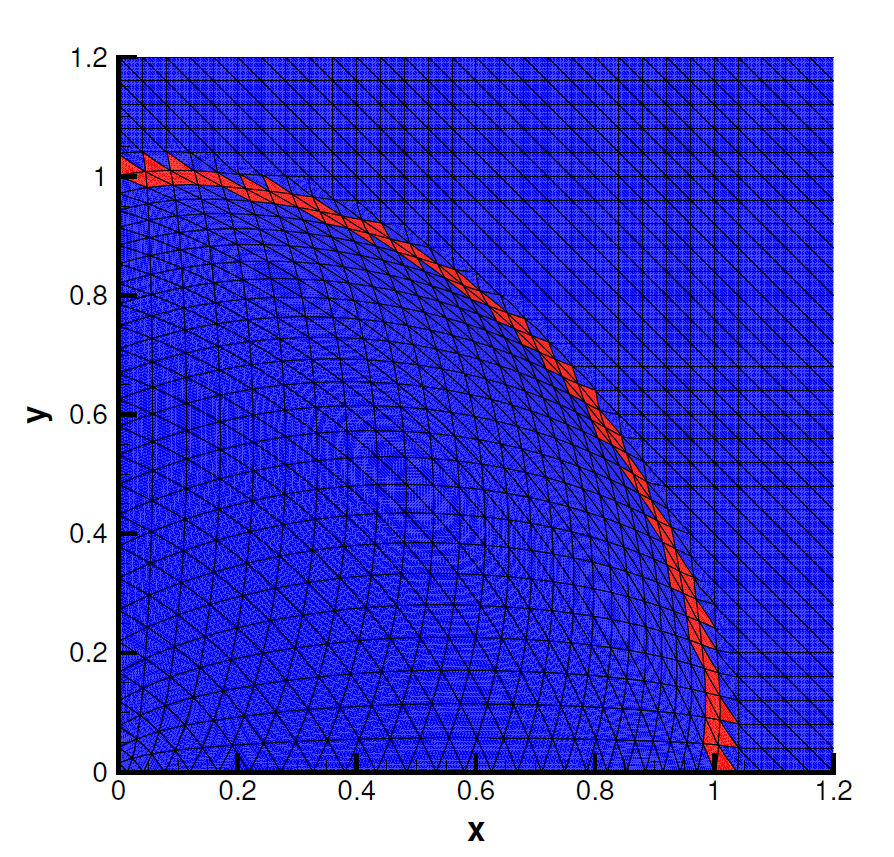}  &           
\includegraphics[width=0.47\textwidth]{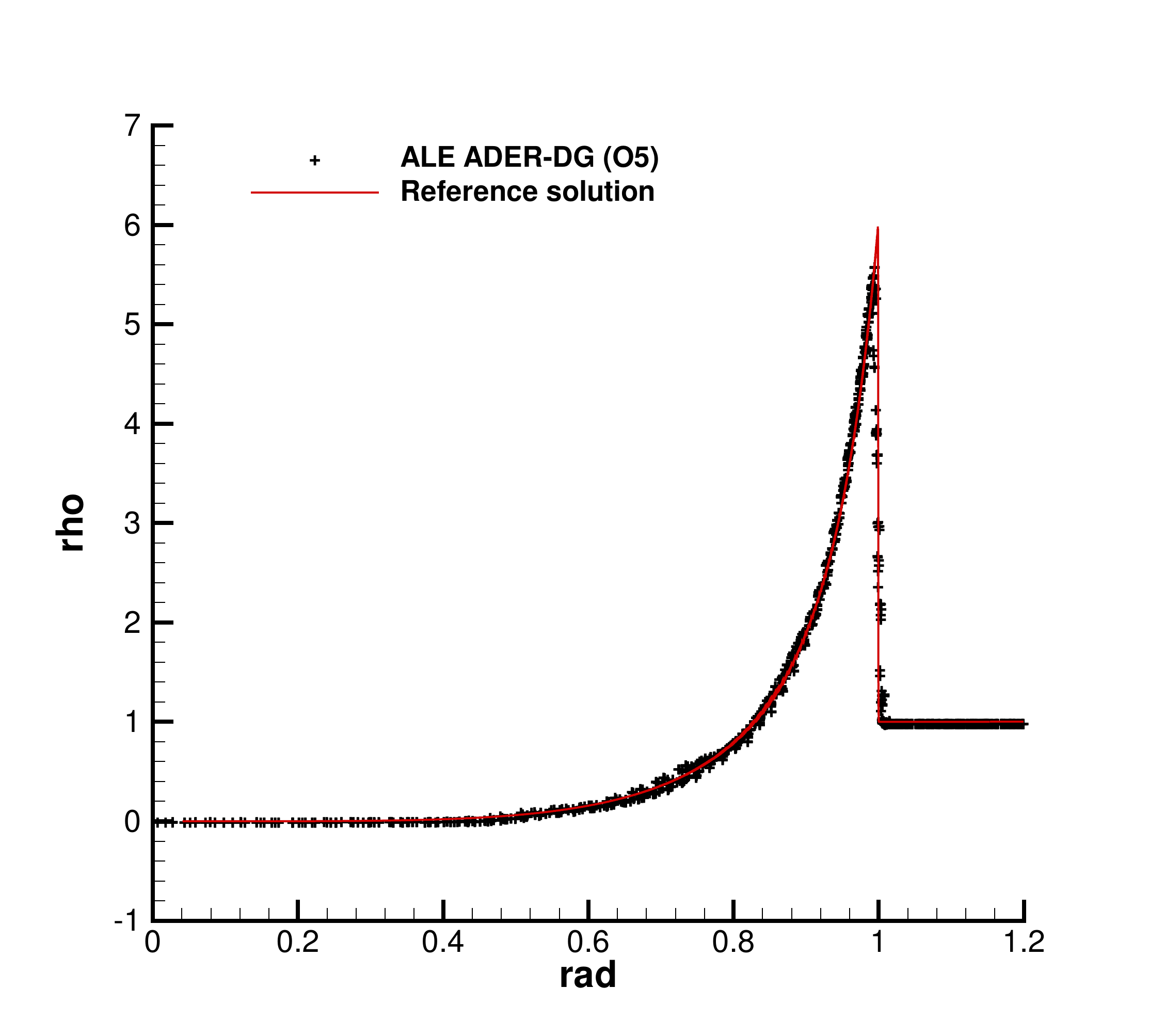} \\           
\end{tabular} 
\caption{Sedov problem. Two-dimensional mesh configuration at the final time $t_f=1$ with sub-cell limiter map (left) and scatter plot of cell density (right). Top: fourth order accurate numerical results. Bottom: fifth order accurate simulation with strong rezoning.} 
\label{fig.Sedov}
\end{center}
\end{figure}

Here, we consider the evolution of a strong shock wave induced by a very high energy deposit, initially located at the origin $\mathbf{O}=(\x)=(0)$ of the computational domain, which is given by $\Omega(0)=[0;1.2]^d$. The mesh is composed by $N_E=30^d$ control volumes, each of them split into two triangles, according to \cite{LagrangeMHD,Lagrange3D}. The Sedov problem constitutes a benchmark in literature \cite{Maire2009,Maire2009b,LoubereSedov3D} since it allows the algorithm to be tested against strong element compressions produced by the diverging shock wave. The initial condition in primitive variable simply reads $\U_0=(1,0,0,0,p_0)$, where the initial pressure is $p_0=10^{-6}$ everywhere except for the cell $c_{or}$ containing the origin of the domain where we assign
\begin{equation}
p_{or} = (\gamma-1)\rho_0 \frac{E_{tot}}{\alpha \cdot V_{or}} \quad \textnormal{ with } \quad 
E_{tot} = 0.244816, 
\label{eqn.p0.sedov}
\end{equation}
where the ratio of specific heats is $\gamma=1.4$ and $E_{tot}$ represents the total energy density. Furthermore, $\alpha$ is a factor which takes into account the cylindrical symmetry, hence becoming $\alpha=4$ for the two-dimensional case. The final time of the simulation is $t_f=1.0$ and the exact solution is a symmetric cylindrical shock wave located at radius $R=\sqrt{\mathbf{x}^2}=1$ with a density peak of $\rho=6$. Figure \ref{fig.Sedov} demonstrates that the fourth order ALE ADER-DG scheme approximates very well the density distribution, although the computational mesh is highly distorted. The sub-cell limiter is active only at the shock front, as expected. To obtain a better quality in the final mesh configuration we have also run the Sedov problem with $N=4$ and a constant relaxation parameter $\omega_k=0.7$ in Eqn. \eqref{eqn.relaxation}, so that the computational mesh is strongly rezoned. The corresponding results are depicted in the bottom panels of Figure \ref{fig.Sedov}.

\subsection{Viscous shock problem} 
\label{sec.VS}
\begin{figure}[!htbp]
\begin{center}
\begin{tabular}{cc} 
\includegraphics[width=0.47\textwidth]{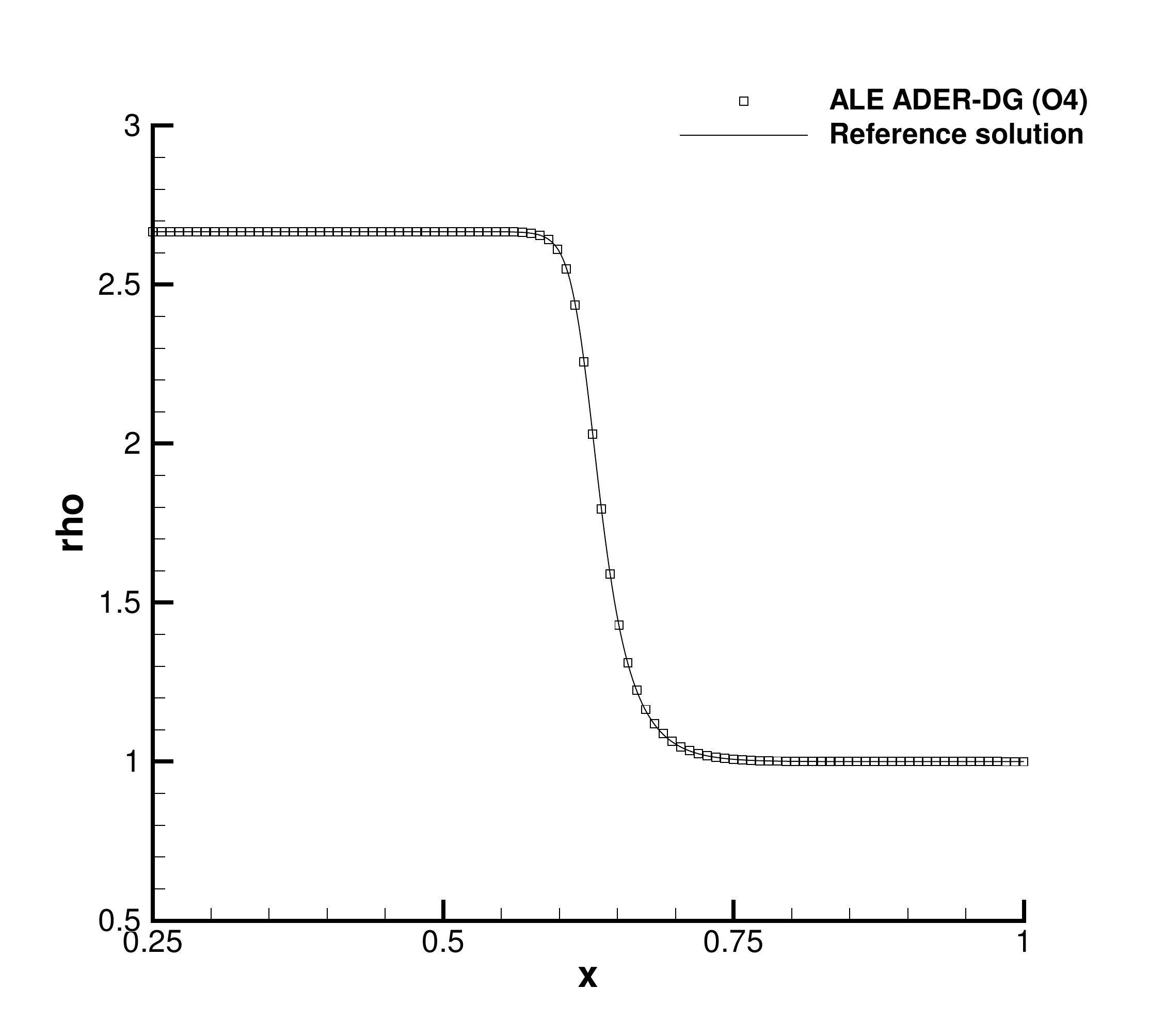}  &           
\includegraphics[width=0.47\textwidth]{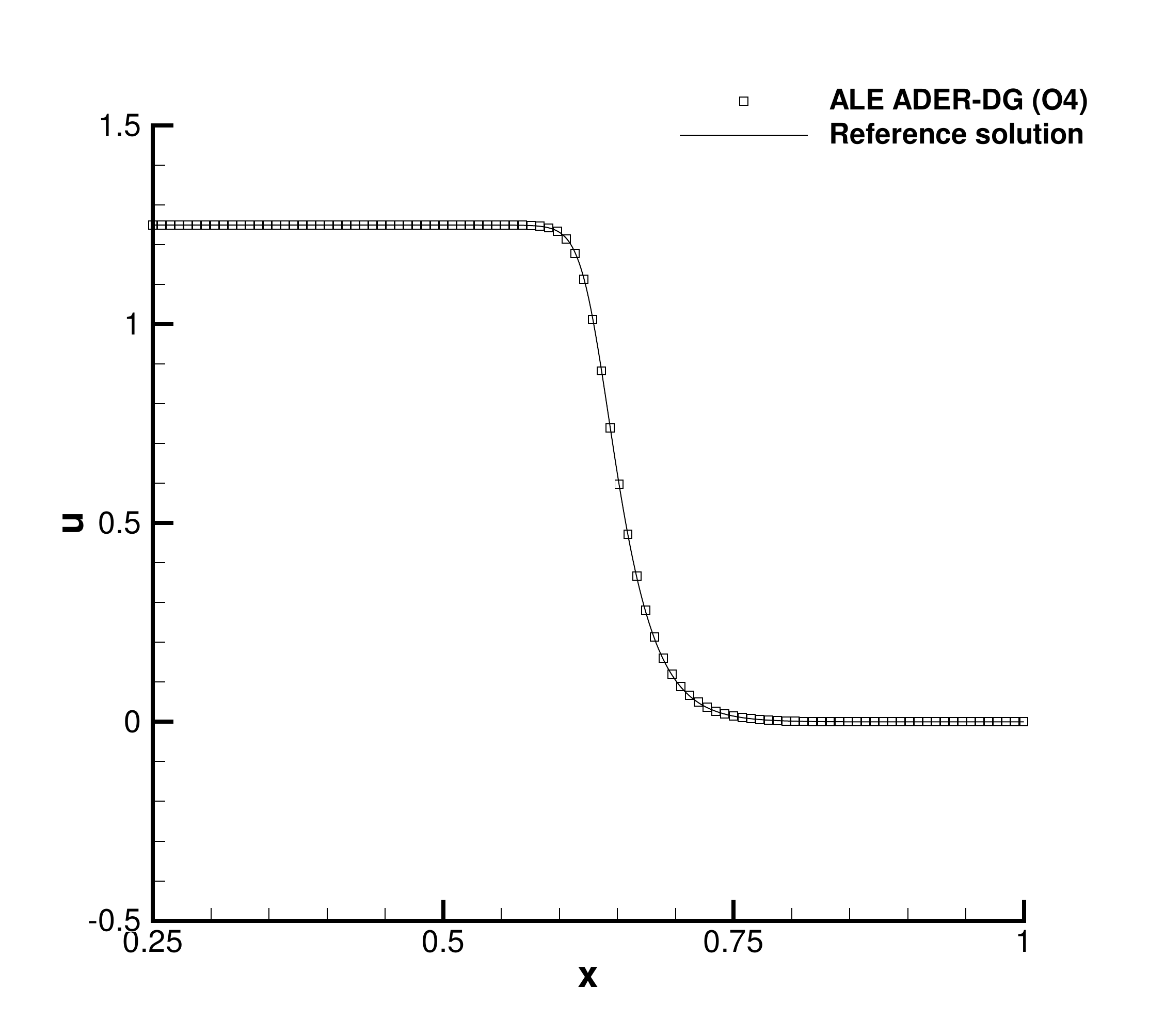} \\ 
\includegraphics[width=0.47\textwidth]{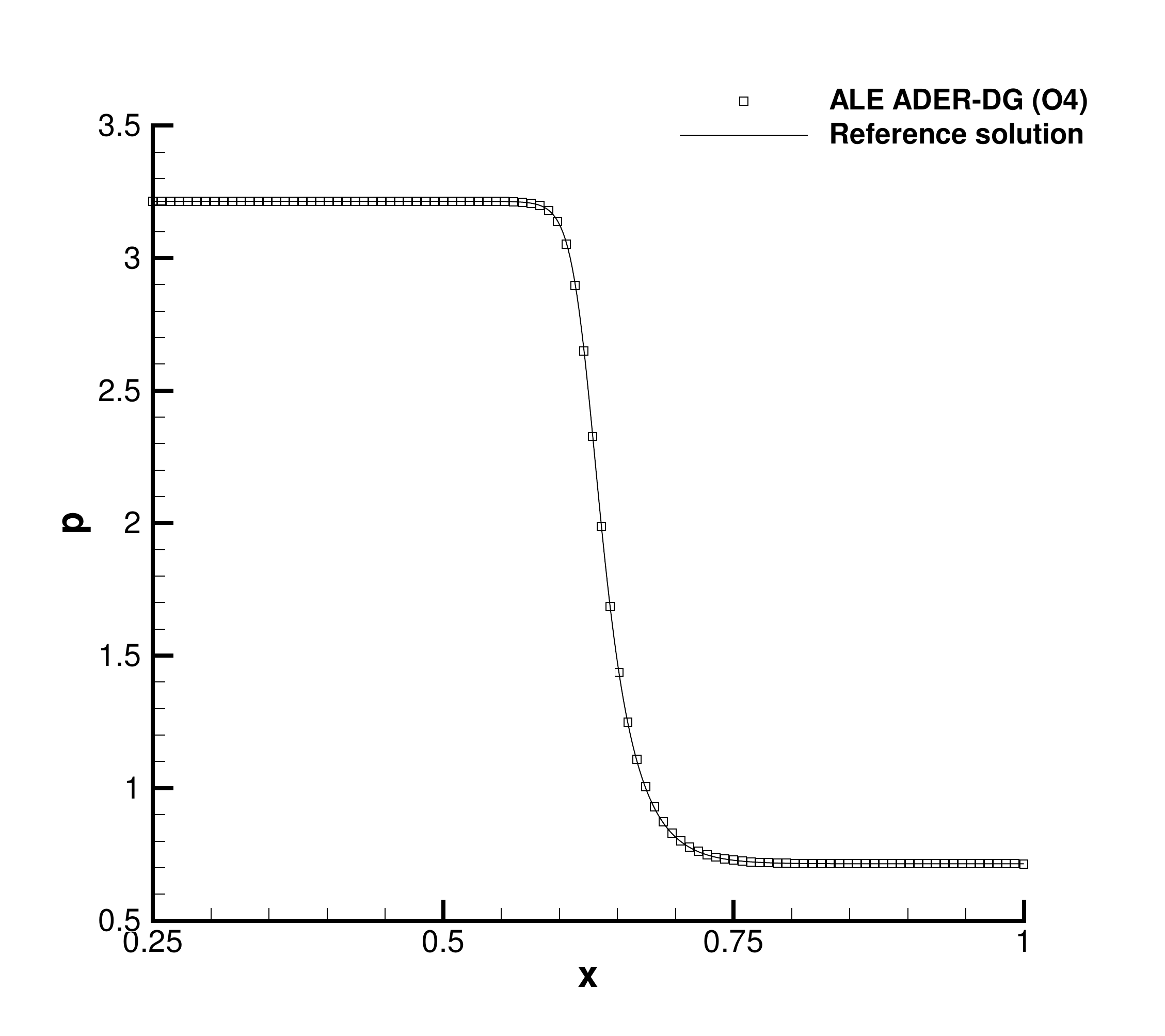}  &           
\includegraphics[width=0.47\textwidth]{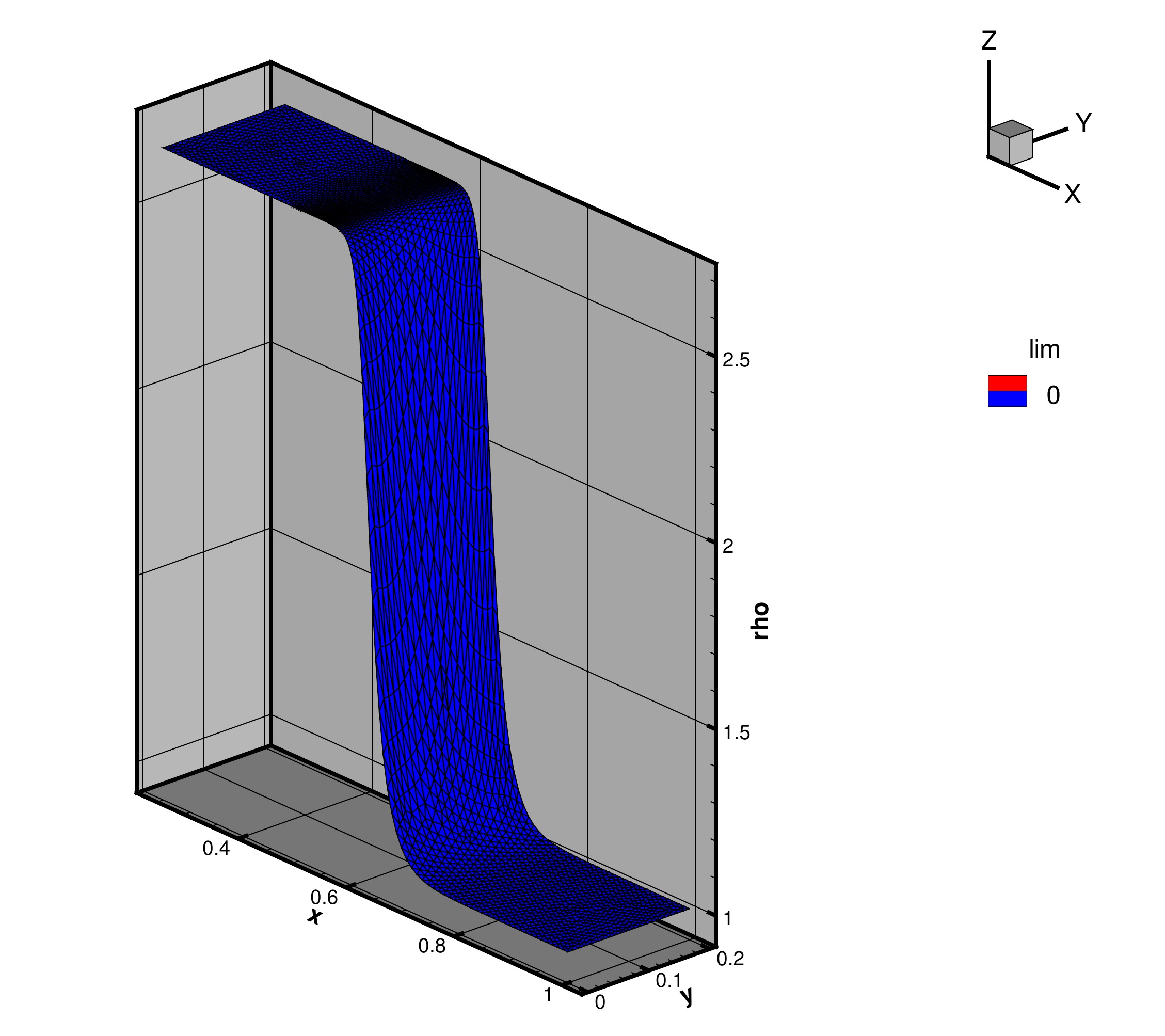} \\           
\end{tabular} 
\caption{Viscous shock problem at final time $t=0.2$. We show the comparison between numerical and analytical solution for density, velocity and pressure as well as a three-dimensional view of the density distribution with the corresponding sub-cell limiter map. } 
\label{fig.VS}
\end{center}
\end{figure}

\begin{figure}[!htbp]
\begin{center}
\begin{tabular}{c} 
\includegraphics[width=0.7\textwidth]{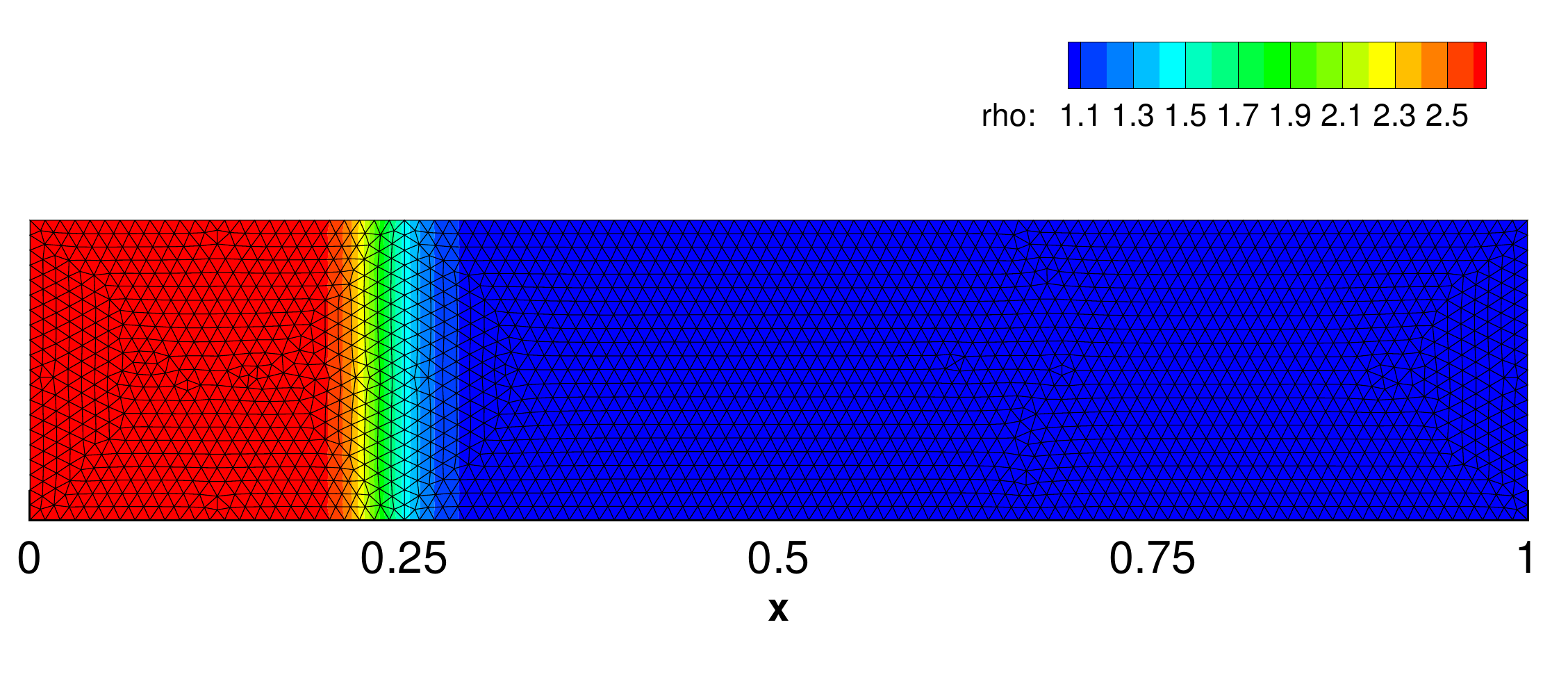} \\           
\includegraphics[width=0.7\textwidth]{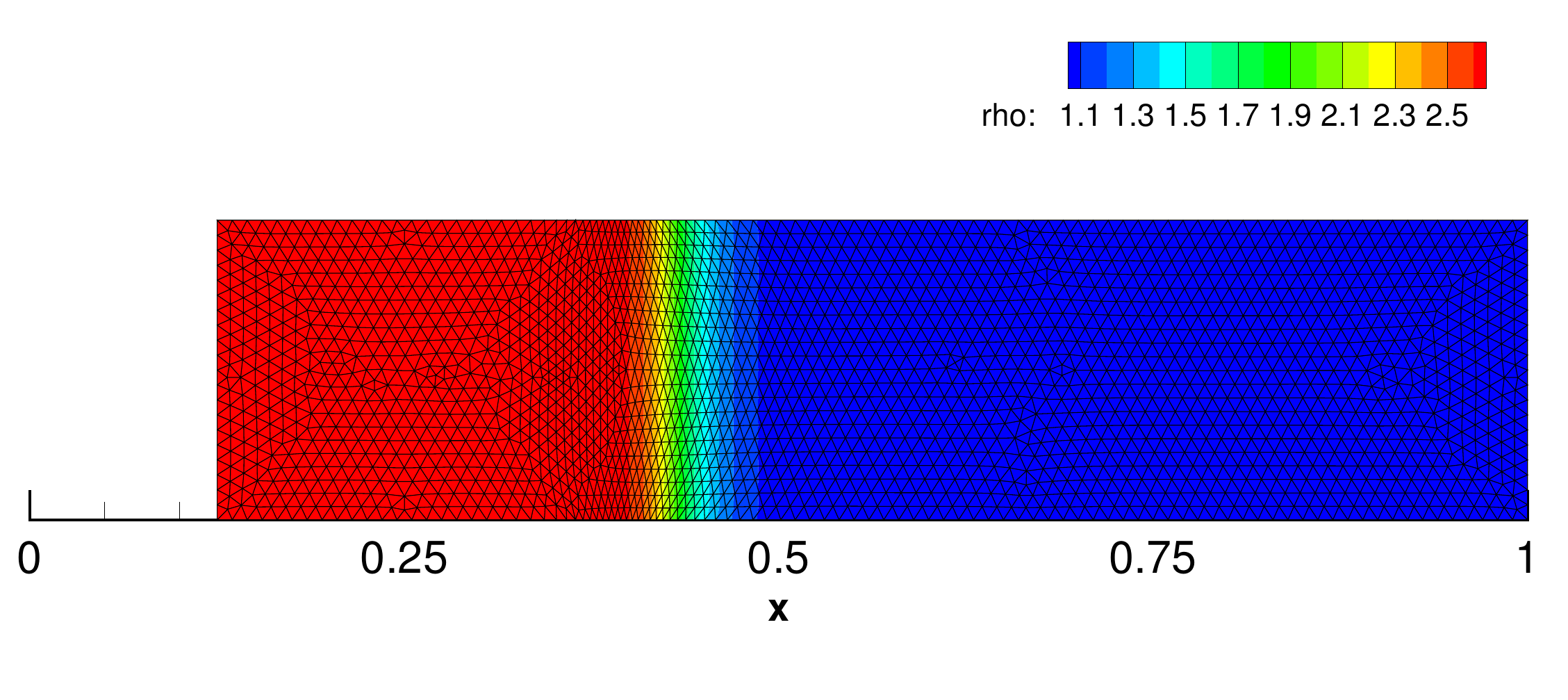} \\ 
\includegraphics[width=0.7\textwidth]{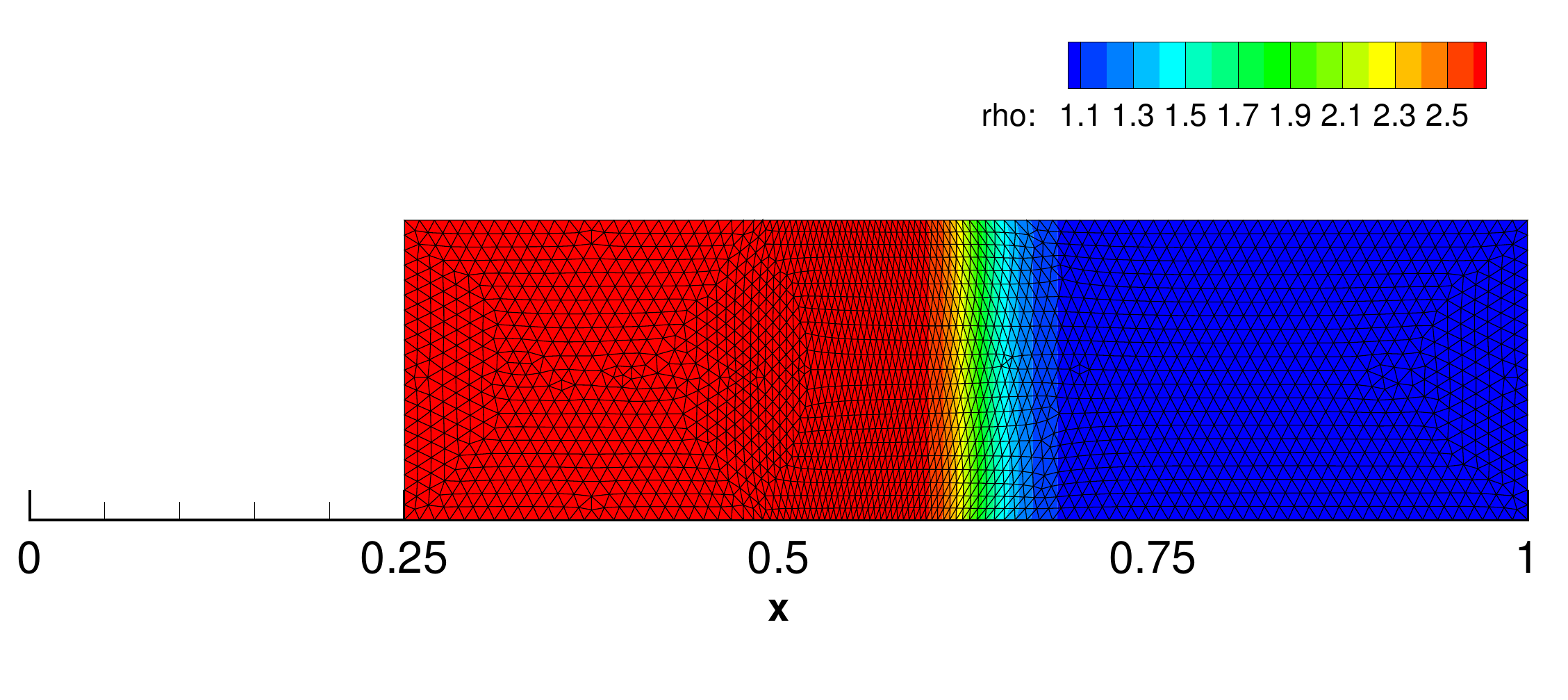} \\          
\end{tabular} 
\caption{Density distribution and mesh configuration for the viscous shock problem at output times $t=0.0$, $t=0.1$ and $t=0.2$. } 
\label{fig.VSevolution}
\end{center}
\end{figure}

The test problem described in the following is concerned with \textit{physical viscosity}, therefore in this case we are solving the compressible Navier-Stokes equations \eqref{eqn.NSflux}. Specifically, an isolated viscous shock wave is propagating into a fluid at rest in the supersonic regime, i.e. with a shock Mach number of $M_s > 1$. The setup of this test problem starts from the analytical solution of the compressible Navier-Stokes equations derived in \cite{Becker1923} for the particular case of a stationary shock wave at Prandtl number $Pr= 0.75$ with constant viscosity. According to \cite{LagrangeHPR,Dumbser_HPR_16}, a constant velocity field $u = M_s c_0$ is superimposed to the previous stationary shock wave solution, hence obtaining a \textit{non-stationary} shock wave initially traveling at $M_s=2$ with a Reynolds number of $Re=100$. The fluid before the shock is assigned with constant density $\rho_0=1$, velocity $u_0=1.25$ and pressure $p_0=1/\gamma$ with $\gamma=1.4$. The physical viscosity is $\mu=2\times 10^{-2}$ and the final time of the simulation is chosen to be $t_f=0.2$. The initial computational domain is the rectangular box $\Omega(0)=[0;1]\times[0;0.2]$, which is discretized by an unstructured computational mesh with characteristic mesh size $h=1/100$, yielding a total number of $N_E=4462$ triangles. Periodic boundaries are imposed in the $y-$direction, a no-slip wall is placed at $x=1$ while the left side of the domain is moved with the local fluid velocity. The shock wave is initially centered at $x=0.25$ and we use the fourth order version of our ALE ADER-DG schemes to run the simulation with the isoparametric approach for the mesh motion. Figure \ref{fig.VS} shows a comparison of the numerical results against the analytical solution, where an excellent matching can be appreciated. Furthermore, the sub-cell limiter is correctly not active in the whole computational domain as expected, since the solution does not involve any discontinuity (the shock structure is fully resolved here). Finally, Figure \ref{fig.VSevolution} depicts the density distribution as well as the mesh configurations at output times $t=0.0$, $t=0.1$ and $t=0.2$.     

\subsection{Taylor-Green vortex} 
\label{sec.TGV}
\begin{figure}[!htbp]
\begin{center}
\begin{tabular}{cc} 
\includegraphics[width=0.47\textwidth]{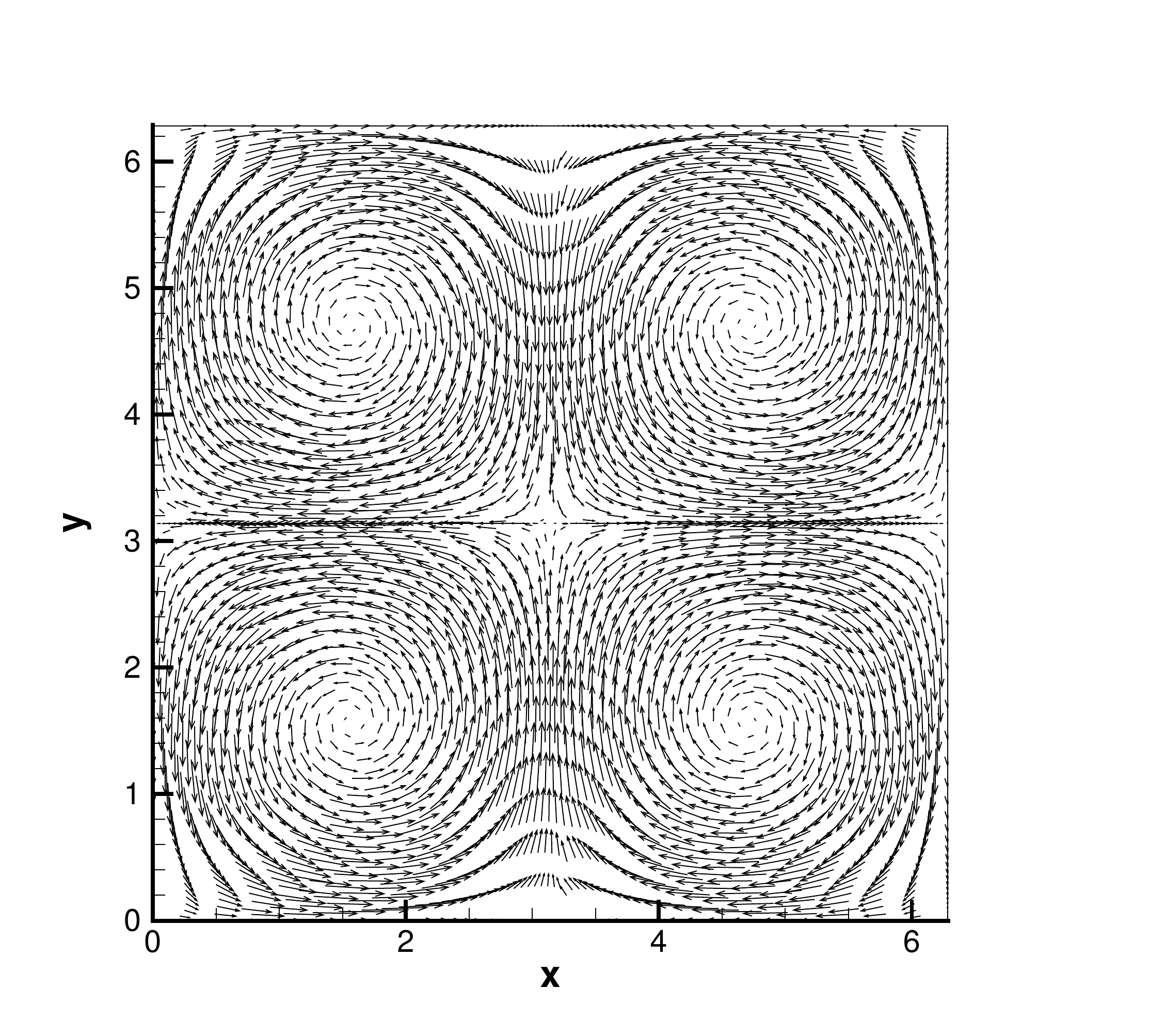}  &           
\includegraphics[width=0.47\textwidth]{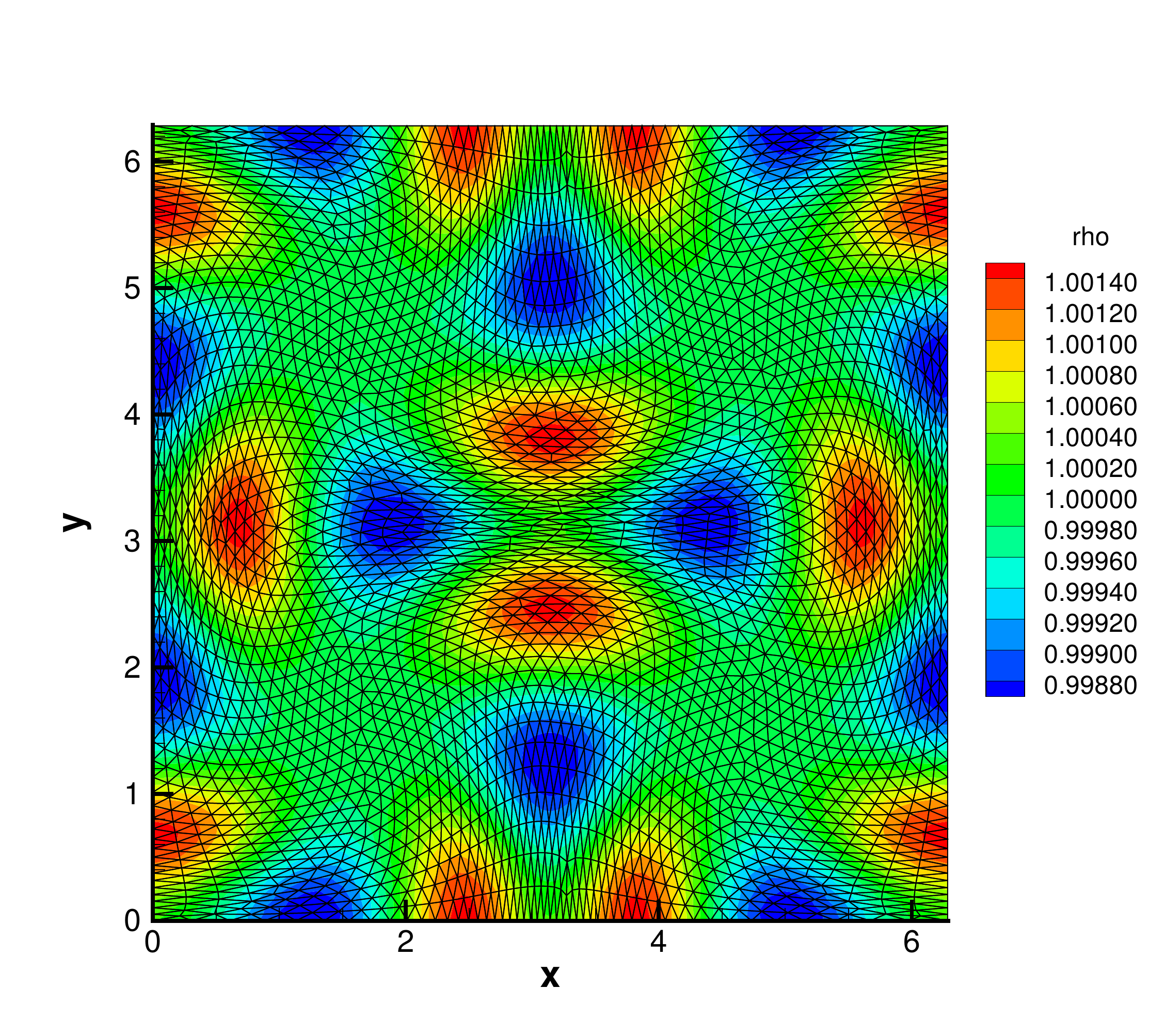} \\ 
\includegraphics[width=0.47\textwidth]{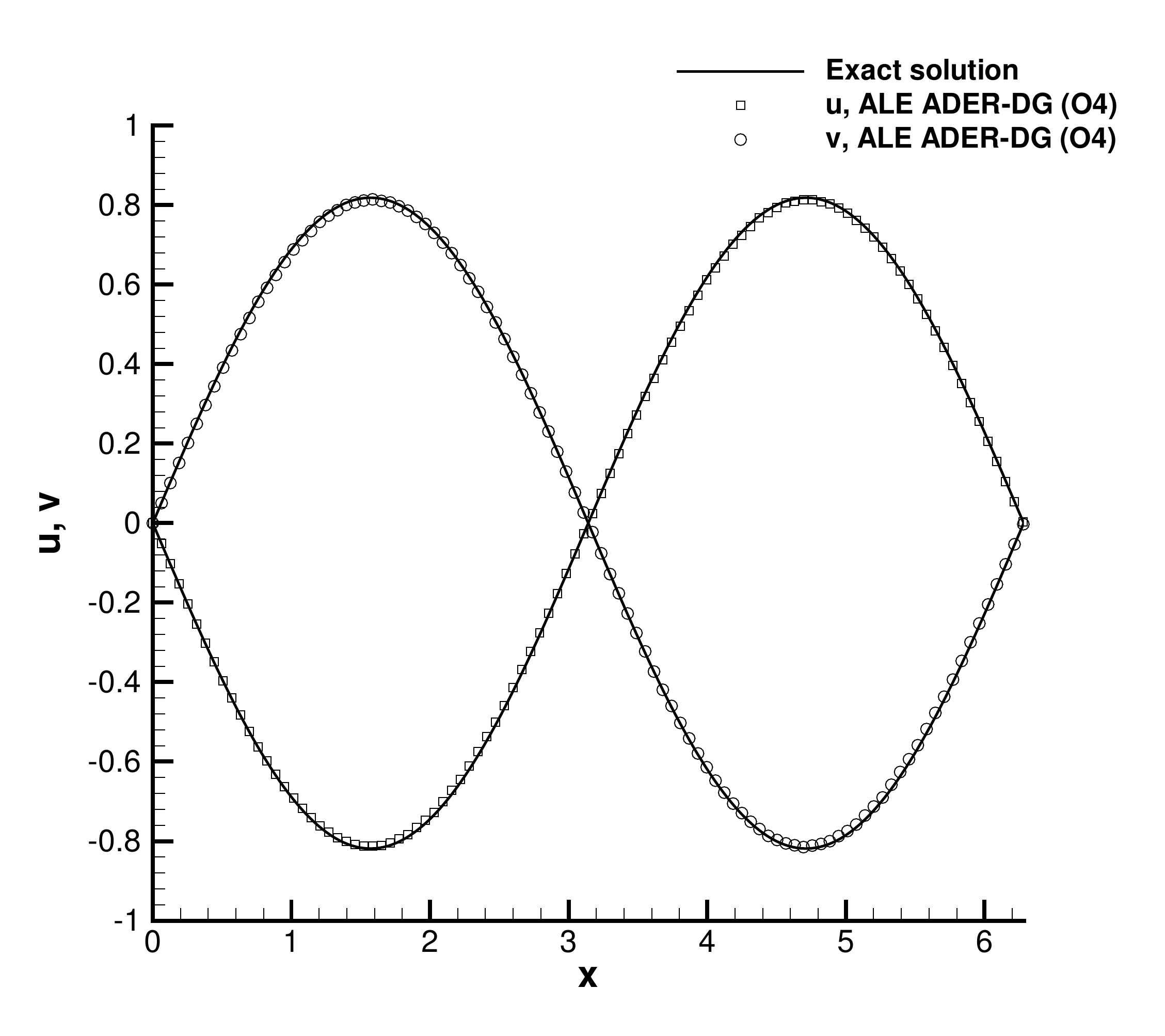}  &           
\includegraphics[width=0.47\textwidth]{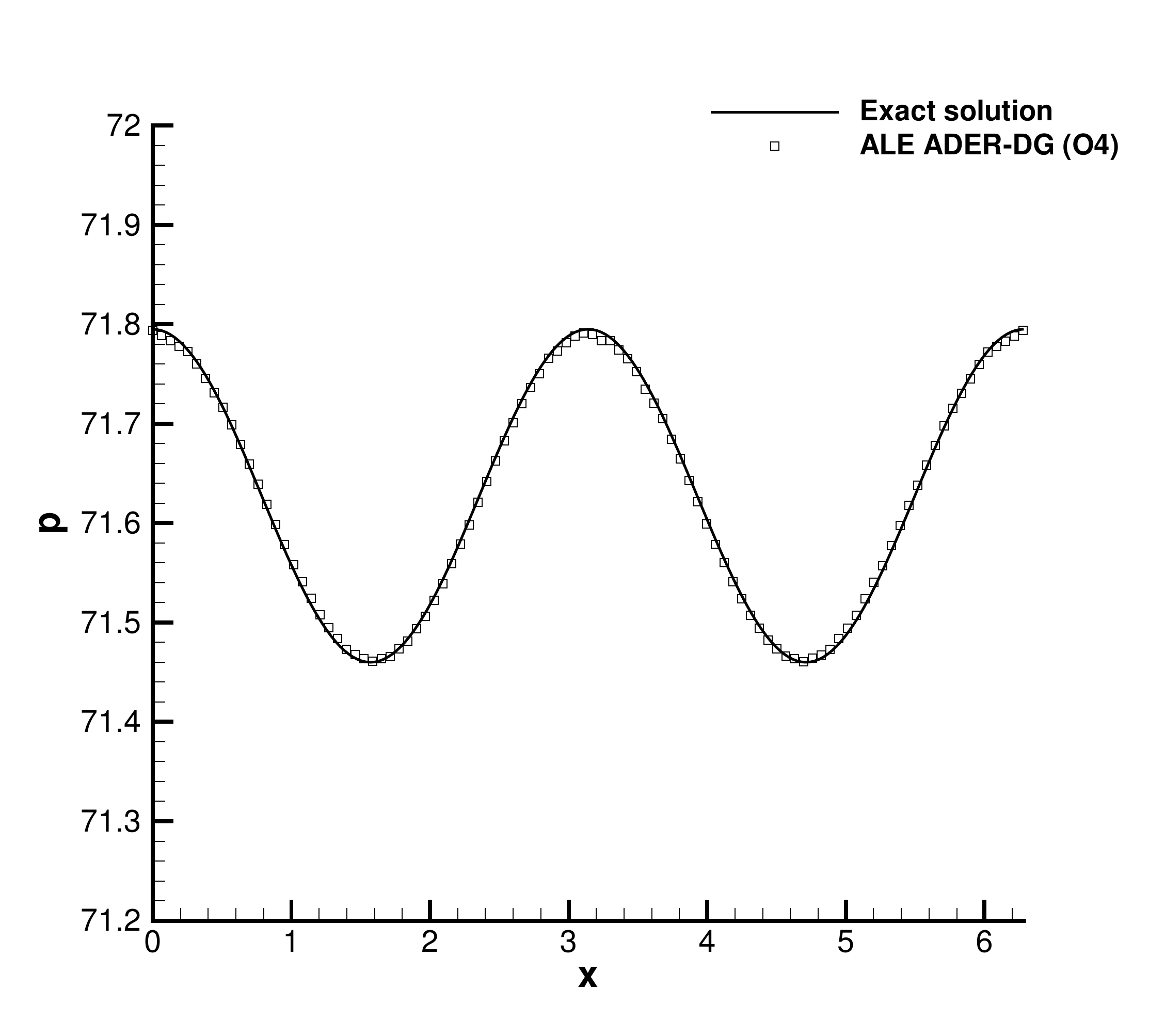} \\           
\end{tabular} 
\caption{Taylor-Green vortex with physical viscosity $\mu=10^{-1}$ at final time $t_f=1.0$. Top: velocity vector field (left) and mesh configuration with density distribution (right). Bottom: comparison between exact and numerical solution for the horizontal velocity components $(u,v)$ (left) and for pressure (right).} 
\label{fig.TGV}
\end{center}
\end{figure}

We solve the Taylor-Green vortex problem on the two-dimensional computational domain $\Omega(0)=[0,2\pi]^2$, where periodic boundaries are set everywhere. The final time of the simulation is $t_f=1.0$ and the mesh is composed by a total number of $N_E=5630$ triangles with characteristic mesh size of $h=2\pi/50$. An exact solution is available solving analytically the two-dimensional incompressible Navier-Stokes equations and it reads
\begin{eqnarray}
    \rho(x,y,t)&=&\rho_0, \nonumber \\
    u(x,y,t)&=&\sin(x)\cos(y)e^{-2\nu t},  \nonumber \\
		v(x,y,t)&=&-\cos(x)\sin(y)e^{-2\nu t}, \nonumber \\
    p(x,y,t)&=& C + \frac{1}{4}(\cos(2x)+\cos(2y))e^{-4\nu t},
\label{eq:TG_ini}
\end{eqnarray}
with the kinematic viscosity $\nu=\frac{\mu}{\rho}$, the density $\rho_0=1$, the ratio of specific heats $\gamma=1.4$ and the initial additive constant for the pressure field $C=100/\gamma$. The analytical solution gives also the initial condition and the physical viscosity is chosen to be $\mu=10^{-1}$. We use a fourth order ALE ADER-DG scheme to carry out the numerical simulation solving the compressible Navier-Stokes equations \eqref{eqn.NSflux}, where the mesh motion is again driven relying on the isoparametric version of our algorithm. Results are depicted in Figure \ref{fig.TGV} and compared against the exact solution, showing an excellent agreement both for velocity and pressure. Since no discontinuities are involved in this test problem, the sub-cell limiter is not active in any cell.

\subsection{Spherical Implosion} 
\label{sec.ICF}
\begin{figure}[!htbp]
\begin{center}
\begin{tabular}{cc} 
\includegraphics[width=0.4\textwidth]{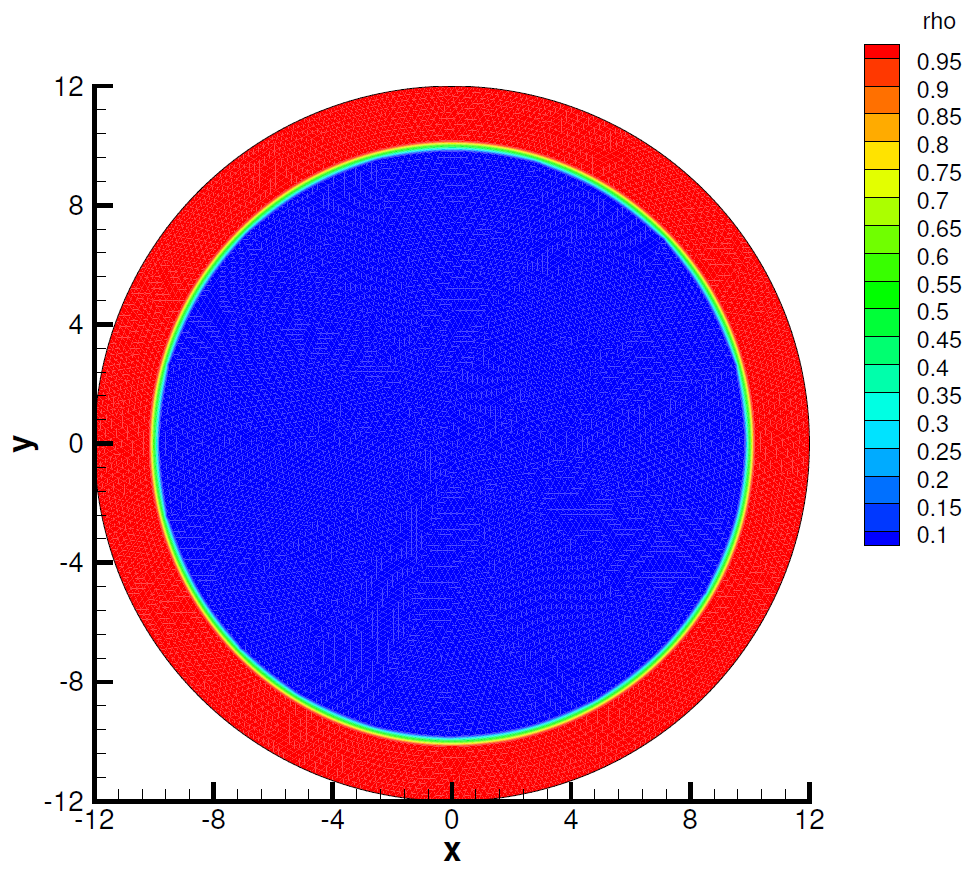} &           
\includegraphics[width=0.4\textwidth]{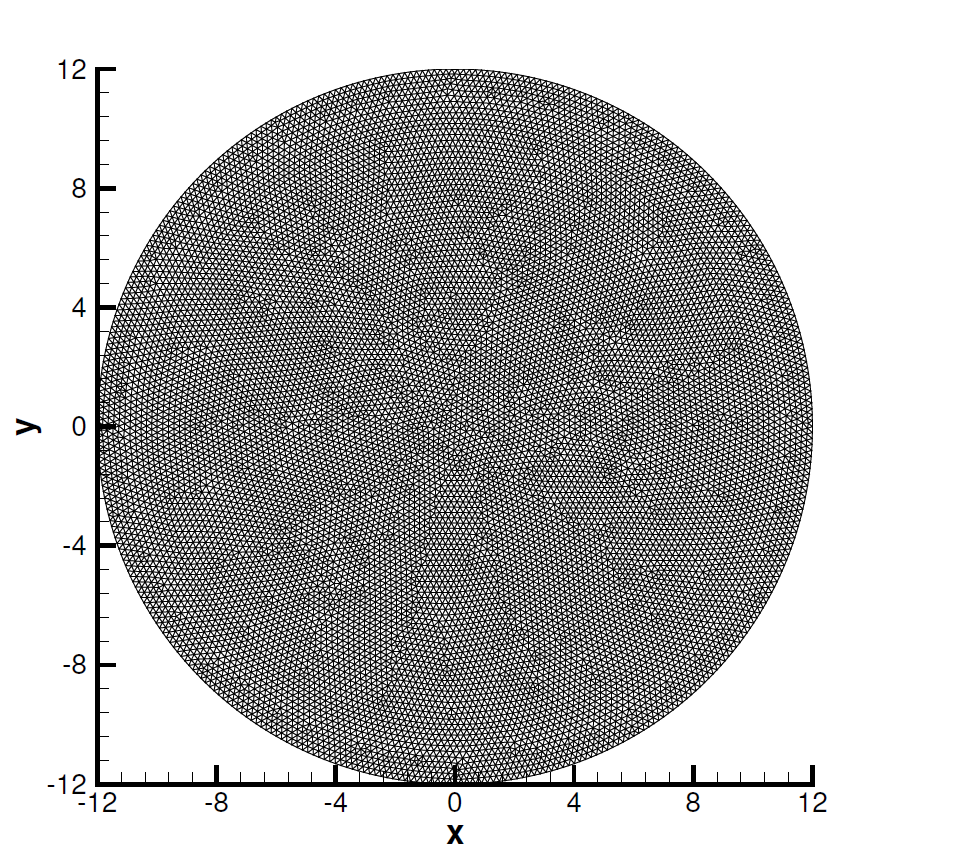} \\
\includegraphics[width=0.4\textwidth]{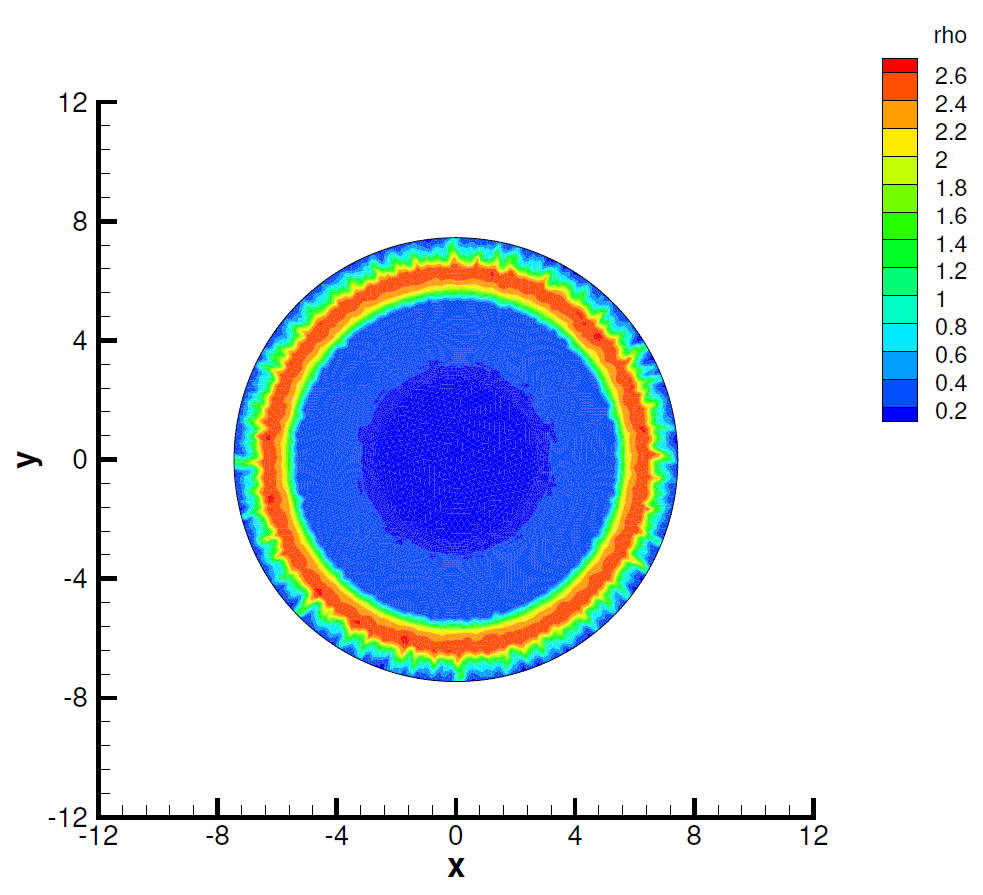} &
\includegraphics[width=0.4\textwidth]{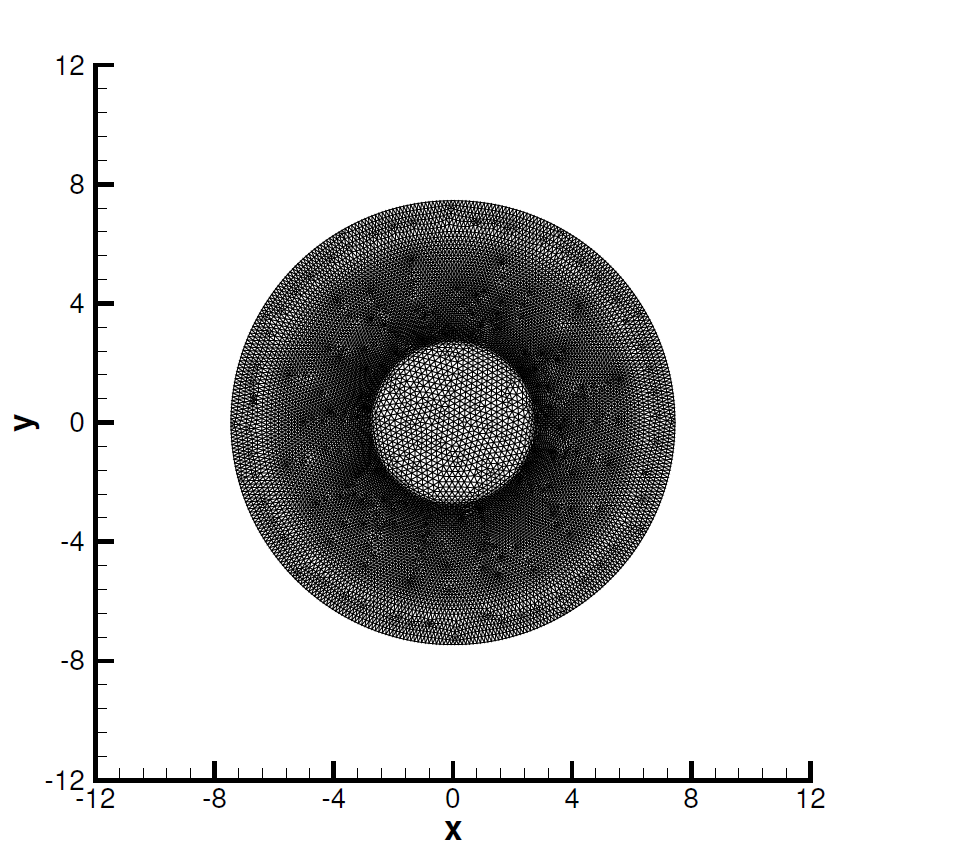} \\           
\includegraphics[width=0.4\textwidth]{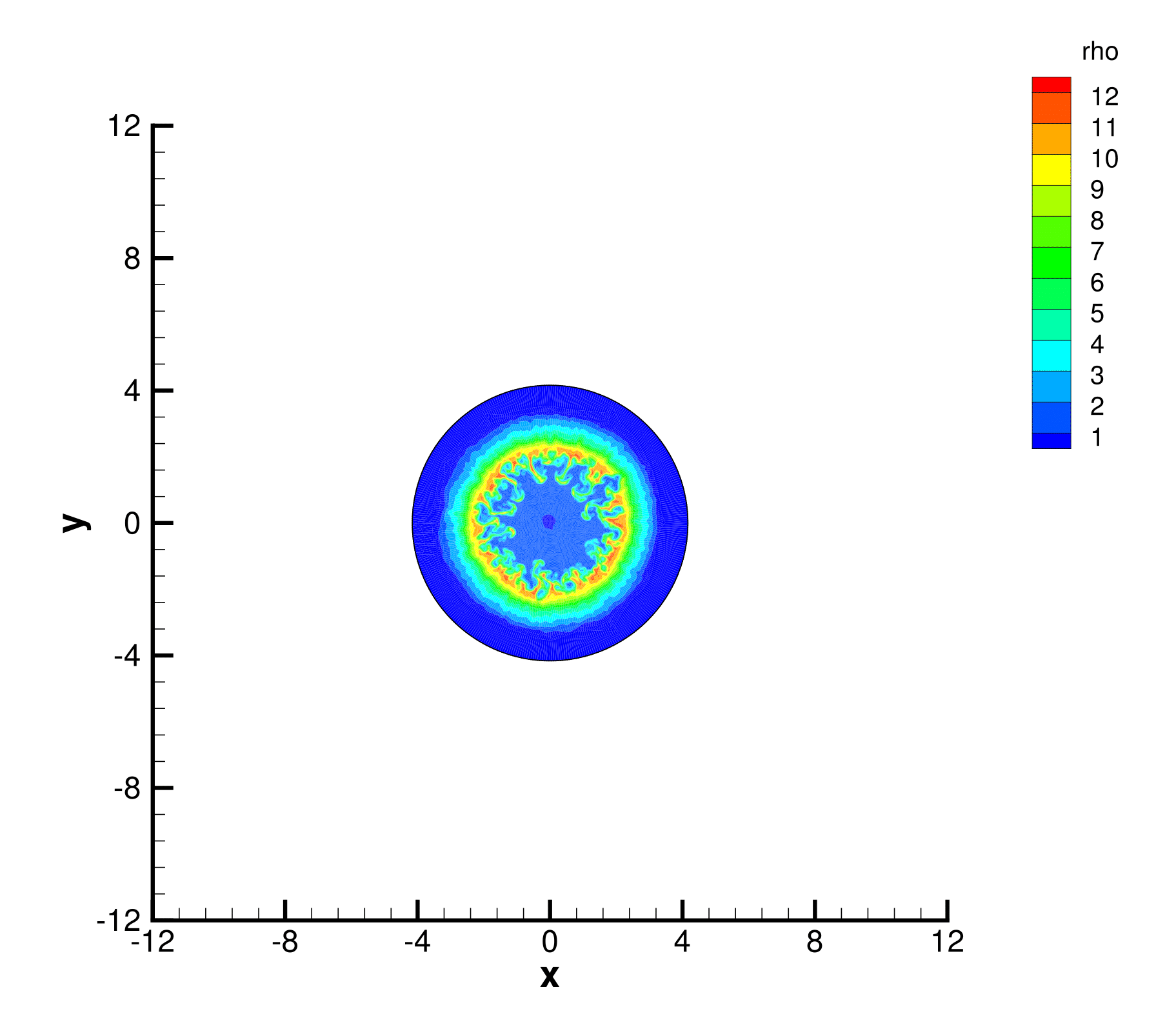} &
\includegraphics[width=0.4\textwidth]{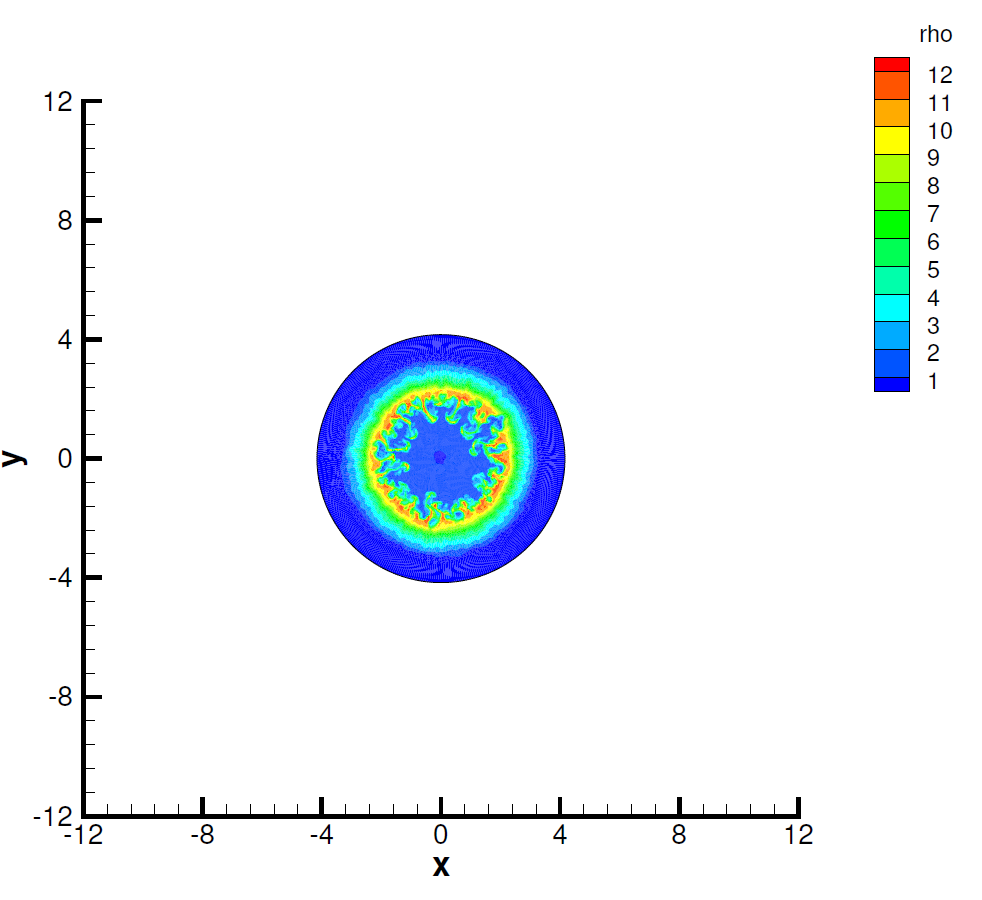} \\          
\end{tabular} 
\caption{Density distribution (left column) and mesh configuration (right column) for the spherical implosion problem at output times $t=0.0$, $t=1.5$ and the final time $t_f=2.77$ with the external radius located at $r_e=4.0$. } 
\label{fig.ICF-B0}
\end{center}
\end{figure}

The last test case describes an implosion, which is quite similar to what happens in Ignition Confinement Fusion (ICF) simulations. The initial computational domain is given by the circle of radius $R=12$, that is split into an internal and an external region at radius $R_s=10$. The inner zone is filled by a light gas with initial density and pressure $(\rho_l,p_l)=(0.05,0.1)$, while in the outer shell there is a heavy fluid with $(\rho_h,p_h)=(1.0,0.1)$. Both fluids are initially at rest and the ratio of specific heats is set to $\gamma_l=\gamma_h=5/3$. On the external boundary we impose the pressure
\begin{equation}
p^*(t) = \left\{ \begin{array}{lcc} 10,       & \textnormal{ if } & t \in [0,0.5] \\ 
                                    12 - 4t,  & \textnormal{ if } & t > 0.5         
                      \end{array}  \right. ,
\label{eqn.p_icf}
\end{equation}
which drives the implosion. Initially, the shell is collapsing towards the center of the domain, while after $t\approx 2.5$ the pressure of the highly compressed light fluid becomes bigger than the one imposed externally, hence leading to an expansion of the shell. The final time of the simulation is chosen in such a way that the external radius is located at $r_e=4$ with the generic radial position $r=\sqrt{x^2+y^2}$. Figure \ref{fig.ICF-B0} shows the density distribution as well as the mesh configuration at output times $t=0.0$, $t=1.5$ and $t_f=2.77$ obtained running a fourth order direct ALE ADER-DG scheme.

\begin{figure}[!htbp]
\begin{center}
\begin{tabular}{cc} 
\includegraphics[width=0.47\textwidth]{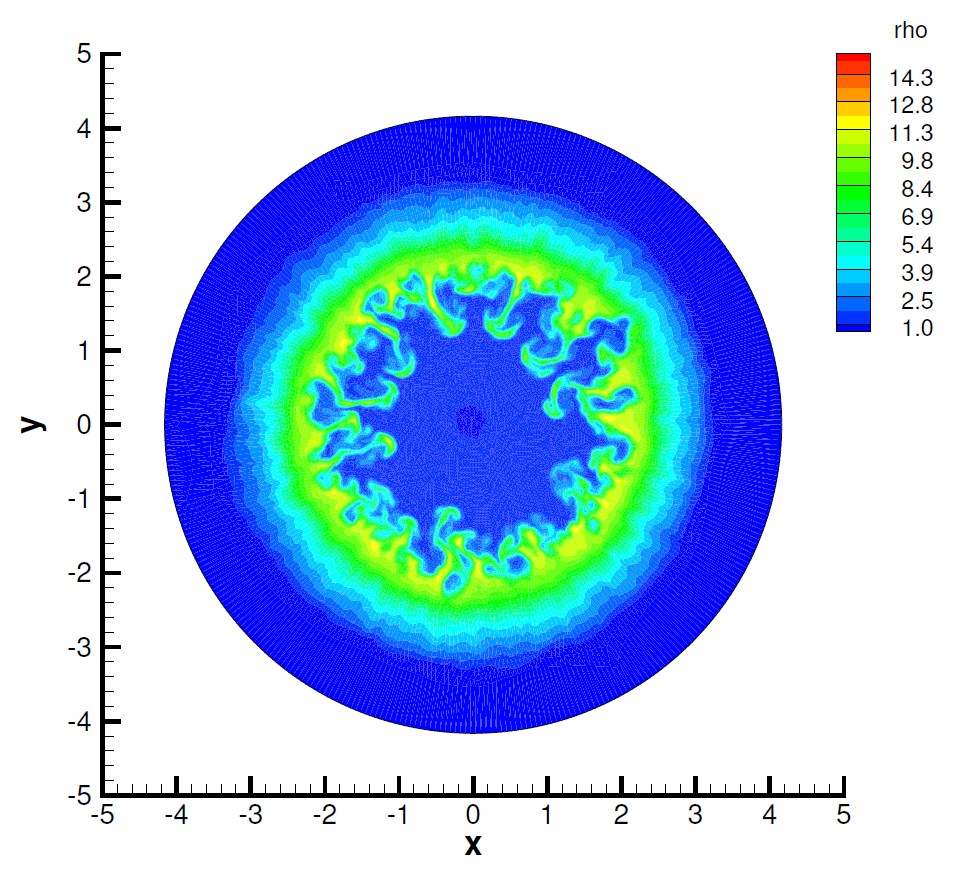} &           
\includegraphics[width=0.47\textwidth]{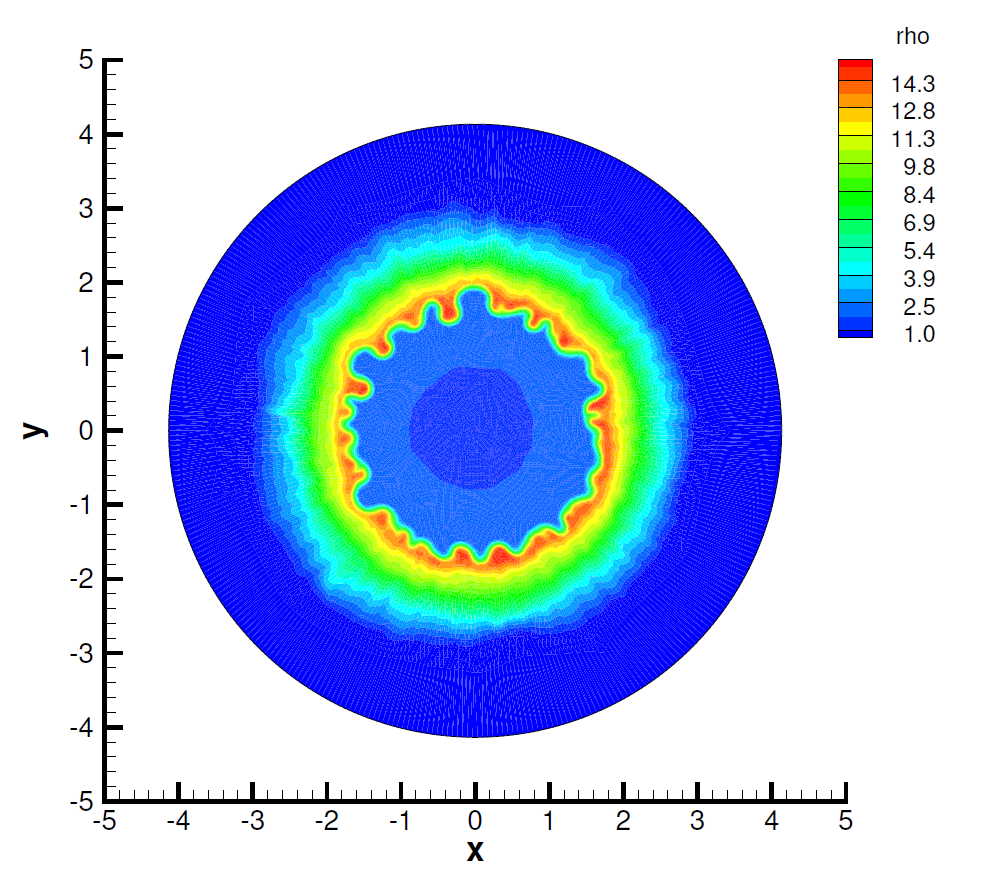} \\
\includegraphics[width=0.47\textwidth]{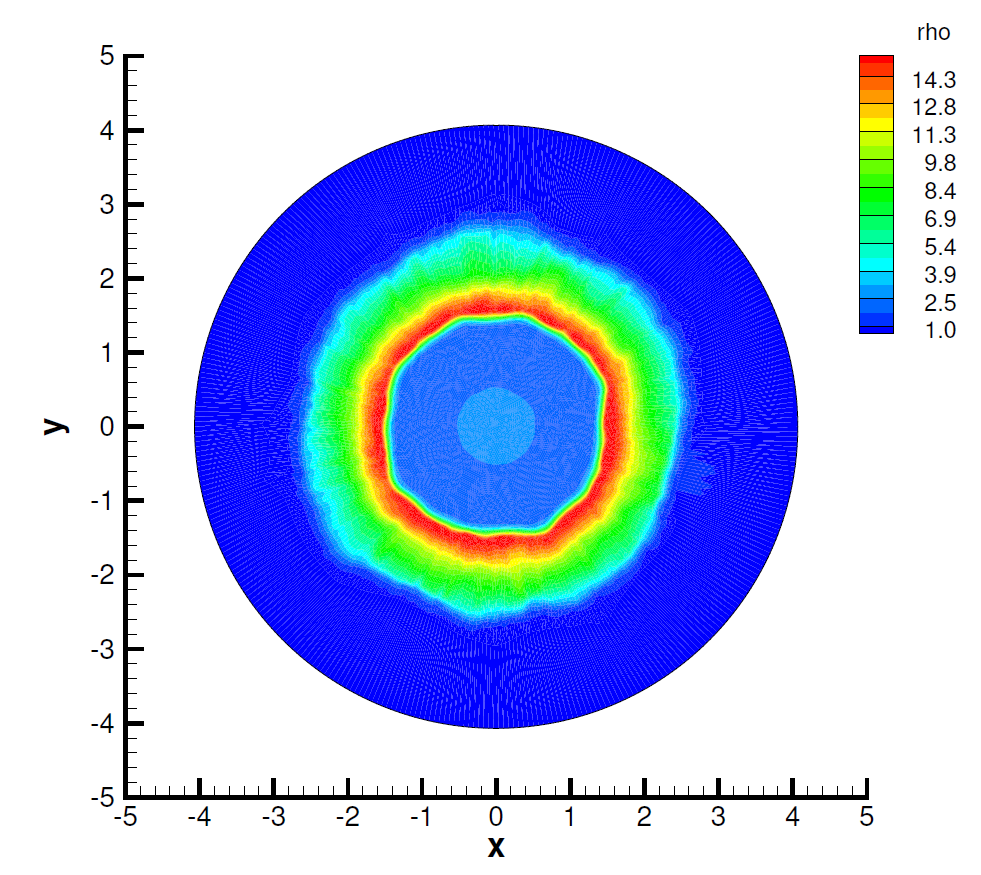} &
\includegraphics[width=0.47\textwidth]{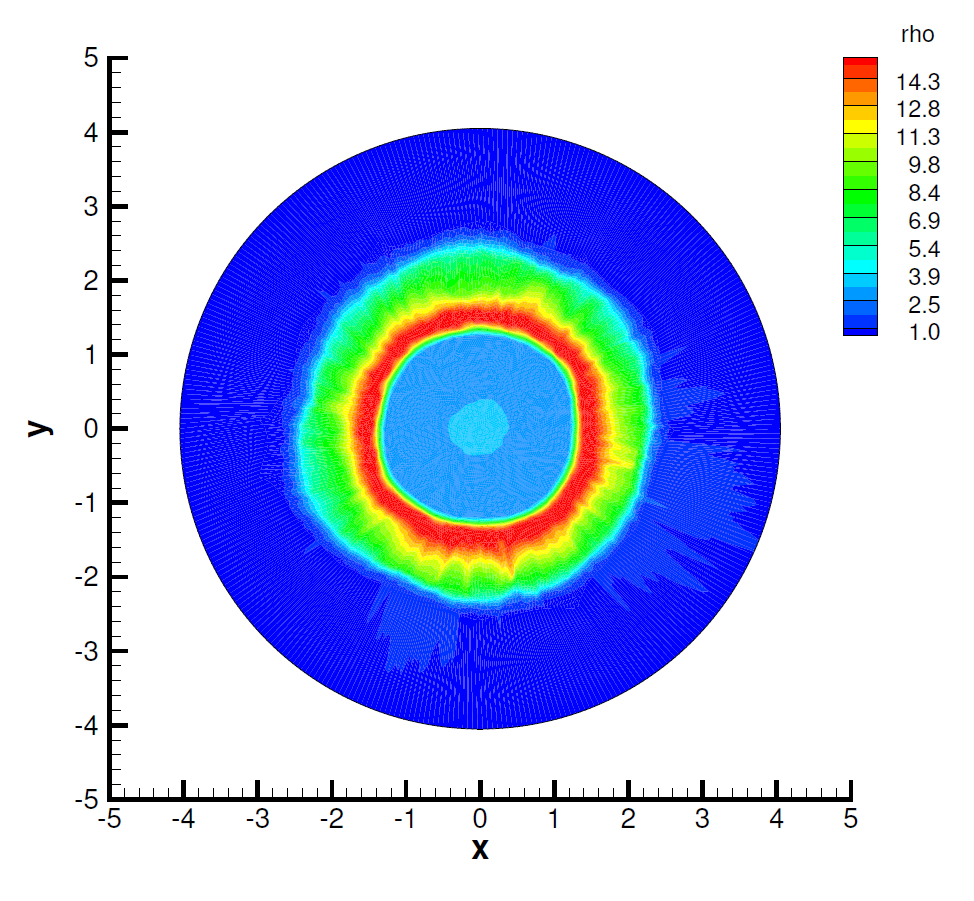} \\           
\end{tabular} 
\caption{Density distribution for the spherical implosion problem solved using the ideal MHD equations with a magnetic field of intensity $B_0$ applied on the horizontal plane $x-y$. From top left to bottom right: $B_0=0$ at output time $t_f=2.77$, $B_0=1$ at output time $t_{f,1}=2.70$, $B_0=2$ at output time $t_{f,2}=2.63$ and $B_0=3$ at output time $t_{f,3}=2.55$. All simulations stop when the external radius of the domain reaches $r_e=4$.} 
\label{fig.ICF-MHD}
\end{center}
\end{figure}

\begin{figure}[!htbp]
\begin{center}
\begin{tabular}{c} 
\includegraphics[width=0.75\textwidth]{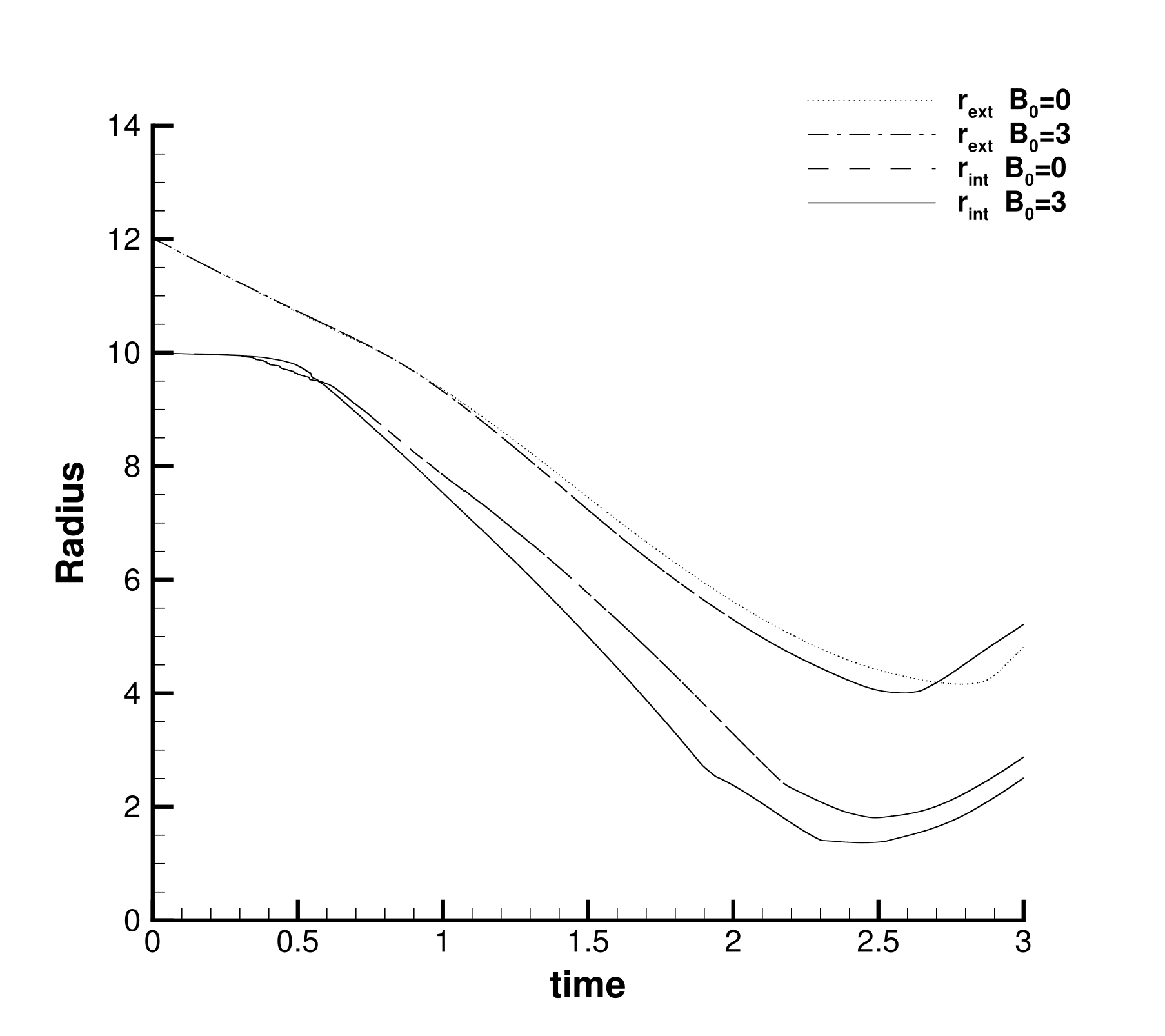}            
\end{tabular} 
\caption{Time evolution of the external radius and the interface position between light and heavy fluid for $B_0=0$ and $B_0=3$ up to the final time $t_f=3.0$. } 
\label{fig.ICF-radii}
\end{center}
\end{figure}

As evident from Figure \ref{fig.ICF-B0}, Rayleigh-Taylor phenomena arise along the interface between light and heavy fluid, generating vortex-like patterns as well as mesh rolling up that destabilize the fluid flow. In order to limit and reduce such instabilities, we apply a magnetic field $\mathbf{B}$ acting on the horizontal plane $x-y$ and we solve the ideal equations for magnetohydrodynamics (MHD) in our direct ALE framework. We refer the reader to \cite{LagrangeMHD,Lagrange3D} for more details on the MHD system and its implementation in the ADER context with moving meshes. In the following, we solve the spherical implosion problem applying a magnetic field of intensity $B_0$ to the fluid, that is
\begin{equation}
\mathbf{B} = (B_x,B_y,B_z) = \boldsymbol{\omega} \times \mathbf{x}, \qquad \boldsymbol{\omega}=(0,0,B_0).
\label{eqn.B_icf}
\end{equation}   
We set $B_0=1$, $B_0=2$, $B_0=3$ and Figure \ref{fig.ICF-MHD} plots the corresponding density distributions obtained at final times $t_{f,1}=2.70$, $t_{f,2}=2.63$ and $t_{f,3}=2.55$, respectively. For comparison purposes we also report the result computed with $B_0=0$ and one can note that the higher is the intensity of the magnetic field the smaller are the Rayleigh-Taylor instabilities, as expected. Moreover, in Figure \ref{fig.ICF-radii} the evolution of the external and the interface radii are shown for the case with $B_0=0$ and $B_0=3$, up to a final time of $t_f=3.0$.  

\begin{figure}[!htbp]
\begin{center}
\begin{tabular}{ccc} 
\includegraphics[width=0.33\textwidth]{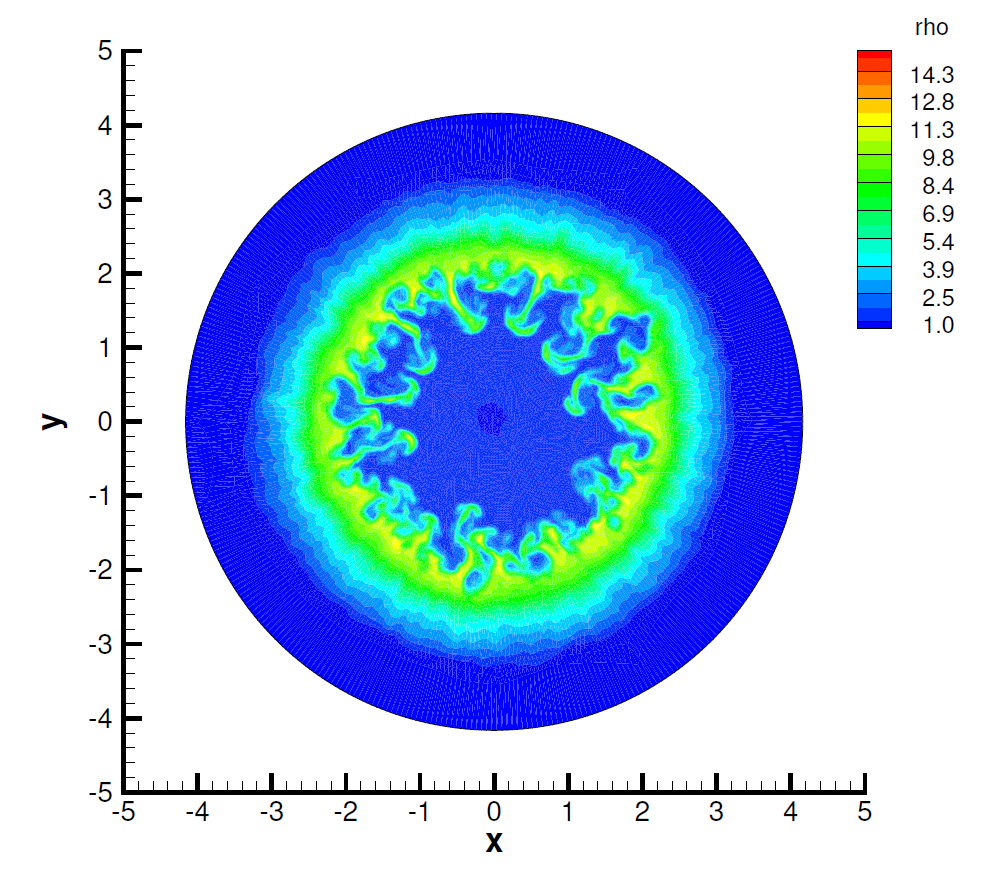} &           
\includegraphics[width=0.33\textwidth]{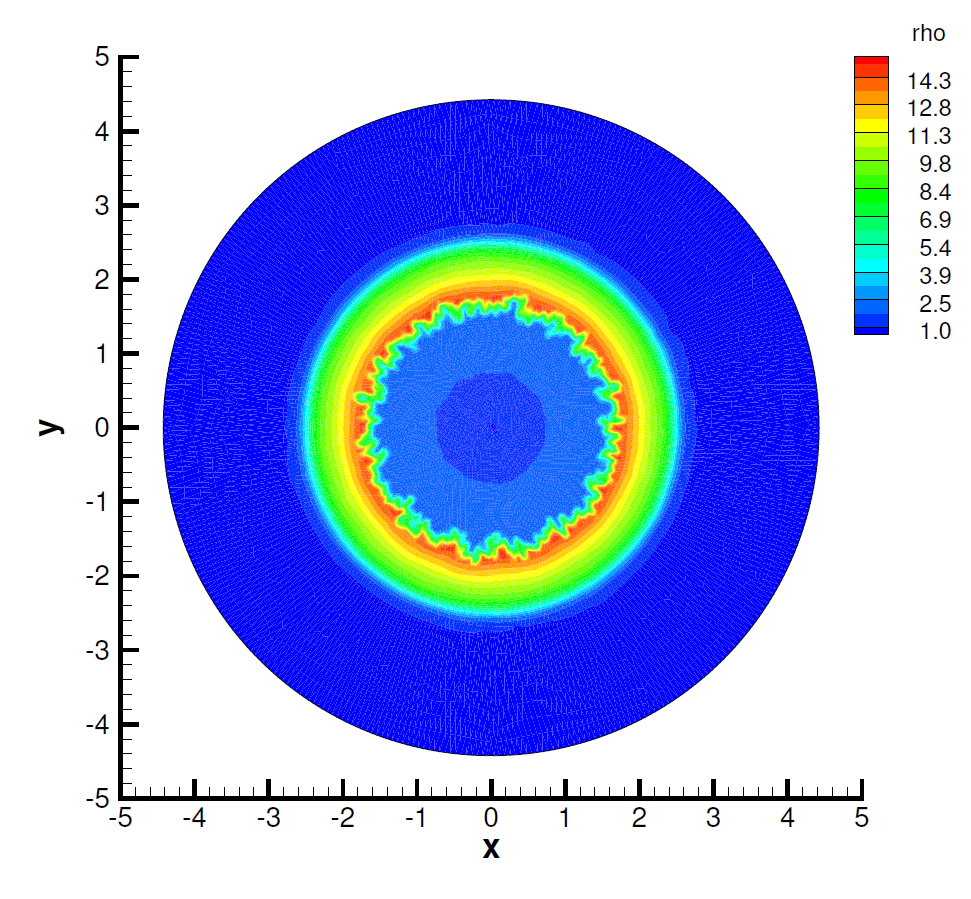} &
\includegraphics[width=0.33\textwidth]{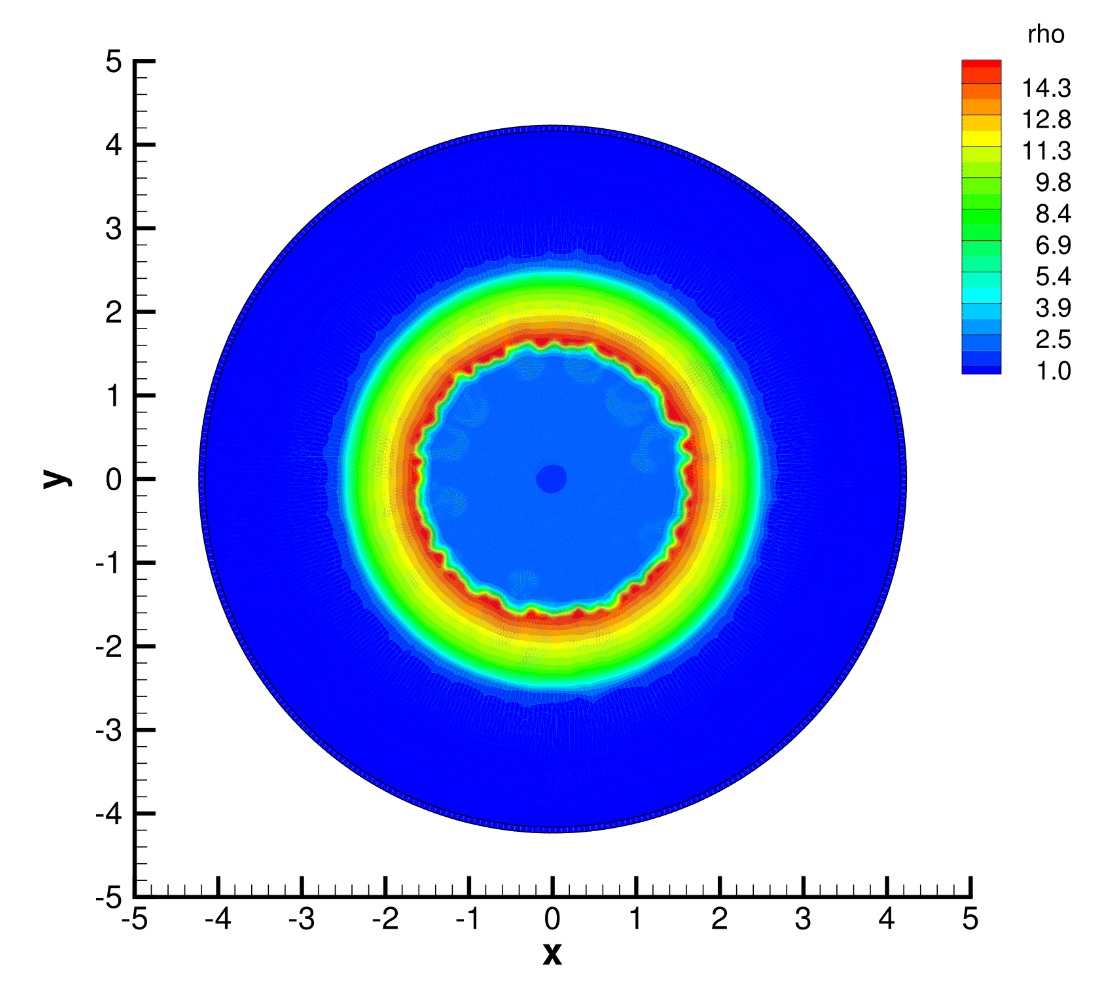} \\
\includegraphics[width=0.33\textwidth]{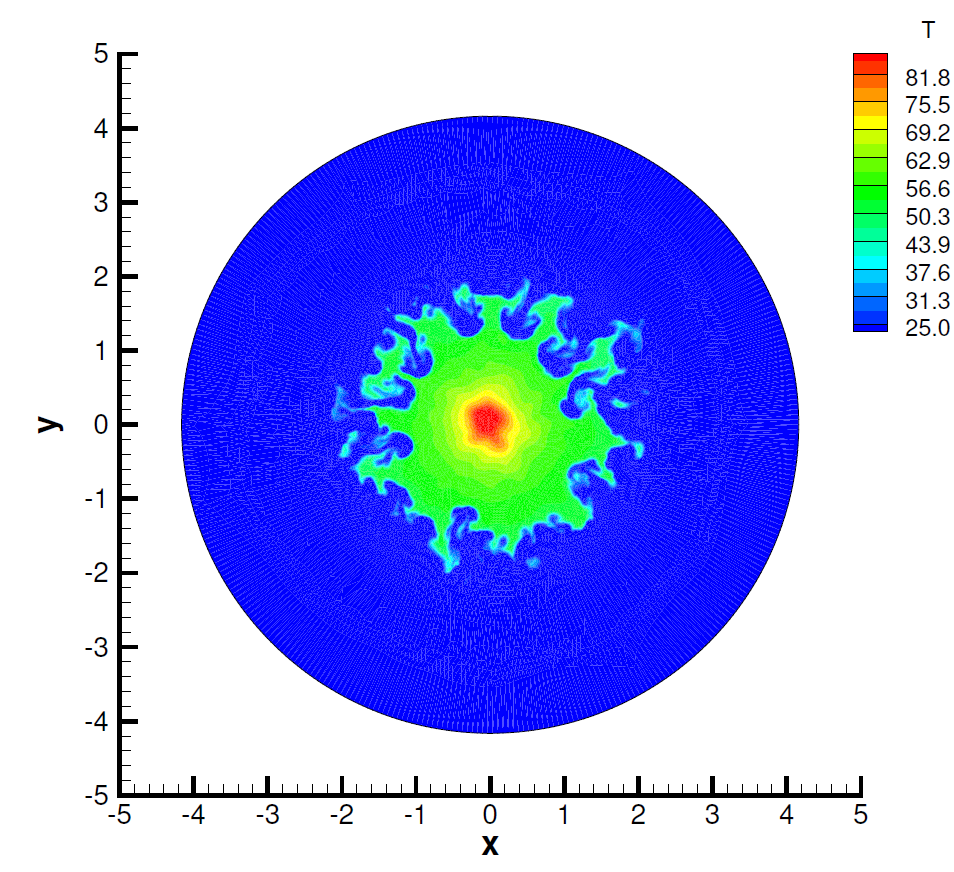} &           
\includegraphics[width=0.33\textwidth]{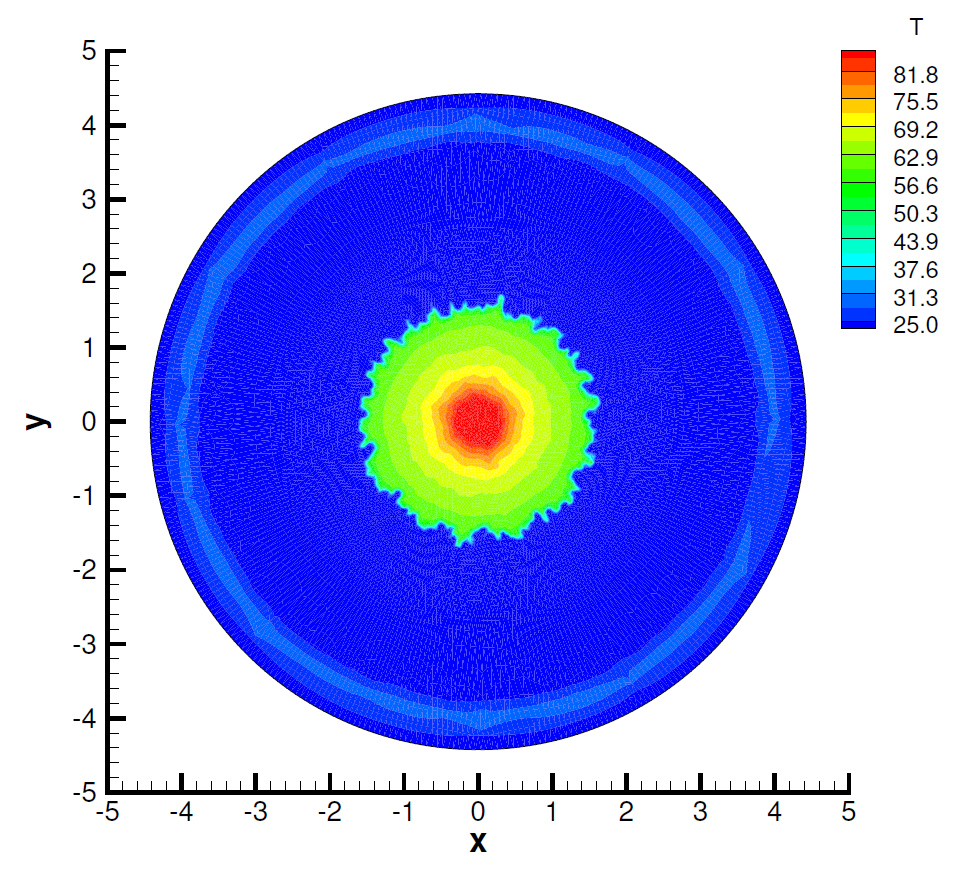} &
\includegraphics[width=0.33\textwidth]{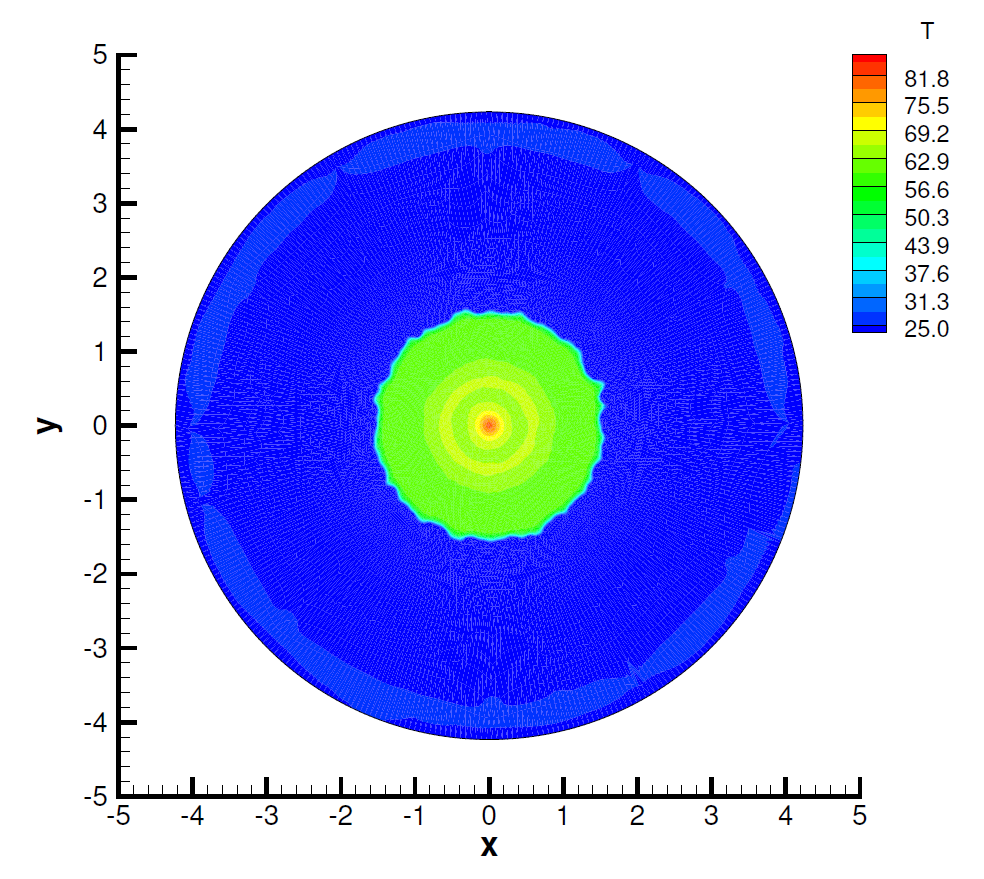} \\
\end{tabular} 
\caption{Density (top row) and temperature (bottom row) distribution at final time $t_f=2.70$ for the spherical implosion problem solved using the Euler equations of compressible gas dynamics (left column), the compressible Navier-Stokes equations with viscosity $\mu=10^{-3}$ (middle column) and the viscous relativistic MHD equations with a magnetic field of intensity $B_0=1$ applied on the horizontal plane $x-y$ and a viscosity coefficient of $\mu=10^{-3}$ (right column).} 
\label{fig.ICF-visc}
\end{center}
\end{figure}

Finally, we include \textit{physical viscosity} in the governing equations, hence we run again a fourth order simulation of this test case with a viscosity coefficient of $\mu=10^{-3}$ until the final time $t_f=2.7$. The results are depicted in Figure \ref{fig.ICF-visc}, where density as well as temperature are shown. We have solved the inviscid Euler equations for compressible gas dynamics \eqref{eulerTerms} (left panels), the compressible Navier-Stokes equations \eqref{eqn.NSflux} (middle panels) and the \textit{viscous resistive} MHD equations \cite{VRMHD} with a Lundquist number of $Lu=10^{3}$ (right panels). One can note that the physical viscosity plays an important role for the stabilization of the fluid.

\section{Conclusions}
\label{sec.concl}

In this paper we have presented a new family of high order ADER Discontinuous Galerkin (DG) finite element schemes in the framework of direct Arbitrary-Lagrangian-Eulerian (ALE) methods on moving unstructured multidimensional meshes. The numerical solution is represented by high order spatial polynomials of degree $N$ in each cell that are evolved in time by a one-step explicit DG scheme, based on a high order space-time predictor computed relying on the ADER methodology. \\
Two different strategies have been developed for moving the mesh in time, namely a piecewise linear decomposition of the control volumes into simplex sub-cells and a curved high order isoparametric approximation of the element geometry. For the sub-nodes lying on an element face, a new nodal solver based on the HLL state is used to evaluate the mesh velocity. The new geometry configuration is directly taken into account in the computation of the fluxes. The proposed explicit one-step ALE ADER-DG scheme is based on a space-time conservation formulation of the governing PDE system, hence satisfying by construction the geometrical conservation law (GCL).

Convergence studies demonstrate the space-time accuracy of the new schemes and a wide range of test cases have been run in order to assess the validity and the robustness of the ALE ADER-DG method. The Euler equations of compressible gas dynamics as well as the compressible Navier-Stokes equations with heat conduction have been considered, solving a set of test problems with strong shock waves and other discontinuities. Finally, a cylindrical implosion problem has been studied and a magnetic field has been applied to the fluid in order to stabilize the Rayleigh-Taylor instabilities arising in this problem. For this purpose the ideal classical and viscous resistive magnetohydrodynamics equations have been used within the framework of the new algorithm illustrated in this paper.

We plan to extend the presented approach to non-conservative systems and stiff source terms in order to apply it to the Godunov-Peshkov-Romenski model of nonlinear hyperelasticity \cite{PeshRom2014,Dumbser_HPR_16,LagrangeHPR}. 

\section*{Acknowledgments}
The presented research has been financed by the European Research Council (ERC) under the European Union's Seventh Framework 
Programme (FP7/2007-2013) with the research project \textit{STiMulUs}, ERC Grant agreement no. 278267. The authors acknowledge 
PRACE for awarding us access to the SuperMUC supercomputer of the Leibniz Rechenzentrum (LRZ) in Munich, Germany.

 \clearpage

\bibliography{biblio}
\bibliographystyle{plain}

\end{document}